\documentclass[11pt]{article}

\usepackage{environ}
\newcommand{\acksection}{\section*{Acknowledgments}}
\NewEnviron{ack}{%
	\acksection
	\BODY } \usepackage[verbose=true,letterpaper,margin=1in]{geometry}

\title{{\bf{Finite-Sample Guarantees for Learning
			Dynamics in Zero-Sum Polymatrix Games}}}
		
		\author{ Fathima Zarin Faizal$\textsuperscript{*}$ \quad Asuman
			Ozdaglar\textsuperscript{\textdagger} \quad Martin
			J. Wainwright\textsuperscript{\textdaggerdbl}\\ Department of
			Electrical Engineering and Computer
			Sciences\textsuperscript{*,\textdagger,\textdaggerdbl}\\ Department
			of Mathematics\textsuperscript{\textdaggerdbl} \\ Massachusetts
			Institute of
			Technology\textsuperscript{*,\textdagger,\textdaggerdbl}
			\\ \texttt{\{fathima,asuman,mjwain\}@mit.edu} \\ }

\usepackage{url} 
\usepackage{booktabs} 
\usepackage{nicefrac} 
\usepackage{microtype} 
\usepackage{lipsum}
\usepackage{amsfonts}
\usepackage{graphicx}
\usepackage{xcolor}
\usepackage{epstopdf}
\usepackage{amsmath,amssymb,amsthm}
\usepackage{algorithm,algpseudocode,amssymb,comment,enumitem,xspace,graphicx,subcaption}
\usepackage{wrapfig} 
\usepackage{hyperref,cleveref}
\usepackage{xurl}
\hypersetup{breaklinks=true}

\usepackage[mathscr]{eucal}

\newcommand{\strat}[1][]{\ensuremath{{\pi^{#1}}}} 
\newcommand{\BStrat}{\boldsymbol{\pi}}  
\newcommand{\simplex}[1]{\ensuremath{\Pi}^{#1}} 
\newcommand{\BigSimplex}{\ensuremath{\boldsymbol{\Pi}}} 
\newcommand{\graph}{\ensuremath{\mathcal{G}}} 
\newcommand{\edges}{\ensuremath{\mathcal{E}}} 
\newcommand{\nodes}{\ensuremath{\mathcal{N}}} 
\newcommand{\dmax}{\ensuremath{d_{\max}(\graph)}}
\newcommand{\bq}{\boldsymbol{q}} 
\newcommand{\NumPlayers}{\ensuremath{N}} 
\newcommand{\actions}{\ensuremath{\mathcal{A}}} 
\newcommand{\BigActions}{\ensuremath{\boldsymbol{\mathcal{A}}}} 
\newcommand{\RMatrix}[2]{\ensuremath{R}^{(#1,#2)}} 
\newcommand{\BRMatrix}{\boldsymbol{R}} 
\newcommand{\cir}{\ensuremath{r}}
\newcommand{\NEstrat}{\overline{\strat}}
\newcommand{\NEBStrat}{\overline{\BStrat}}
\newcommand{\AvgPayoff}[1]{\ensuremath{q^{#1}}}
\newcommand{\BAvgPayoff}{\ensuremath{\boldsymbol{q}}}
\newcommand{\action}{\ensuremath{a}}
\newcommand{\BAction}{\ensuremath{\boldsymbol{a}}}
\newcommand{\regBR}[1][i]{\ensuremath{{\sigma_{\tau}^{#1}}}}

\newcommand{\pii}{\strat[i]} \newcommand{\pmi}{\pi^{-i}}
\newcommand{\nashg}{\operatorname{NG}}

\newcommand{\bigdelta}{\boldsymbol{\Delta}}
\newcommand{\vkl}{\lyap{V}_{\scriptstyle {\mbox{KL}}}}
\newcommand{\valt}{\lyap{V}_{\text{\mbox{alt}}}}
\newcommand{\myexp}{\ensuremath{\eta}}
\newcommand{\koff}{\ensuremath{k_0}}
\newcommand{\offset}{\ensuremath{\nu}}
\newcommand{\pihat}{\ensuremath{\widehat{\pi}}}

\newcommand{\real}{\ensuremath{\mathbb{R}}}
\newcommand{\infosetuplower}{minimal information\xspace}

\newcommand{\constYoung}{r}
\newcommand{\uplowratio}{\ensuremath{\xi}}
\newcommand{\timescalesep}{\ensuremath{c_{\alpha,\beta}}}
\newcommand{\nablaVBoundPartialFirst}{\ensuremath{S_1}}
\newcommand{\nablaVBoundPartialSecond}{\ensuremath{S_2}}
\newcommand{\Btderror}[1]{\boldsymbol{w}_{#1}}

\newcommand{\SmallEdgeFrac}{\rho}
\newcommand{\SmallEdgeScaling}{\omega}

\newcommand{\lyap}[1]{\ensuremath{\mathcal{#1}}}
\newcommand{\TotLyap}{\ensuremath{\lyap{T}}}
\newcommand{\TDLyap}{\ensuremath{\lyap{W}}}
\newcommand{\stratLyap}{\ensuremath{\lyap{V}}}
\newcommand{\StratLyap}{\stratLyap}

\newcommand{\StratLyapVal}[1]{\ensuremath{V_{#1}}}
\newcommand{\stratLyapVal}[1]{\ensuremath{\StratLyapVal{#1}}}
\newcommand{\TotLyapVal}[1]{\ensuremath{T_{#1}}}

\newcommand{\TDLyapVal}[1]{\ensuremath{W_{#1}}}

\newcommand{\R}{\mathbb{R}} 
 \newcommand{\N}{\mathbb{N}}
 
\newcommand{\A}{\mathcal{A}}

 \newcommand{\E}{\mathbb{E}}

\newcommand\numberthis{\addtocounter{equation}{1}\tag{\theequation}}

\definecolor{OliveGreen}{cmyk}{0.91,0.0, 0.99,0.1} 

\definecolor{DarkPurple}{rgb}{0.7,0,0.7}

\usepackage{scalerel}

\newcommand{\smallsub}[2]{#1_{\scaleto{#2}{4pt}}}
\newcommand{\tinysub}[2]{#1_{\scaleto{#2}{3pt}}}
\newcommand{\medsub}[2]{#1_{\scaleto{#2}{6pt}}}

\newcommand{\Vinit}{V_1}
\newcommand{\HarrisLyap}{\ensuremath{\smallsub{\lyap{V}}{\textrm{H}}}}

\newcounter{lcc}
\newcommand{\newlcc}[1]{\ifthenelse{\equal{#1}{}}{\stepcounter{lcc}\ensuremath{c_{\arabic{lcc}}}}{\refstepcounter{lcc}\label{#1}\ensuremath{c_{\thelcc}}}}

\newcommand{\matsnorm}[2]{|\!|\!| #1 | \! | \!|_{{#2}}}
\newcommand{\vecnorm}[2]{\| #1\|_{#2}}

\newcommand{\opnorm}[1]{\ensuremath{\matsnorm{#1}{{2}}}}

\newcommand{\inprod}[2]{\ensuremath{\langle #1 , \, #2 \rangle}}

\newcommand{\Kgood}{\ensuremath{\medsub{K}{\mbox{good}}}}

\newcommand{\myunderparagraph}[1]{\vspace*{0.05in} \noindent \underline{{#1}}}

\newcommand{\Amax}{\ensuremath{A_{\max}}}

\newcommand{\lammax}{\ensuremath{\tinysub{\lambda}{\mbox{max}}}}

\newcommand{\Valt}{\valt}

\newcommand{\Gfun}{\ensuremath{\mathcal{G}}}
\newcommand{\EntFun}{\ensuremath{H}}

\newcommand{\myparagraph}[1]{\noindent {\bf{#1.}}}

\newcommand{\Amat}{\ensuremath{\mathbf{A}}}
\newcommand{\Mmat}{\ensuremath{\mathbf{M}}}

\newcommand{\PayOff}{\ensuremath{\bar{R}}}

\usepackage[giveninits=true,
backref=true,
backend=biber,
style=alphabetic,
sortlocale=de_DE,
natbib=true,
doi=false,
isbn=false,
url=false]{biblatex}

\newtheorem{theorem}{Theorem} \newtheorem{prop}{Proposition}
 \newtheorem{corollary}{Corollary}
\newtheorem{lemma}{Lemma} 
 
\newcommand{\defn}{\ensuremath{: \, =}}

\begin{document}

\maketitle

\begin{abstract}
We study best-response type learning dynamics for zero-sum polymatrix
games under two information settings. The two settings are
distinguished by the type of information that each player has about
the game and their opponents' strategy. The first setting is the full
information case, in which each player knows their own and their
opponents' payoff matrices and observes everyone’s mixed
strategies. The second setting is the \infosetuplower case, where
players do not observe their opponents' strategies and are not aware
of any payoff matrices (instead they only observe their realized
payoffs). For this setting, also known as the radically uncoupled case
in the learning in games literature, we study a two-timescale learning
dynamics that combine smoothed best-response type updates for strategy
estimates with a TD-learning update to estimate a local payoff
function. For these dynamics, without additional exploration, we
provide polynomial-time finite-sample guarantees for convergence to an
$\epsilon$-Nash equilibrium.
\end{abstract}


\section{Introduction}

Game theory is used to model interactions between two or more players,
with the assumption that each player behaves rationally so as to
maximize their individual payoffs subject to available information.  A
standard solution concept is that of a Nash equilibrium: it can be
justified as the steady-state outcome of the learning process of
players acting in their self-interest. This motivates the study of
various learning dynamics by which players choose an action while
learning the strategy of their opponents and adjust their play
accordingly. Key questions associated with a given set of learning
dynamics are stability, convergence to equilibria, and in a more
refined analysis, the rate of convergence to an equilibrium.

Among the most natural and well-studied forms of learning dynamics is
fictitious play (FP), first introduced by
Brown~\cite{Brown49,brownfp1951}. It is a simple and interpretable
form of dynamics: players use information from previous rounds of play
to estimate the opponent's strategy, and then play a best-response
action based on this estimate. For two-player zero-sum games,
Robinson~\cite{robinson1951} established the asymptotic convergence of
FP. This was followed by a number of papers that studied its
convergence properties for different classes of games;
see~\Cref{SecRelated} below for further details.  Classical FP assumes
that each player observes their opponent's actions and knows the
payoff function, which can be viewed as a full information
setting. While this is a useful benchmark, players may have limited
information about their opponents' play or even their own payoff
functions. In this paper, in addition to studying the full information
setting, we also study the case where players do not observe their
opponents' play; they only have access to their realized payoffs at
each step, a case which we refer to as the \infosetuplower setting.

We examine these issues within the class of zero-sum polymatrix games,
which is a generalization of the class of two-player zero-sum games to
the multiplayer setting. Players interact with each other in a
pairwise manner, and an underlying interaction graph captures each
player's payoff's dependencies on the other players. A pairwise matrix
game is played on each edge and players choose a single strategy for
each of the pairwise games they play in their neighborhood. Each
player's payoff is the sum of the payoffs they receive from each of
the pairwise games they play. The zero-sum constraint ensures that the
sum of the payoffs of all the players equals zero, i.e., there is no
total flux of payoffs in or out of the system. Note that the matrix
game played on each edge \emph{need not} be zero-sum; each edge-based
game being zero-sum is a special case of the more general class of
zero-sum polymatrix games.

Zero-sum polymatrix games are used to model scenarios where pairwise
interactions are dominant. For instance, consider a competitive
resource allocation problem where multiple countries attempt to
allocate their defense budget in order to obtain maximum control of
shared resources. The existence of an edge between two countries in an
underlying interaction graph $(\nodes,\edges)$ indicates whether or
not those countries share a border. The action taken by each country
in the polymatrix game would be to divide their defense budget among
each of the pairwise interactions in which they are involved. On a
particular edge $(i,j) \in \edges$ between two countries $i$ and $j$,
the country that allocates more budget for edge $(i,j)$ receives a
payoff of $r_{i,j}$ which represents the utility obtained from
controlling the shared resource on that edge; the country that
allocates less budget for edge $(i,j)$ receives a payoff of 0. While
not a zero-sum game, this is a constant-sum polymatrix game, i.e., the
total sum of payoffs across all countries for each possible budget
allocation that all countries choose is equal to $\sum_{(i,j) \in
  \edges} r_{i,j}$. By subtracting this quantity from the payoff
matrices of each player, this can be transformed into a zero-sum
polymatrix game.

Despite much work on (asymptotic) convergence guarantees for
normal-form games, there are relatively few results on the iteration
complexity of best-response dynamics.  The iteration complexity, for a
given tolerance level $\epsilon>0$, refers to the number of rounds
$K(\epsilon)$ required to obtain strategies that form an
$\epsilon$-optimal Nash equilibrium. Of particular relevance is the
recent work of Chen et al.~\cite{zaiwei2023}, which considered a
smoothed best-response type dynamics for two-player zero-sum
normal-form and Markov games.  They established explicit bounds on the
iteration complexity, but in the absence of an additional exploration
device---such as mixing with the uniform distribution---their bounds
on $K(\epsilon)$ scale exponentially in $(1/\epsilon)$.


\subsection{Our contributions}

In this paper, we study smoothed best-response type dynamics for
zero-sum polymatrix games.  We analyze their behavior in both a
\emph{full information} setting, as well as in the
\emph{\infosetuplower} setting.  In the latter setting, also referred
to as the radically uncoupled setting or the bandit setting, each
player observes only their own realized payoff at each round, and has
no other information about either the game or their opponents' play.
Focusing on smoothed best-response dynamics without any modifications
to encourage additional exploration, we prove that the number of
iterations $K(\epsilon)$ required to converge to $\epsilon$-optimal
Nash equilibrium scales polynomially in the ratio $1/\epsilon$ for
both information settings.  To the best of our knowledge, this is the
first known polynomial-time guarantee for best-response type dynamics
in the \infosetuplower setting without the introduction of additional
exploration in the learning dynamics.

To be clear, our focus in this paper is on the properties of natural
best-response type dynamics, and \emph{not} on developing
alternative---and possibly more efficient---algorithms for computing
Nash equilibria.  By searching more broadly in the space of updates,
it can be possible to obtain faster iteration complexities (e.g., see
the papers~\cite{AoCenChi2022asynchronous,cai2016polymatrix} for
results of this type). However, best response dynamics are a classical
and arguably natural form of interaction between a collection of
agents.  By providing convergence guarantees for best response, we
give evidence for the emergence of Nash equilibria as efficiently
achievable steady-state behavior of a collection of agents who each
act greedily at each step while playing a zero-sum polymatrix game.

In more detail, our first main result applies to the simpler case of
full information, where players can perform smoothed best-response
updates based on knowledge of the opponents' mixed strategy. We
introduce a modified version of Lyapunov functions used in previous
work~\cite{harris1998,zaiwei2023}, and use it to prove that the
dynamics converge to an $\epsilon$-Nash equilibrium in $K(\epsilon)
\asymp (1/\epsilon^2)$ iterations.  In contrast to previous
analysis~\cite{zaiwei2023}---which led to exponential dependence in
$(1/\epsilon)$---our modification has desirable smoothness properties
that allow us to establish the claimed polynomial scaling in
$1/\epsilon$.

We then turn to the \infosetuplower setting, in which players do not
observe their opponents' play, but instead only observe their realized
payoffs at each round.  For this more challenging setting, we analyze
a simple best-response type learning dynamics involving updates on two
timescales.  First, on the faster timescale, each player maintains and
updates an estimate of their local payoff function or
$q$-value---i.e., the average payoff as a function of their
actions. The local payoff function carries information about the
opponents' strategy. At the same time, on a slower timescale, the
players also modify their strategies by forming a smoothed
best-response based on the estimated $q$-values. The updates of the
$q$-values are in the spirit of TD learning~\cite{sutton1988learning},
and involve an adaptive learning rate~\cite{lesliecollins2005} that
ensures unbiased estimates of the local payoff function. When recast
as a form of stochastic approximation \cite{Borkar2008}, this unbiased
property means that the estimated $q$-values are updated using
zero-mean noisy estimates of the underlying payoff function.

However, a major challenge is that the noise variance explodes as the
player's strategies approach the boundary of the probability
simplex. We overcome this variance explosion---without any
modifications to the updates---by explicitly tracking how quickly the
strategies approach the boundary. Combined with careful design of
Lyapunov functions with favorable smoothness properties, we prove an
upper bound on the iteration complexity that scales as
$(1/\epsilon)^{8 + \offset}$, where the offset $\offset > 0$ can be
chosen arbitrarily close to zero at the price of growth in constant
pre-factors.

Finally, for both information settings, our bounds on the iteration
complexity show a dependence on the underlying graph and the pairwise
games. In the worst case, it scales with the maximum degree of the
graph but depending on the graph structure and the nature of the
edge-wise games, it can exhibit a much milder dependence.


\subsection{Related work}
\label{SecRelated}

So as to put our contributions in context, we now turn to a discussion
of past work on best-response type dynamics and
the \infosetuplower setting.  To be clear, this overview is far from
comprehensive: we limit our discussion to the literature most closely
related to our work.

There is a long line of work on the convergence of fictitious play
(FP) for various classes of finite normal-form games; all of the
classical work assumes knowledge of the payoff matrices.  In more
detail, following the introduction of fictitious
play~\cite{Brown49,brownfp1951}, Robinson~\cite{robinson1951}
established the convergence of both the payoffs and strategies of
discrete-time FP applied to two-player zero-sum games.  Later work by
Shapiro~\cite{shapiro1958note} proved the payoffs converge at least as
quickly as $k^{-1/(|\A^1| + |\A^2| -2)}$, where $\A^1, \A^2$ are the
action sets of each player.  Karlin~\cite{karlin1959mathematical}
conjectured that the actual convergence rate was $k^{-1/2}$, but it
was later shown that this conjecture does not hold for general
tie-breaking rules~\cite{daskalakis2014counter,abernethy2021fast}.
Miyasawa~\cite{miyasawa1961convergence} established FP convergence for
two-person non-zero-sum games in which players have at most two
actions and satisfy a notion of non-degeneracy; other classes known to
have convergence properties under FP dynamics include common interest
games~\cite{monderer1996fictitious} and weighted-potential
games~\cite{monderer1996potential}.  On the negative side,
Shapley~\cite{Shapley1963} provided a two-person non-zero-sum game
with three actions for which FP fails to converge; moreover, other
counterexamples have been given for certain coordination
games~\cite{foster1998nonconvergence}.

Harris~\cite{harris1998} studied the continuous-time version of
fictitious play, and used a Lyapunov function argument to show that it
converges at the rate $(1/t)$ for two-player zero-sum games.  Our
analysis makes use of a modification of the Lyapunov function from
this paper.  Smoothed versions of fictitious play---in which actual
best-responses are replaced by smoothed versions---have been studied
in various papers~\cite{fudenbergkreps1993, hofbauersandholm2002}; we
also analyze a smoothed version in this paper. In our work, we focus
on best-response type dynamics based on fictitious play since they are
a well-studied model of learning for myopic agents. We note that
recent literature has studied various other dynamics to compute
equilibria of different classes of games, including gradient descent,
mirror descent and its variants
(e.g.,~\cite{zeng2022regularized,sokota2022unified}); regret
matching~\cite{hart2000simple}; extragradient versions of
multiplicative updates~\cite{cen2021fast}; as well as various
regret-based and online learning methods
(e.g.,~\cite{Cesa-Bianchi_Lugosi_2006,daskalakis2011,rakhlin2013optimization,syrgkanis2015fast,chen2020hedging}).

There is also a more recent and evolving line of work on the
\infosetuplower setting for the class of two-player zero-sum
games. Leslie and Collins~\cite{lesliecollins2005} introduced the use
of adaptive stepsizes to estimate the average payoffs of each player,
and established its asymptotic convergence.  Our algorithm includes
updates of this type, and the adaptive stepsizes ensure a key
unbiasedness property.  In addition, inspired by the
papers~\cite{leslie2003convergent,zaiwei2023}, we also make use of a
\textit{doubly-smoothed} best-response to update the strategies where
in addition to a smoothed best-response, we use a learning rate to mix
with the current strategy estimate. Our set-up has connections with
the entropy-regularized algorithms studied in the
paper~\cite{bakhtin2022mastering}. Their analysis is predicated upon
availability of the exact payoff functions, as well as the opponent's
actions, whereas by contrast, we assume only availability of the
random payoff resulting from the unobserved actions. The
paper~\cite{zaiwei2023} also analyzes this form of minimal
information; while their main focus is on stochastic games, they also
establish a bound on iteration complexity for two player zero-sum
matrix games.  However, they do not make use of adaptive
stepsizes~\cite{lesliecollins2005}, and provide guarantees relative to
$\epsilon$-Nash equilibrium up to a smoothing bias. Incorporating the
smoothing bias means that the end-to-end guarantees grow exponentially
in the inverse tolerance $(1/\epsilon)$, whereas one of our main
contributions is to provide schemes (and analysis) with polynomial
scaling.

Some recent work~\cite{yangcai2023,ouhamma2023learning} has also
provided guarantees with polynomial dependence on $(1/\epsilon)$ for
the minimal information setting, but for schemes that depart from the
usual best-response learning dynamics. In the
paper~\cite{yangcai2023}, exploration is enforced by explicitly
limiting how quickly strategies are allowed to approach the boundary
of the probability simplex; the work~\cite{zaiwei2023} also notes that
mixing with a uniform distribution can alleviate these issues. On the
other hand, the analysis in the paper~\cite{ouhamma2023learning}
regularizes using the Tsallis entropy, resulting in a non-standard
strategy update.

It is also possible to exploit linear programming methods so as to
compute the equilibria of zero-sum polymatrix games (cf. the
papers~\cite{Janovskaja1968,bregman1987methods,bregman1998separable,cai2016polymatrix}).
Leonardos et al.~\cite{leonardos2021exploration} studies the
continuous-time version of our full information dynamics.  In other
related work in the full information setting, Ao et
al.~\cite{AoCenChi2022asynchronous} provide finite-time guarantees for
zero-sum polymatrix games, in particular by analyzing the finite-time
convergence of the optimistic multiplicative weights (OMWU) method.
They show that the rate of convergence of their algorithm scales with
the maximum degree of the graph.


\section{Background and problem set-up}

We now provide background along with a more precise set-up of the
problem. We begin in~\Cref{SecZS} with notation and the basic
formalism of zero-sum polymatrix games. In~\Cref{SecDynamics}, we
discuss various types of learning dynamics, including the two sets of
updates (for full and \infosetuplower respectively) that we study in
this paper.

\subsection{Zero-sum polymatrix games}
\label{SecZS}

A zero-sum polymatrix game with $\NumPlayers$ players is defined by an
underlying graph $\graph=(\nodes,\edges)$ with $\nodes = \{1, \ldots,
\NumPlayers\}$. Each node $i \in \nodes$ corresponds to a player who
has access to a finite set $\actions^i$. Let $\BigActions =
\otimes_{i=1}^n \actions^i$ represent the joint action set and let
$\Amax=\underset{i \in \nodes}{\max} \ {|\actions^i|}$ be the maximum
number of actions across all the players. Each edge $(i,j)$ of the
graph is associated with a pairwise game between players $i$ and $j$,
with the matrix $\RMatrix{i}{j} \in \R^{\vert \actions_i \vert \times
  \vert \actions_j \vert}$ representing the payoff matrix of player
$i$ in the matrix game played by players $i$ and $j$. If there is no
edge between these two players, i.e., $(i,j) \notin \edges$, then
$\RMatrix{i}{j}$ is the zero matrix of appropriate dimensions.  We
assume throughout that the payoff matrices are normalized so that
$\max_{\action^i, \action^j} | \RMatrix{i}{j}(\action^i,\action^j) |
\leq 1 $. We compile the pairwise payoff matrices into a block matrix
$\BRMatrix$ whose $(i,j)^{\text{th}}$ block is given by
$\RMatrix{i}{j}$.

At each discrete time instant, all players simultaneously choose a
strategy to apply across all of their pairwise interactions. Players
may randomize their choice of actions, so for each player $i \in
\nodes$, we denote their set of possible strategies by $\simplex{i}$,
corresponding to the probability simplex supported on $\actions^i$. We
use $\BigSimplex \defn \otimes_{i=1}^n \simplex{i}$ to denote the set
of all possible joint strategies that can be played at each time
instant. For any $i \in \nodes$, we use $\strat^{-i}$ to denote the
collection of all strategies indexed by players $j \in \nodes
\backslash \{i \}$; it belongs to the Cartesian product of simplices
denoted by $\simplex{-i}$.  For any strategy profile $\strat^{-i} \in
\simplex{-i}$, we denote the vector of expected payoffs of player $i$
for each action they play by $\AvgPayoff{i}(\strat^{-i}) \in
\R^{|\actions^i|}$, and note the relationship
\begin{align*}
  \AvgPayoff{i}(\strat^{-i} ) = \sum_{j \in \nodes \setminus \{i\}}
  \RMatrix{i}{j} \strat[j] = \sum_{j \in \nodes} \RMatrix{i}{j}
  \strat[j],
\end{align*}
i.e., the sum of the expected payoffs from the matrix games that $i$
plays on each of the edges connected to $i$. At times, we use
$\BAvgPayoff(\BStrat) = (\AvgPayoff{i}(\strat^{-i}))_{i \in \nodes}$
to refer to the collection of average payoffs of each player under
$\BStrat$. The zero-sum constraint on the polymatrix game means that
the sum of the payoffs across all players equals zero---that is
\begin{align}
\label{EqnZeroSum}
\sum_{i \in \nodes} {(\strat[i])}^{\top} \AvgPayoff{i}(\strat^{-i}) =
\BStrat^\top \BRMatrix \BStrat = 0 \qquad \mbox{for any $\BStrat \in
  \BigSimplex$.}
\end{align}
A joint strategy $\NEBStrat \in \BigSimplex$ is said to be a
\textit{Nash equilibrium} if
\begin{align}
\label{EqnNashPoly}
      {(\NEstrat^i)}^{\top} {\AvgPayoff{i}({\NEstrat}^{-i} )}\geq
      {(\strat^{i})}^{\top} {\AvgPayoff{i}({\NEstrat}^{-i} )}, \quad
      \mbox{for all players $i \in \nodes$, and $\strat^{i} \in
        \simplex{i}$.}
\end{align}
Nash's existence theorem \cite{nash1950equilibrium} ensures the
existence of a Nash equilibrium in any finite $\NumPlayers$-player
game. When $\NumPlayers=2$, the class of zero-sum polymatrix games
reduces to the class of two-player zero-sum games, in which case the
existence of a Nash equilibrium also follows from Von Neumann's
minimax theorem~\cite{von_neumann_minmax}. If for some $\epsilon > 0$
the inequalities~\eqref{EqnNashPoly} hold up to an $\epsilon$-additive
relaxation for some $\NEBStrat_{\epsilon}$, then
$\NEBStrat_{\epsilon}$ is called an \textit{$\epsilon$-Nash
  equilibrium}.


\paragraph{Nash gap}
Let us now introduce a measure of the closeness to a Nash equilibrium
(or NE for short).  For any strategy $\BStrat \in \bigdelta$, the
\emph{Nash gap} is defined as
\begin{align}
\label{EqnNashGap}
\nashg(\BStrat) = \sum_{i\in \nodes} \left\{ \underset{\hat{\strat}
  \in \simplex{i}}{\max} \left( \hat{\strat} - \strat[i] \right)^{\top}
\AvgPayoff{i}(\strat^{-i}) \right\}.
\end{align}
\noindent Observe that any mixed strategy $\BStrat \in \bigdelta$ with
Nash gap bounded as $\nashg(\BStrat) \leq \epsilon$ is guaranteed to
be an \emph{$\epsilon$-Nash equilibrium}, in the sense that it
satisfies the NE inequalities~\eqref{EqnNashGap} up to an additive
offset of $\epsilon$. Our goal is to show that natural best-response
type dynamics that have been studied in the literature can be used to
recover an $\epsilon$-Nash equilibrium in time polynomial in
${1}/{\epsilon}$ even when players cannot observe the strategies used
by the other players.

\subsection{Learning dynamics}
\label{SecDynamics}

We now describe and provide intuition for the different learning
dynamics we analyze in this paper. We begin with the smoothed
best-response and then discuss a form of best-response dynamics for
the full and minimal information settings (see~\ref{AlgFull}
and~\Cref{AlgPartial}).


\subsubsection{The smoothed best-response}
\label{SecSmoothedBR}

Recall the definition of $\AvgPayoff{i}(\strat^{-i})$ as the vector of
expected payoffs of player $i$ when all the players except for player
$i$ choose their strategy to be $\strat^{-i}$. The \emph{best-response
function} for player $i$ maps a strategy profile $\strat^{-i}$ of
players other than $i$ to a strategy for player $i$ via
\begin{subequations}
\begin{align}
\label{EqnBestResponse}
  \operatorname{br}^i(\strat^{-i}) & \in \arg \max_{\hat{\strat} \in
    \simplex{i}} \ {\hat{\strat}}^{\top}
  \AvgPayoff{i}(\strat^{-i}).
\end{align}
In this paper, we study a smoothed variant of the best-response that
is indexed by a \emph{regularization parameter} $\tau>0$. For a
strategy $\strat[i] \in \simplex{i}$, we define the \emph{Shannon
entropy} function
  \begin{align}
    \label{EqnDefnShannon}
    \EntFun(\strat[i]) \defn - \sum_{\action^i \in \actions^i}
    \strat^i(\action^i) \log \strat^i(\action^i).
  \end{align}
Now suppose that rather than the pure best-response, player $i$
instead plays the \textit{$\tau$-regularized best-response}
\begin{align}
\label{eqn:perturbed_best_response}
\regBR(\strat^{-i}) & \defn \arg \max_{\hat{\strat} \in \simplex{i}}
\big\{{\hat{\strat}}^{\top} \AvgPayoff{i}(\strat^{-i}) + \tau
\EntFun(\hat{\strat})\big\},
\end{align}
where $\tau > 0$ is a smoothing parameter.  In discrete choice theory,
for a given $\strat^{-i}$, $\regBR(\strat^{-i})$ is the strategy that
player $i$ would play if their payoffs had been perturbed by
Gumbel-distributed noise before choosing the action with the maximum
payoff. Such a perturbation to the payoffs can be used to model an
outsider's uncertainty about player $i$'s unobserved
preferences~\cite{anderson1992discrete}.  Also called the logit choice
model, the $\tau$-regularized best
response~\eqref{eqn:perturbed_best_response} has been studied
extensively in the learning-in-games literature, particularly in
connection with the fictitious play paradigm, where players
iteratively update their strategies based on estimated opponent
behavior~\cite{hofbauersandholm2002,fudenbergkreps1993}.

From an algorithmic perspective, there are several benefits associated
with using such a smoothed best-response as well.  The update rule is
unique: the optimization problem~\eqref{eqn:perturbed_best_response}
admits the closed-form expression
\begin{align}
\regBR(\strat^{-i})(\action^i) =
\frac{e^{\AvgPayoff{i}(\strat^{-i})(\action^i)/\tau}}{\sum_{a \in
    \actions^i} e^{\AvgPayoff{i}(\strat^{-i})(a)/\tau} } \qquad
\mbox{for each $\action^i \in \actions^i$.}
\end{align}
\end{subequations}
Consequently, for any $\tau > 0$, the strategy assigns strictly
positive mass to each action, thereby ensuring exploration.  Moreover,
as $\tau \rightarrow 0^+$, the smoothed best-response
$\regBR(\strat^{-i})$ converges to the best-response
$\operatorname{br}(\strat^{-i})$.

Using a smoothed best-response leads to a natural relaxation of a Nash
equilibrium.  We say that a mixed strategy $\BStrat $ is a
\textit{$\tau$-regularized Nash equilibrium} if
\begin{align}
\label{EqnTauRegNash}
\pii = \regBR(\strat^{-i}), \ i \in \nodes.
\end{align}
The existence of a $\tau$-regularized NE follows from the Brouwer
fixed-point theorem.\footnote{In fact, for each $\tau>0$, there is a
unique strategy satisfying condition~\eqref{EqnTauRegNash}; see
Proposition~\ref{PropLyapKL} in the sequel for details.}  Moreover, it
can be shown that a $\tau$-regularized equilibrium is an approximate
Nash equilibrium in a precise sense: given any tolerance $\epsilon >
0$, setting $\tau \leq \epsilon/(\NumPlayers \log \Amax)$ ensures that
any $\tau$-regularized Nash equilibrium has Nash
gap~\eqref{EqnNashGap} at most $\epsilon$, i.e., a $\tau$-regularized
Nash equilibrium is an $\epsilon$-Nash equilibrium.  This notion of a
$\tau$-regularized equilibrium plays a central role in our dynamics
and analysis. \\

We now turn to the two classes of smoothed best-response dynamics that
we study.

\subsubsection{Dynamics for full information}

We begin by analyzing smoothed best-response updates in the simpler
full information setting. In this case, players have access to the
payoff matrices on every edge, and observe the other players' mixed
strategies at all times.  The classical best-response dynamics in
continuous time, given by $\dot \strat^i_t \in
\operatorname{br}(\strat^{-i}_t) - \strat^i_t$ \mbox{for each $i \in
  \nodes$.}  The set of Nash equilibria of the underlying game
coincides with the equilibrium points of these dynamics. Since each
player has access to their opponents' previous strategies in this
information setting, player $i$ employs the best-response type update
\mbox{$\strat^i_{k+1} = \strat^i_{k} + \frac{1}{k+1} \big (
  \operatorname{br}^i(\AvgPayoff{i}(\strat^{-i}_k)) - \strat[i]_k
  \big)$,} where $\{\beta_k\}_{k \geq 1} \subset (0, 1)$ is a sequence
of stepsizes. This can be interpreted as a discretization of the
continuous-time best-response dynamics.

For the reasons outlined in~\Cref{SecSmoothedBR}, we study the closely
related $\tau$-smoothed best response dynamics given by
\begin{align}
  \label{EqnFullPolicy}
  \BStrat_{k+1} & = \BStrat_k + \beta_k \left
  (\sigma_{\tau}(\BStrat_k) - \BStrat_k \right ) \qquad \mbox{where
    $\sigma_{\tau}(\BStrat_k) = (\regBR(\strat^{-i}_k))_{i \in
      \nodes}$,}
\end{align}
where $\tau > 0$ is an algorithmic parameter.  These updates
correspond to a damped version of the usual operator power
method\footnote{For sufficiently large values of $\tau$, the operator
$\regBR[]$ is a contraction, in which case convergence of the update
$\BStrat_{k+1} = \regBR[](\BStrat_k)$ follows easily.  However, given
that we use $\tau$-regularized NE as an approximation to NE, our focus
is the more challenging small $\tau$-regime, in which this contractive
properly need not hold.}  attempting to find a $\tau$-regularized Nash
equilibrium~\eqref{EqnTauRegNash} and hence can be considered an
approximation of a greedy best-response. These updates also have an
interpretation involving players' beliefs as in stochastic fictitious
play (see the paper~\cite{hofbauersandholm2002} for more details). We
give a pseudocode specification of the full information dynamics
in~\Cref{AlgFull}.

\begin{algorithm}
\caption{Learning dynamics for full information}
\label{AlgFull}
\begin{algorithmic}
  \Require $K$, $\strat[i]_1 \sim \text{Unif}(\actions^i)$,
  temperature $\tau$ and stepsizes $\big \{ \beta_i\big
  \}_{i=1}^\infty$ \For{$k=1,\ldots,K$} \State $\BStrat_{k+1} =
  \BStrat_k + \beta_k \left (\sigma_{\tau}(\BStrat_k) - \BStrat_k
  \right ) $
\EndFor
\Ensure $\BStrat_{K+1}$
\end{algorithmic}
\end{algorithm}


\subsubsection{Dynamics for \infosetuplower}

In the \infosetuplower setting, players only have access to the
realized payoffs of their own actions. Notably, they do not observe
the opponents' strategies or actions and are not privy to any of the
payoff matrices. If each player $j$ chooses action $\action^j \in
\actions^j$, player $i$ observes $\sum_{j \in \nodes}
\RMatrix{i}{j}(\action^i,\action^j)$, the sum of the payoffs from each
of the pairwise games that they play.  Consequently, it is not
possible to directly compute the exact expected payoff vector for
every joint strategy.

In view of these restrictions, a natural solution (e.g., see the
papers~\cite{lesliecollins2005,sayin2021decentralized,zaiwei2023})
involves each player estimating the expected payoff vector for the
strategy chosen by their opponents in the previous round. In
particular, if action $\BAction$ is chosen according to the joint
strategy $\BStrat$, then (conditionally on $\BStrat$) the random
variable $\tfrac{\mathbb{I}(\action^i = \action
  )}{\strat[i](\action^i)}\sum_{j \in \nodes}
\RMatrix{i}{j}(\action^i,\action^j)$ is an unbiased estimate of the
average payoff $\AvgPayoff{i}(\strat[-i])(\action)$ for each $\action
\in \actions^i$.  We use this idea to sequentially build
approximations of the expected payoff vectors $\BAvgPayoff(\BStrat_k)$
for the joint strategy $\BStrat_k$ employed in round $k$. Note that
conditioned on past history, each player plays an independent strategy
at each step in this setting.

For each player $i \in \nodes$ and action $\action^i \in \actions^i$,
let $e^i(\action^i) \in \R^{|\actions^i|}$ be the standard basis
vector with a one in the position indexed by $\action^i \in
\actions^i$.  Consider the updates
\begin{align}
 \label{EqnTDUpdateSingle}
q^i_{k+1} & = q^i_k - \alpha_k \:
\frac{{e}^i(\action^i_k)}{\strat^i_k(\action^i_k)} \big ( \sum_{j \in
  \nodes} \RMatrix{i}{j}(\action^i_k, \action^{j}_k) -
q^i_k(\action^i_k )\big ) \qquad \mbox{for $k = 1, 2, \ldots$,}
\end{align}
where $\alpha_k > 0$ denotes a positive stepsize and $\action^j_k$ is
the action chosen by player $j \in \nodes$ according to their strategy
$\strat^j_k$. This is an asynchronous update rule where only the
component corresponding to the action $\action^i_k$ is updated. Note
that in order to compute this update, the only additional information
required by player $i$ is the payoff $\sum_{j \in \nodes}
\RMatrix{i}{j}(\action^i_k, \action^{j}_k)$ received by playing action
$\action^i_k$. We re-iterate that player $i$ does \emph{not} require
knowledge of the payoff matrices nor the actions of their opponent.

Equivalently, the update~\eqref{EqnTDUpdateSingle} can be rewritten as
\begin{align*}
q^i_{k+1}
      \; = \; q^i_k - \alpha_k \: \frac{E^i(\action^i_k)
      }{\pi_k^i(\action^i_k)} \Big \{ \sum_{j \in \nodes}
      \RMatrix{i}{j} \boldsymbol{e}^j(\action^j_k) - q^i_k \Big \},
\end{align*}
where $E^i(a^i_k) \defn
{e}^i(\action^i_k){{e}^i(\action^i_k)}^{\top}$, i.e., the matrix with
zeros in all entries except for a one in diagonal entry
$(\action^i_k,\action^i_k)$. By concatenating the iterates of all
players as $\BAvgPayoff_k = (q^i_k)_{i \in \nodes}$, we have the
combined update rule
\begin{subequations}
\begin{align}
 \label{EqnTDUpdate}
\BAvgPayoff_{k+1} & = \BAvgPayoff_k - \alpha_k \:
\frac{\boldsymbol{E}(\BAction_k)}{\BStrat_k(\BAction_k)} \big (
\BRMatrix \ \boldsymbol{e}(\BAction_k) - \BAvgPayoff_k \big ) \qquad
\mbox{for $k = 1, 2, \ldots$,}
\end{align}
where $\BAction_k = (\action^i_k)_{i \in \nodes}$ denotes the action
profile sampled from $\BStrat_k$, i.e., for each $i \in \nodes$,
$\action^i_k$ is sampled according to $\strat^i_k$,
$\boldsymbol{e}(\boldsymbol{\action}_k) = (e^i(\action^i_k))_{i \in
  \nodes}$ is the concatenation of the action vectors of all players
and
$\boldsymbol{E}(\boldsymbol{\action}_k)/\BStrat_k(\boldsymbol{\action}_k)$
is a diagonal block matrix whose $i^{\text{th}}$ diagonal block is
$E^i(\action^i_k)/\strat^i_k(\action^i_k)$ for every $i \in \nodes$.

Since the update~\eqref{EqnTDUpdate} is in the spirit of TD
Learning~\cite{sutton1988learning,zaiwei2023}, we refer to it as a TD
update.  The rescaling of the stepsize by $\strat[i](\action^i_k)$ for
each player $i$ in this update was originally proposed by Leslie and
Collins~\cite{lesliecollins2005}.  It can be understood as a form of
importance reweighting designed to ensure that $\E \big [ q^i_{k+1} -
  q^i_k \ \vert \ \BStrat_k, q^i_k \big ] = \alpha_k ( q^i_k -
\AvgPayoff{i}(\strat^{-i}_k))$, so that in expectation, the updates
are synchronous, and for a fixed constant strategy $\BStrat_k \equiv
\tilde{\BStrat}$, the TD updates have $\BAvgPayoff(\tilde{\BStrat})$
as a fixed point. We often refer to the iterates of the TD updates as
$q$-values.

\begin{algorithm}
\caption{Learning dynamics for the \infosetuplower setting}
\label{AlgPartial}
\begin{algorithmic}
  \Require $K, \,\BAvgPayoff_1=\mathbf{0} \in \otimes_{i \in \nodes}
  \R^{|\A^i|}, \, \strat^i_1 \sim \text{Unif}(\A^i), \tau,
  \timescalesep, \big \{ \beta_i\big \}_{i=1}^\infty$
  \For{$k=1,\ldots,K$}
  \State Play $\action^i_k \sim \strat^i_k$ independently for $i \in
  \nodes$
  \State $\BStrat_{k+1} = \BStrat_k + \beta_k \left
  (\sigma_{\tau}(\BAvgPayoff_k) - \BStrat_k \right )$
  \State $\BAvgPayoff_{k+1} = \BAvgPayoff_k - \alpha_k \:
  \frac{\boldsymbol{E}(\BAction_k)}{\BStrat_k(\BAction_k)} \big (
  \BRMatrix \ \boldsymbol{e}(\BAction_k) - \BAvgPayoff_k \big )$
  \EndFor
  \Ensure $\BStrat_{K+1} $
\end{algorithmic}
\end{algorithm}
In parallel with the TD updates for the expected payoff vectors, the
players update their strategies according to the smoothed
best-response dynamics
\begin{align}
\label{EqnPartialStrategy}
\BStrat_{k+1} & = \BStrat_k + \beta_k \big (\regBR[](\BAvgPayoff_k) -
\BStrat_k \big ) \qquad \mbox{where $\regBR(\BAvgPayoff_k) =
  (\regBR(\AvgPayoff{i}_k))_{i \in \nodes}$.}
\end{align}
\end{subequations}
The smoothed best-response plays an important role here: it ensures
that each action is played infinitely often as is required for TD
convergence~\cite{sutton1988learning}.  As noted earlier, smoothed
best-response updates ensure that this property
holds. See~\Cref{AlgPartial} for a pseudocode description of our
dynamics for the \infosetuplower setting.

We note that the TD updates combined with the smoothed best-response
dynamics form a coupled two-timescale stochastic approximation scheme
(e.g.,~\cite{Borkar2008}).  In such schemes, the iterates of the
update with the larger stepsize (often referred to as the ``faster''
timescale) has a fixed point that depends on the current value of the
iterates on the ``slower'' timescale, which evolves with a smaller
stepsize. The larger stepsize enables the faster-timescale iterates to
effectively ``track'' their moving fixed point. Similarly, the
strategies are updated slower than the $q$-values so that the TD
updates can track their fixed point. We achieve this timescale
separation by setting $\beta_k = \timescalesep \alpha_k$ for some
scalar $\timescalesep \in (0, 1)$ to be specified in our analysis.


\subsection{Lyapunov function}
\label{SecLyap}

In order to analyze the learning dynamics in~\Cref{AlgFull}
and~\Cref{AlgPartial}, we make use of a novel Lyapunov function for
the strategies, given by
\begin{align}
\label{EqnNovelLyapunov}
\stratLyap(\BStrat) & \defn \sum_{i=1}^{\NumPlayers}
\underset{\hat{\strat} \in \simplex{i}}{\max} \left\{
         {\hat{\strat}}^\top \AvgPayoff{i}(\strat^{-i}) + \tau
         \EntFun(\hat{\strat})\right\} ,
\end{align}
where $H$ denotes the Shannon entropy function~\eqref{EqnDefnShannon}.
Note that the Lyapunov function $\stratLyap$ has a natural
game-theoretic interpretation: it corresponds to the sum of the
average payoffs when each player chooses the best-response to the
entropy-regularized payoffs.  It can be viewed as a regularized
version of the function used by Harris~\cite{harris1998} to prove
finite-time guarantees for continuous-time zero-sum games. The main
purpose of our Lyapunov function~\eqref{EqnNovelLyapunov} is to
provide a way to bound the Nash gap from above in our proofs; more
precisely, they are related by the inequality
\begin{align}
\label{eqn:nash_bound}
\nashg(\BStrat) \leq \stratLyap(\BStrat).
\end{align}
The Lyapunov function $\stratLyap$ differs from other Lyapunov
functions used in the literature~\cite{harris1998,zaiwei2023} in some
crucial ways. Its smoothness properties enable us to prove a
polynomial-time finite-sample guarantee (as opposed to the
exponential-time guarantee in Chen et al.~\cite{zaiwei2023}). We
discuss this issue in further detail in the supplementary
Section~\ref{sec:lyapunov_remarks}.

In addition, our proof also involves the
\emph{$q$-Lyapunov function}
\begin{align}
	\label{EqnQLyap}
	\TDLyap(\BStrat, \bq) & \defn \sum_{i\in\nodes} \big \| q^i -
	\AvgPayoff{i}(\strat^{-i}) \big \|_2^2,
\end{align}
which we use to track the quality of $q^i$ as an estimate of
$\AvgPayoff{i}(\strat^{-i})$.  The \emph{total Lyapunov function}
$\TotLyap$ for our analysis is given by the sum 
\begin{align}\label{EqnTotLyap}
	\TotLyap(\BStrat,
	\bq) \defn \stratLyap(\BStrat) + \TDLyap(\BStrat, \bq).
\end{align}
\section{Main results for zero-sum polymatrix games}

In this section, we state our results for zero-sum polymatrix games
for both information settings. \Cref{SecFullPoly} provides guarantees
for the full information dynamics (cf.~\Cref{AlgFull}) whereas
in~\Cref{SecMinPoly}, we study the minimal information two-timescale
dynamics (cf.~\Cref{AlgPartial}). In stating and proving these
results, al variables of the form $c_i, i \in \N$ denote numerical
constants that are independent of both the game and the learning
dynamics. Our results can also be specialized for the class of
two-player zero-sum games, with more detail given in the supplementary
Section~\ref{SecFull2P}.


\subsection{Full information}
\label{SecFullPoly}

In this section, we provide finite-sample guarantees for the full
information dynamics specified
in~\Cref{AlgFull}. \Cref{ThmFullInfoPoly} provides an upper bound on
the Nash gap for three choices of stepsize schedules. In stating these
claims, we make use of the shorthand $\Vinit \defn
\stratLyap(\BStrat_1)$ for the initial value of the Lyapunov function
$\stratLyap$ from equation~\eqref{EqnNovelLyapunov}, where $\BStrat_1$
is the initial set of mixed strategies.

\begin{theorem}[Nash gap finite-sample guarantees]
\label{ThmFullInfoPoly}
For any $\tau > 0$, consider the full information dynamics
(\Cref{AlgFull}) initialized with $\BStrat_1$.  Then the Nash gap
after $K$ iterations is bounded as:
 \begin{enumerate}[label=(\alph*), leftmargin=1.8em]
 \item For a constant stepsize sequence $\beta_k \equiv \beta \in
   (0,1)$, we have
   \begin{align*}
        \nashg(\BStrat_{K+1}) \leq (1-\beta)^K \StratLyapVal{1} + \tau
        \NumPlayers \log \Amax + \frac{\NumPlayers
          \opnorm{\BRMatrix}^2 \beta }{\tau}.
   \end{align*}
\item For the inverse linear stepsize sequence $\beta_k =
  \tfrac{\beta}{k}$ for some $\beta \in (1,2]$, we have
\begin{align*}
  \nashg(\BStrat_{K + 1}) & \leq \frac{\StratLyapVal{1}}{(K+1)^\beta}
  + 8 \NumPlayers \tau \log \Amax + \frac{4 \NumPlayers
    \opnorm{\BRMatrix}^2 \beta^2}{\tau (\beta-1) K} .
\end{align*}
\item For the inverse polynomial stepsize
sequence $\beta_k = \tfrac{\beta}{(k+\koff)^\myexp}$ with
$\beta,\myexp \in (0,1)$ and $\koff \geq \lceil \left (
\tfrac{1-\myexp}{\beta}\right )^{1/\eta} \rceil$, we have
\begin{multline*}
	\nashg(\BStrat_{K+1}) \leq \exp \big( - \frac{\beta}{1-\myexp} \big(
	(K+\koff+1)^{1-\myexp} - (1+\koff)^{1-\myexp} \big)\big )
	\stratLyapVal{1} \\
	+ \NumPlayers \tau \log \Amax+ \frac{2 \NumPlayers
		\opnorm{\BRMatrix}^2 \beta}{\tau (K + \koff)^\myexp}.
\end{multline*}
 \end{enumerate}
\end{theorem}

For all three stepsizes, the first term in the upper bound involves
the initial Lyapunov value $\StratLyapVal{1}$, and so reflects the
rate at which the algorithm ``forgets'' its initialization as it
converges. The second term in each bound scales linearly in $\tau$,
and corresponds to a form of bias introduced by the players using a
$\tau$-regularized best-response instead of an actual
best-response. The third term in each upper bound scales with
$(1/\tau)$, which is a measure of how smooth the $\tau$-regularized
best-response is.

For the purposes of interpretation, it is useful to derive bounds on
on the iteration complexity of the procedure.  For a given level
$\epsilon > 0$, the \emph{iteration complexity} $K(\epsilon)$ is the
minimum number of iterations required to ensure that
$\nashg(\BStrat_{K(\epsilon)+1}) \leq \epsilon$. In order to obtain
explicit bounds on $K(\epsilon)$, we choose the temperature parameter
$\tau$ and stepsizes so as to ensure that after $K(\epsilon)$ rounds,
each of the three terms in the upper bound in~\Cref{ThmFullInfoPoly}
is at most \mbox{$\frac{\epsilon}{3}$}.  By doing so, we obtain the
following:
\begin{corollary}
\label{CorFullInfoPoly}
Consider the full information dynamics (\Cref{AlgFull}) initialized
with \mbox{$\tau = \newlcc{} \epsilon/(\NumPlayers\log \Amax)$}.  Then
the iteration complexity $K(\epsilon)$ is bounded as follows:
\begin{subequations}
\begin{enumerate}[label=(\alph*)]
\item For the constant stepsizes $\beta_k \equiv \beta \defn
  \epsilon^2/(\newlcc{} \NumPlayers^2 \opnorm{\BRMatrix}^2
  \log\Amax)$, we have
  \begin{align*}
    K(\epsilon) \leq \frac{\newlcc{} \NumPlayers^2
      \opnorm{\BRMatrix}^2\log\Amax}{\epsilon^2}\log \big(
    \frac{\StratLyapVal{1}}{\epsilon} \big).
  \end{align*}
 \item For the inverse linear decay $\beta_k = \beta/k$ for some
   $\beta \in (1,2]$, we have
 \begin{equation*}
   K(\epsilon) \leq \frac{\newlcc{} \NumPlayers^2 \opnorm{\BRMatrix}^2
     \beta^2 \log\Amax \StratLyapVal{1}}{(\beta-1)\epsilon^2}.
 \end{equation*}
\item 	For the inverse polynomial stepsize
sequence $\beta_k = \tfrac{\beta}{(k+\koff)^\myexp}$ with
$\beta,\myexp \in (0,1)$ and $\koff \geq \lceil \left (
\tfrac{1-\myexp}{\beta}\right )^{1/\eta} \rceil$, we have
\begin{align*}
	K(\epsilon) \leq \big ( k_0^{1-\myexp} + \frac{(1-\myexp)}{\beta} \log
	\frac{\stratLyapVal{1}}{\epsilon} \big )^{1/(1-\myexp)} + \big (\frac{
		\newlcc{} \NumPlayers^2 \opnorm{\BRMatrix}^2 \beta
		\log\Amax}{\epsilon^2} \big )^{1/\myexp} .
\end{align*}
\end{enumerate}
\end{subequations}
\end{corollary}

Focusing on the triple $(\NumPlayers, \opnorm{\BRMatrix}, \epsilon)$,
our theory guarantees that smooth best-response dynamics has iteration
complexity bounded as
$\mathcal{O}(\opnorm{\BRMatrix}^2\NumPlayers^2/\epsilon^2)$.  As
mentioned previously, if one allows for different types of algorithms,
it can be possible to obtain different scalings of the iteration
complexity.  For example, Ao et al.~\cite{AoCenChi2022asynchronous}
studied updates based on a multiplicative weights (MW) update, and
analyzed its behavior in terms of the maximum error over all
players, as opposed to the sum of errors~\eqref{EqnNashGap} in our analysis.  For this
error metric and algorithm, they obtained an iteration complexity
scaling as as $\mathcal{O}(\dmax
\matsnorm{\BRMatrix}{\max}/\epsilon)$, where
$\matsnorm{\BRMatrix}{\max}$ is the maximum absolute element of the
matrix $\BRMatrix$, and $\dmax$ is the maximum degree of the graph.
Since $\dmax \leq \NumPlayers$ and $\matsnorm{\BRMatrix}{\max} \leq
\opnorm{\BRMatrix}$, this MW iteration complexity---albeit for a
different error metric---represents a quadratic improvement over that
of smoothed best-response.  Apart from the error metrics (which differ
by a factor of $\NumPlayers$ in an extreme case), we attribute this
gap to the linear nature of our best-response updates as opposed to
their multiplicative updates, which are known to work well with KL
divergence-based Lyapunov functions.

It is notable that our bounds scale (quadratically) in the spectral
norm $\opnorm{\BRMatrix}$, a quantity which depends on subtle ways on
both the graph structure and the pairwise games on each edge.  In this
way, our analysis affords some insight into the type of zero-sum
polymatrix games for which it is easier to converge to Nash
equilibria. We discuss the structural properties of the parameter
$\opnorm{\BRMatrix}$ in further detail in~\Cref{SecRprops}. 


\subsection{Minimal information} \label{SecMinPoly}

We now turn to the analysis of the updates in~\Cref{AlgPartial} that
apply to the minimal information setting.  Since this is a stochastic
algorithm, our bounds apply to the \emph{iteration complexity}
$K(\epsilon)$ defined by the minimum number of rounds required to
ensure that $\E \big[\nashg(\BStrat_{K(\epsilon) + 1}) \big] \leq
\epsilon$.  We again show that the iteration complexity scales
polynomially in $(1/\epsilon)$--in this case, we can guarantee a
scaling of the order $(1/\epsilon)^{8 + \offset}$ for an exponent
parameter $\offset > 0$ that can be chosen arbitrarily close to zero.
The price of taking $\offset \rightarrow 0^+$ manifests in the growth
of certain pre-factors; we use functions of the form $g(\offset)$ and
variants thereof to indicate terms of this type. We make use of the shorthand $\TotLyapVal{1}$ for the initial value of the Lyapunov function
$\TotLyap(\BStrat_1,\bq_1)$ from equation~\eqref{EqnTotLyap}, where $\BStrat_1$
is the initial set of mixed strategies and $\bq_1$ is the initial estimate of $\BAvgPayoff(\BStrat_1)$. 

Our result applies to~\Cref{AlgPartial} where the temperature and the
timescale separation constant are set as
\begin{equation}
\label{eqn:tau_c_initPoly}
\tau = \frac{g_{\tau}(\offset)\epsilon}{ \NumPlayers \log \Amax} \quad
\text{and} \quad \timescalesep = \frac{g_{\alpha,\beta}(\offset)
  \tau^3}{\opnorm{\BRMatrix}^2} \quad \text{, respectively,}
\end{equation}
and the function $g_{\alpha, \beta}$ satisfies the scaling
$g_{\alpha,\beta}(\offset) \rightarrow 0^+$ as $\offset \rightarrow
0^+$.  Our guarantees holds in terms of a triple of functions $(g_1,
g_2,g_3)$ such that $\max_{j = 1, 2, 3 } g_j(\offset) \rightarrow 0^+$
as $\offset \rightarrow 0^+$.
\begin{theorem}
\label{ThmMinInfoPoly}
Consider the minimal information dynamics (\cref{AlgPartial})
initialized with the parameters~\eqref{eqn:tau_c_initPoly}.  Then for
any $\offset > 0$ and in each of the following cases, the iteration
complexity $K(\epsilon)$ satisfies the following upper bounds:
\begin{enumerate}[label=(\alph*)]
 \item For the constant stepsize $\beta_k \equiv \beta \defn
   \frac{g_1(\offset)\epsilon^{8+\offset}}{\Amax^6\NumPlayers^8\opnorm{\BRMatrix}^6}$,
we have
 \begin{align*}
K(\epsilon) \leq
\frac{\Amax^6\NumPlayers^8\opnorm{\BRMatrix}^6}{g_2(\offset)\epsilon^{8+\offset}}
\log \big ( \frac{3 \TotLyapVal{1}}{\epsilon}\big ).
 \end{align*}
\item For the inverse polynomial stepsize $\beta_k =
  \frac{\beta}{(k+\koff)^\myexp}$ for some exponent $\myexp \in
  (0,1)$, offset \mbox{$\koff = \big \lceil
    (\frac{2\myexp}{\beta})^{1/(1-\myexp)} \big \rceil$,} and
  $\beta=\frac{g_1(\offset)\epsilon^{8+\offset}}{
    \Amax^6\NumPlayers^8\opnorm{\BRMatrix}^6}$, we have
\begin{align*}
K(\epsilon) & \leq  \Big \{
\frac{(1-\myexp) \Amax^6 \NumPlayers^8
  \opnorm{\BRMatrix}^6}{g_3(\offset)\epsilon^{8+\offset}} \log \big (
\frac{3 \TotLyapVal{1}}{\epsilon}\big ) \Big \}^{\frac{1}{1 -
    \myexp}}.
\end{align*}
\end{enumerate}
\end{theorem}

As discussed previously, it is key that all the guarantees
in~\Cref{ThmMinInfoPoly} scale polynomially in $(1/\epsilon)$---in
particular, as $(1/\epsilon)^{8 + \offset}$ along with additional
logarithmic factors.  To the best of our knowledge, these are the
first known polynomial guarantees for the standard best-response
dynamics in this setting.  One of the main challenges in establishing
this polynomial scaling is controlling the variance of the
$q$-updates.  This variance depends on the probability of the actions
chosen in each round, and choosing actions based on softmax response
leads to probabilities that are \emph{exponentially small} in $\tau$.
Consequently, a naive approach yields an exponential dependence on
$1/\tau$, and hence---since our choice of $\tau$ is proportional to
the target accuracy $\epsilon$---an exponential dependence on
$(1/\epsilon)$.  Notably, the analysis in some past
work~\cite{zaiwei2023} exhibits this type of exponential growth. In
our analysis, we resolve this issue by initializing the dynamics away
from the boundary, and then choosing the stepsize parameter in a
way that allows us to control how quickly the iterates approach the
boundary of the probability simplex. See~\Cref{SecProofThmMinInfo} for
the details of this argument.

The stepsize choices in~\Cref{ThmMinInfoPoly} depend on the functions
$g_1$ and $g_2$ of $\offset$, which have explicit expressions given in
the proof of~\Cref{ThmMinInfoPoly} in~\Cref{SecProofThmMinInfo}.  The
stepsizes also depend on the action sizes, number of players, and the
spectral norm $\opnorm{\BRMatrix}$.  The latter parameter is global in
nature, requiring each player to be aware of the full payoff structure
of the polymatrix game.  This requirement can be relaxed by bounding
$\opnorm{\BRMatrix}$ from above: for example, since the block matrix
$\BRMatrix$ is formed by concatenating the payoff matrices
$\RMatrix{i}{j}$ for game $(i,j)$, we have the bound
\begin{align*}
  \opnorm{\BRMatrix} \leq \max_{i \in \nodes} \sum_{j \in \nodes}
  \opnorm{\RMatrix{i}{j}} \; \leq \; \dmax \max_{i, j}
  \opnorm{\RMatrix{i}{j}} \; \leq \; \Amax \dmax
  \matsnorm{\BRMatrix}{\max},
\end{align*}
where $\matsnorm{\BRMatrix}{\max}$ denotes the maximum absolute entry
of the matrix $\BRMatrix$.

\subsection{Studying the spectral norm}
\label{SecRprops}
The spectral norm $\opnorm{\BRMatrix}$ captures the dependence of the
rate of convergence on the underlying graph and the pairwise
games. Studying the parameter $\opnorm{\BRMatrix}$ can reveal answers
to various qualitative questions of interest.  For instance, what
types of graphs and pairwise games admit the fastest convergence to
Nash equilibria for best-response dynamics?

In order to address this question, let us consider the special class of zero-sum polymatrix
games in which each edge is associated with the \emph{same} two-player
zero-sum game.  Given a graph $\graph = (\nodes, \edges)$, define a
skew-symmetric weighted adjacency matrix $\Amat \in \{-1, 0,1
\}^{\NumPlayers \times \NumPlayers}$ where the $(i,j)^{\text{th}}$ entry is non-zero only if $(i,j) \in \edges$.
Letting $\PayOff \in \real^{m \times m}$ be a symmetric matrix of
payoffs, we associate with edge $(i,j)$ the zero-sum game defined by
the payoff matrix $A_{ij} \PayOff$ for player $i$, and $-A_{ij}
\PayOff$ for \mbox{player $j$.}

A useful property of this zero-sum polymatrix game is that the block
matrix $\BRMatrix$ can be written as $\BRMatrix = \Amat \otimes
\PayOff$, where $\otimes$ is the Kronecker product.  It follows that
$\opnorm{\BRMatrix} = \opnorm{\Amat} \opnorm{\PayOff}$, so that this
decomposition separates out the dependency on the game payoff
$\PayOff$ and the graph structure $\Amat$.   We now discuss these
two factors in turn.

\subsubsection{Graph dependence}

Recall the matrix norms \mbox{$\matsnorm{\Amat}{\infty} \defn
  \matsnorm{\Amat^\top}{1}$} and $\matsnorm{\Amat}{1} \mbox{\defn}
\max \limits_{j \in \nodes} \sum_{i \in \nodes} |A_{ij}|$. We always
have the elementary upper bound $\opnorm{\Amat} \leq \max \big \{ \matsnorm{\Amat}{1},
\; \matsnorm{\Amat}{\infty} \big \} \leq \dmax.$
For graphs with particular structure, this bound can be tightened:
we now derive bounds $\opnorm{\Amat}$ that reveal its exact dependence
on $\dmax$ for two sub-classes of graphs.

\begin{wrapfigure}{l}{0.3\textwidth}
	\centering
	\begin{subfigure}[b]{\linewidth}
		\centering
		\includegraphics[width=0.8\linewidth]{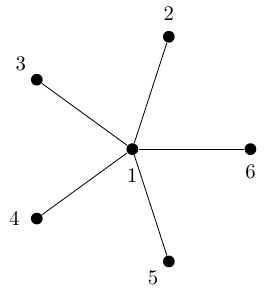}
		\caption{Star-connected graph}
		\label{FigStar}
	\end{subfigure}
	
	\vspace{0.5em} 
	
	\begin{subfigure}[b]{\linewidth}
		\centering
		\includegraphics[width=0.8\linewidth]{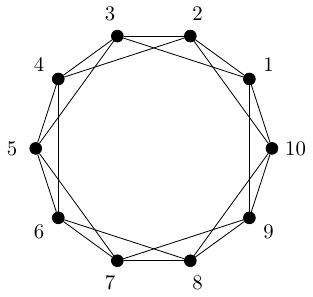}
		\caption{4-regular ring graph}
		\label{FigRing}
	\end{subfigure}
	
	\caption{Graph structures.}
	\label{fig:stacked_graphs}
\end{wrapfigure}

\myunderparagraph{Star-connected graph.} In this type of graph, a
central node---say node $1$---is connected to all the other nodes,
yielding a graph with maximum degree $\dmax = \NumPlayers - 1$.  With
node $1$ as the central node, only the first row and first column of
$A$ have non-zero entries.  \Cref{FigStar} provides an illustration of
a star graph with $\NumPlayers = 6$ nodes; an example of the
associated skew-symmetric adjacency matrix is given by
\begin{align*}
\Amat & \defn \left [ \begin{matrix} 0 & 1 & -1 & -1 & 1 & 1 \\ -1 & 0
    & 0 & 0 & 0 & 0 \\ 1 & 0 & 0 & 0 & 0 & 0 \\ 1 & 0 & 0 & 0 & 0 & 0
    \\ -1 & 0 & 0 & 0 & 0 & 0 \\ -1 & 0 & 0 & 0 & 0 & 0
    \end{matrix} \right ].
\end{align*}
Introducing the shorthand $\Mmat \defn \Amat^\top \Amat$, we have
$\Mmat_{11} = \NumPlayers-1$ and $\Mmat_{1j} = \Mmat_{j1} = 0$ for all
$j \neq 1$. Therefore, the matrix $\Mmat$ has $\NumPlayers-1$ as an
eigenvalue. The submatrix of $\Mmat$ formed by excluding the first
row and column has rank 1 with the only non-zero eigenvalue being
equal to $\NumPlayers-1$. It follows that $\opnorm{\Amat} =
\sqrt{\NumPlayers -1 }$, i.e., $\opnorm{\BRMatrix}$ scales as
$\sqrt{\dmax}$.

\myunderparagraph{$k$-regular ring graph.} In a $k$-regular ring graph
(see~\Cref{FigRing} for an example), every node is connected to $k$
other nodes. The weighted adjacency matrix $\Amat$ in this case is
circulant, i.e., the entire matrix can be generated by shifting the
first row. It is possible to explicitly write down the imaginary
eigenvalues of $\Amat$ and show that $\opnorm{\Amat}$ grows linearly
with $k$, i.e., grows linearly with $\dmax$
(see~\Cref{SeckRegularRMatrix} for details).

To summarize, we see that for graphs with a fixed number of players
$\NumPlayers$, the rate of convergence to approximate Nash equilibria
of the best-response type dynamics we analyze does scale with the
maximum degree of the graph $\dmax$. A similar trend is documented in
the paper~\cite{AoCenChi2022asynchronous} using an algorithm based on
the multiplicative weights method; they observed that the iteration 
complexity for their algorithm grows with $\dmax$. Our bounds exhibit
a finer dependence on $\dmax$ through the quantity
$\opnorm{\BRMatrix}$ in the sense that it enables us to compare
different types of graphs. For instance, based on the preceding
analysis, our bounds indicate that for a large enough number of
players $\NumPlayers$, a $k$-regular graph with $k >
\sqrt{\NumPlayers-1}$ converges more slowly than a star-connected
graph with the same number of nodes, which in turn converges more
slowly than a $k$-regular graph with $k<\sqrt{\NumPlayers-1}$.

\subsubsection{Dependence on pairwise games}

Finally, we turn to the effect of the underlying pairwise games on the
rate of convergence.  Let $\Amat$ be the skew-symmetric adjacency
matrix of a $k$-regular ring graph, and let $\PayOff$ be a fixed
payoff matrix.  Each node $i \in \nodes$ is associated with $k$ payoff
matrices: for some scalar $\SmallEdgeFrac \in [0,1]$, suppose that at
most $\SmallEdgeFrac k$ of these payoff matrices are set equal to
$\PayOff$, whereas the remaining payoffs are set equal to
$\SmallEdgeScaling \PayOff$ for a weight parameter $\SmallEdgeScaling
\in [0, 1]$.

The parameter $\SmallEdgeFrac$ controls the number of ``strong''
pairwise matrix games, whereas the parameter $\SmallEdgeScaling$ is a
measure of the ratio between weak and strong games.  Using the
techniques described previously, the spectral norm
$\opnorm{\BRMatrix}$ of the block matrix $\BRMatrix$ can be bounded as
\begin{align*}
\opnorm{\BRMatrix} & \leq (\SmallEdgeFrac k) \opnorm{\PayOff} + k (1 -
\SmallEdgeFrac) \opnorm{\SmallEdgeScaling \PayOff} \; = \; k
\opnorm{\PayOff} \big \{ \SmallEdgeFrac + (1- \SmallEdgeFrac)
\SmallEdgeScaling \big \}.
\end{align*}
Thus, the dependence of $\opnorm{\BRMatrix}$ on the maximum degree $k$
can be made arbitrarily small by suitable choices of the two
parameters $(\SmallEdgeFrac, \SmallEdgeScaling)$.


\section{Proofs of main results}
\label{SecProofs}

This section is to the proofs of our main results, with
Subsections~\ref{SecProofThmFullInfo} and~\ref{SecProofThmMinInfo}
devoted to the proofs of Theorems~\ref{ThmFullInfoPoly}
and~\ref{ThmMinInfoPoly}, corresponding to full information and
minimal information cases, respectively.

\subsection{Proof of~\Cref{ThmFullInfoPoly}}
\label{SecProofThmFullInfo}

It is convenient to introduce the shorthand \mbox{$\StratLyapVal{k}
  \defn \stratLyap(\BStrat_k)$} for the value of the Lyapunov function
at iteration $k$.  All three sub-claims in~\Cref{ThmFullInfoPoly} are
derived via a drift inequality for $\stratLyap$: in particular, we
claim that
\begin{align}
\label{EqnDriftStratLyap}
\StratLyapVal{k+1} & \leq (1-\beta_k) \StratLyapVal{k} + \beta_k
\tau \NumPlayers \log \Amax + \tfrac{\NumPlayers
\opnorm{\BRMatrix}^2 \beta_k^2 }{\tau}, \qquad \mbox{for
$k=1,2,\ldots$}.
\end{align}
See~\Cref{SecDriftProof} for the proof of this claim.

\subsubsection{From drift inequality to Nash gap}
Solving the inequality~\eqref{EqnDriftStratLyap} for each type of stepsize gives an upper bound on $\StratLyapVal{K+1}$. Via
inequality~\eqref{eqn:nash_bound}, these bounds translate to upper
bounds on the Nash gap, and hence the iteration complexity.  Here we
provide this argument for the constant stepsizes; see
Section~\ref{AppFullStep} for the remaining cases.

In the constant stepsize case, we have $\beta_k \equiv \beta \in
(0,1)$ for all iterations $k = 1, 2, \ldots$.  For this choice,
iterating the bound~\eqref{EqnDriftStratLyap} $K$ times
yields
\begin{align*}
\StratLyapVal{K+1} & \leq (1-\beta)^K \StratLyapVal{1} + \tau
\NumPlayers \log \Amax + \frac{4\NumPlayers \opnorm{\BRMatrix}^2
\beta }{\tau}.
\end{align*}
For a target error $\epsilon \in (0,1)$, setting $\tau
=\opnorm{\BRMatrix}\sqrt{\beta/\log\Amax}$ yields the bound
\begin{align*}
\StratLyapVal{K+1} & \leq (1-\beta)^{K} \StratLyapVal{1} + {2
\NumPlayers \opnorm{\BRMatrix} (\log \Amax)^{3/2}} \sqrt{\beta}.
\end{align*}
Thus, if we set $\beta = \epsilon^2/(16 \NumPlayers^2
\opnorm{\BRMatrix}^2 (\log\Amax)^3)$, we can conclude that it suffices
to take $K(\epsilon) = \frac{16 \NumPlayers^2 \opnorm{\BRMatrix}^2
  (\log\Amax)^3}{\epsilon^2}\log \big(
\frac{\StratLyapVal{1}}{\epsilon} \big)$ iterations in order to ensure
that $\StratLyapVal{K(\epsilon) + 1} \leq \epsilon$, which in turn
results in the Nash gap \mbox{being less than $\epsilon$.}

\subsubsection{Proof of drift inequality~\eqref{EqnDriftStratLyap}}
\label{SecDriftProof}
The following auxiliary result plays a key role in this proof (as
well as that of~\Cref{ThmMinInfoPoly}).
\begin{lemma}
  \label{LemVprops}
For the Lyapunov function $\stratLyap$ from
equation~\eqref{EqnNovelLyapunov}:
  \begin{enumerate}[label=(\alph*),leftmargin=1.8em]
  \item It is twice continuously differentiable, and the
    $\ell_2$-operator norm $\opnorm{\cdot}$ of its Hessian is
    uniformly bounded as
    \begin{subequations}
      \begin{align}
        \label{EqnStratLsmooth}
\opnorm{\nabla^2 \stratLyap(\BStrat)} & \leq L \defn
\frac{\opnorm{\BRMatrix}^2}{\tau} \qquad \mbox{for all $\BStrat \in
  \bigdelta$.}
      \end{align}
    \item \mbox{For all $\BStrat \in \bigdelta$,} we have the bound
      \begin{align}
\label{EqnStratNegativeDrift}
\sum_{i=1}^{\NumPlayers} \big \langle \nabla_{\strat[i]}
\stratLyap(\BStrat), \sigma_{\tau}(\AvgPayoff{i}(\strat^{-i})) -
\strat[i] \big \rangle &\leq - \stratLyap(\BStrat) + \NumPlayers \tau
\log \Amax.
      \end{align}
    \end{subequations}
  \end{enumerate}
\end{lemma}
\noindent See~\Cref{SecLemVprops} for the proof of this claim. \\

\noindent Using this auxiliary result, we now prove the
claim~\eqref{EqnDriftStratLyap}.  From the Hessian
bound~\eqref{EqnStratLsmooth} on $\StratLyap$, we have
\begin{align*}
  \underbrace{\stratLyap(\BStrat_{k+1})}_{\StratLyapVal{k+1}} & \leq
  \underbrace{\stratLyap(\BStrat_k)}_{\StratLyapVal{k}} +
  \inprod{\nabla \stratLyap(\BStrat_k)}{\BStrat_{k+1} -
    \BStrat_k}  + \frac{L}{2} \vecnorm{\BStrat_{k+1} - \BStrat_k}{2}^2.
\end{align*}
Using H\"{o}lder's inequality,
we can upper bound $\vecnorm{\BStrat_{k+1} - \BStrat_k}{2}^2$
by
\begin{align}\label{EqnStratOneStep}
  \beta_k^2 \sum_{i=1}^\NumPlayers
  \vecnorm{\sigma_{\tau}(\AvgPayoff{i}(\strat^{-i}_k) ) -
    \strat[i]_k}{2}^2 \; \leq \; \beta_k^2 \sum_{i=1}^\NumPlayers
  \|\sigma_{\tau}(\AvgPayoff{i}(\strat^{-i}_k) ) - \strat[i]_k\|_1^2
  \|\sigma_{\tau}(\AvgPayoff{i}(\strat^{-i}_k) ) -
  \strat[i]_k\|_\infty^2.
\end{align}
Since both $\sigma_{\tau}(R^i \pmi_k )$ and $\pii_k$ belong to the
probability simplex, we have $\|\sigma_{\tau}(R^i \pmi_k ) -
\pii_k\|_\infty \leq 1$ and $\|\sigma_{\tau}(R^i \pmi_k ) - \pii_k\|_1
\leq 2$, whence $\vecnorm{\BStrat_{k+1} - \BStrat_k}{2}^2 \leq 4
\NumPlayers \beta_k^2$.  Recalling that \mbox{$\BStrat_{k+1} -
  \BStrat_k = \beta_k \big(\sigma_{\tau}( \BAvgPayoff(\BStrat_k) ) -
  \BStrat_k \big)$,} observe that we have $\inprod{\nabla
  \stratLyap(\BStrat_k)}{\BStrat_{k+1} - \BStrat_k} = \beta_k \Big \{
\sum_{i=1}^{\NumPlayers} \big \langle \nabla_{\strat[i]}
\stratLyap(\BStrat_k), \sigma_{\tau}(\AvgPayoff{i}(\strat^{-i}_k)) -
\strat[i]_k \big \rangle \Big \}$. Putting together the pieces yields
the bound
\begin{align}
  \label{EqnFullInfoVDrift}
  \StratLyapVal{k+1} & \leq \StratLyapVal{k} + \beta_k \Big \{
  \sum_{i=1}^{\NumPlayers} \big \langle \nabla_{\strat[i]}
  \stratLyap(\BStrat), \sigma_{\tau}(\AvgPayoff{i}(\strat^{-i})) -
  \strat[i] \big \rangle \Big \} + \frac{4 \NumPlayers
    \opnorm{\BRMatrix}^2}{\tau} \beta_k^2.
\end{align}
It remains to bound the first-order term on the right-hand side of
inequality~\eqref{EqnFullInfoVDrift}, and we do so using
inequality~\eqref{EqnStratNegativeDrift} from~\Cref{LemVprops}.
Substituting this bound into inequality~\eqref{EqnFullInfoVDrift} and
re-arranging yields the claimed drift
inequality~\eqref{EqnDriftStratLyap}.

\subsection{Proof of~\Cref{ThmMinInfoPoly}}
\label{SecProofThmMinInfo}

We now turn to the proof of~\Cref{ThmMinInfoPoly} that applies to the
minimal information setting. We proceed in two parts: first, we
establish a drift inequality for the Lyapunov functions defined
in~\Cref{SecLyap} subject to a condition on how far away the strategy
vectors $\BStrat_k$ are from the boundary of the simplex. This is done
to control the variance of the $q$-updates. Solving this drift
inequality gives us a bound on the iteration complexity subject to the
aforementioned condition. In the second part of the proof, we show
that by choosing certain quantities related to the dynamics
appropriately, this condition will be satisfied.

Recall the definition~\eqref{EqnNovelLyapunov} of the
strategy-tracking Lyapunov function $\StratLyap$, the $q$-Lyapunov
function $\TDLyap$~\eqref{EqnQLyap} and the total Lyapunov function
$\TotLyap$~\eqref{EqnTotLyap}. Since our goal is to bound the Nash
gap, it is useful to make note of the sandwich relation
\begin{align}
  \label{eqn:nash_bound_M}
  \nashg(\BStrat) \stackrel{(i)}{\leq} \stratLyap(\BStrat) \;
  \stackrel{(ii)}{\leq} \TotLyap(\BStrat, \bq),
\end{align}
where inequality (i) follows from the bound~\eqref{eqn:nash_bound} and
inequality (ii) follows from the definition of the total Lyapunov
function and the non-negativity of $\TDLyap(\BStrat, \bq)$.  Thus, in
order to bound the Nash gap, it suffices to upper bound
$\TotLyap(\BStrat, \bq)$. Also recall that our algorithm generates a
sequence $\{\BStrat_q, \bq_k \}_{k \geq 1}$ of strategy and $q$-value
pairs.  For $k = 1, 2, \ldots$, we introduce the shorthand notation
$\StratLyapVal{k} \defn \StratLyap(\BStrat_k)$, $\TDLyapVal{k} \defn
\TDLyap(\BStrat_k, \bq_k)$, and $\TotLyapVal{k} \defn
\TotLyap(\BStrat_k, \bq_k)$ for the values of the three Lyapunov
functions as a function of the iteration number $k$.


\subsubsection{Drift inequality for total Lyapunov function $\TotLyap$} \label{SecLemTotDrift}

The first step in proving~\Cref{ThmMinInfoPoly} is to establish a drift
inequality (or recursive bound) on the expected values $\E
\TotLyapVal{k}$ of the total Lyapunov function $\TotLyap$ over
iterations $k$.  The main challenge in this step is that the variance
of the $q$-updates explodes exponentially in $(1/\tau)$ as the
strategy vectors $\BStrat_k$ approach the boundary of the probability
simplex.  So as to avoid this explosion, we need to track carefully
the rate at which the iterates approach the boundary.  For a tolerance
parameter $\delta \in (0,1)$, we say that the iterates $\{\BStrat_k
\}_{k \geq 1}$ are \textit{$\delta$-good} up to time $K$ if we have
\begin{align}
  \label{EqnDeltaGood}
  \min_{\action^i \in \A^i} \BStrat^i_k(\action^i) \geq \delta \qquad
  \mbox{for all iterations $k = 1, \ldots, K$, and nodes $i \in
    \nodes$.}
\end{align}

\begin{lemma}
\label{LemCoupledDrift}
Suppose that the iterates $\{\BStrat_k\}_{k \geq 1}$ are $\delta$-good
up to time $K$.  Then for each $k = 1, \ldots, K$, we have
\begin{multline}\label{EqnTotLyapDrift}
\E\TotLyapVal{k+1} \leq \big (1-\beta_k \big (1-2 \constYoung \big )
\big )\E \StratLyapVal{k} + \Big ( 1-\alpha_k +\alpha_k^2 \tfrac{\Amax
\NumPlayers }{\delta}+ \beta_k\tfrac{3\opnorm{\BRMatrix}^2}{
\constYoung \tau^3} \Big ) \E \TDLyapVal{k} \\+ 2\NumPlayers \tau
\beta_k \log \Amax + \tfrac{8
\NumPlayers\opnorm{\BRMatrix}^2}{\tau} \beta_k^2 + \tfrac{4
\NumPlayers^2 \Amax \opnorm{\BRMatrix}^2}{\delta} \alpha_k^2 .
\end{multline}
for any choice of scalar $\constYoung \in (0, \frac{1}{2})$.
\end{lemma}
\noindent See Section~\ref{SecLemCoupledDrift} for the proof.  \\

A few remarks are in order: note that
inequality~\eqref{EqnTotLyapDrift} contains a term that grows as
$(1/\delta)$ in terms of the distance to the boundary. This term
arises because the variance of the TD updates explodes at the
boundary.  For this reason, we need to track the distance to the
boundary as a function of the iteration number.  Note that the bounds
also involve a parameter $\constYoung \in (0, 1/2)$, and we use the
freedom in choosing this quantity in a later part of the argument. \\

Our next step is to use~\Cref{LemCoupledDrift} to derive a recursive
bound on the expected value $\E \TotLyapVal{k}$. To control the error terms on the right-hand side of equation~\eqref{EqnTotLyapDrift}, our analysis involves the
temperature and stepsize parameter specifications
\begin{subequations}
	\begin{align}
		\label{EqnKeyTauBeta}
		\tau(\constYoung) \defn \frac{(1-2\constYoung) \epsilon}{6 \NumPlayers \log \Amax
		}, \quad \mbox{and} \quad \beta(\epsilon,\constYoung,\delta) & \defn
		\frac{(1-2\constYoung) \timescalesep(\constYoung,\tau(\constYoung))^2
			\delta \epsilon}{15 \Amax^2 \NumPlayers^2 \opnorm{\BRMatrix}^2}.
	\end{align}
The right-hand side of the bound~\eqref{EqnTotLyapDrift} contains contraction factors in
front of $\E \StratLyapVal{k}$ and $\E \TDLyapVal{k}$; we would like
to relate these contraction factors so as to form a single term
involving the sum $\E \TotLyapVal{k} = \E \StratLyapVal{k} + \E
\TDLyapVal{k}$.  We do so via a careful choice of the timescale
separation constant that relates the two stepsizes via
$\frac{\beta_k}{\alpha_k} = \timescalesep$:
\begin{align}
  \label{EqnDefnTimescaleSep}
  \timescalesep(\constYoung,\tau) \defn
  \frac{\constYoung(1-\constYoung) \: \tau^3}{4 \opnorm{\BRMatrix}^2}.
\end{align}
\end{subequations}
This choice enables us to match the contraction terms for
$\E\StratLyapVal{k}$ and $\E\TDLyapVal{k}$. While the
specification for $\tau$ is only relevant later while characterizing the iteration complexity, the choice for
$\beta(\epsilon,\constYoung,\delta)$ ensures that $\alpha_k =
\frac{\beta_k}{\timescalesep(\constYoung,\tau)}$ satisfies the upper
bound $\alpha_k \leq \frac{\timescalesep\delta}{\Amax
	\NumPlayers} \leq \frac{\constYoung\delta}{\Amax
	\NumPlayers}$, using the fact that $\timescalesep(\constYoung,\tau) \leq \constYoung$. Then
\begin{align*}
	1 - \alpha_k + \alpha_k^2 \tfrac{\Amax \NumPlayers }{\delta} & \leq 1
	- \alpha_k (1 - \constYoung).
\end{align*}
The choice of the timescale separation constant also ensures that
\begin{align*}
	\frac{\timescalesep(\constYoung,\tau)}{1-\constYoung} =
	\frac{1}{\frac{4 \opnorm{\BRMatrix}^2}{\constYoung \tau^3}} \leq
	\frac{1}{ 1-2\constYoung + \frac{3 \opnorm{\BRMatrix}^2}{\constYoung
			\tau^3}},
\end{align*}
and hence
\begin{align*}
	- \alpha_k(1-\constYoung)+ \beta_k\tfrac{3\opnorm{\BRMatrix}^2}{
		\constYoung \tau^3} = \beta_k \big
	(-\frac{1-\constYoung}{\timescalesep(\constYoung,\tau)} + \frac{
		3\opnorm{\BRMatrix}^2 }{\constYoung \: \tau^3} \big ) \leq -\beta_k
	(1-2 \constYoung),
\end{align*}
enabling us to match the contraction terms for $\StratLyapVal{k}$ and
$\TDLyapVal{k}$.

From our choice~\eqref{EqnDefnTimescaleSep}, we are also guaranteed to have
\begin{align*}
 \tfrac{8 \NumPlayers\opnorm{\BRMatrix}^2}{\tau} \beta_k^2 \leq
 \tfrac{\NumPlayers^2 \Amax \opnorm{\BRMatrix}^2}{\delta}
 \alpha_k^2.
\end{align*}
Combining these relations with the bound~\eqref{EqnTotLyapDrift} yields
\begin{align}
\label{EqnTotalDrift}
\E\TotLyapVal{k+1} \leq \big ( 1 - \beta_k  ( 1- 2\constYoung ) \big ) \E
\TotLyapVal{k} + 2\NumPlayers \tau \beta_k \log \Amax + \tfrac{5 \NumPlayers^2 \Amax \opnorm{\BRMatrix}^2}{\delta} \alpha_k^2.
\end{align}
Note that this drift inequality is valid as long as the iterates are
$\delta$-good (cf. equation~\eqref{EqnDeltaGood}).

Given the drift inequality \eqref{EqnTotalDrift}, our next step is to
bound the iteration complexity.  Concretely, for each $\epsilon > 0$,
we let $K(\epsilon,\delta,\constYoung)$ denote the number of
$\delta$-good iterations of~\Cref{AlgPartial} that are required to
ensure that $\E \big
[\nashg(\BStrat_{K(\epsilon,\delta,\constYoung)+1}\big] \leq
\epsilon$.  The following result, proved in
Section~\ref{SecLemLowerBound}, gives upper bounds on this iteration
complexity:
\begin{lemma}
\label{LemNGBoundPartial}
For any $\constYoung \in (0, \frac{1}{2})$,
initialize \cref{AlgPartial} with the timescale separation constant
and the temperature and stepsize parameters from equations
~\eqref{EqnKeyTauBeta} and~\eqref{EqnDefnTimescaleSep}.  Then the
iteration complexity $K(\epsilon,\delta,\constYoung)$ is bounded as
follows:
  \begin{enumerate}[label=(\alph*),leftmargin=1.8em]
  \item For the constant stepsize $\beta_k \equiv \beta(\epsilon,
    \constYoung,\delta)$, we have
\begin{subequations}
  \begin{align}
\label{EqnKrequiredConstant}
K(\epsilon,\delta,\constYoung) & \leq \big \lceil \frac{\log \big (
  \frac{\epsilon}{3 \TotLyapVal{1} }\big )}{\log\big( 1
  -\beta(\epsilon, \constYoung,\delta) \big (1-2\constYoung \big )
  \big )} \big \rceil .
 \end{align}
\item For the inverse polynomial stepsize $\beta_k = \frac{
  \beta(\epsilon,\constYoung,\delta)}{(k+\koff)^\myexp}$ for some
  exponent $\myexp \in (0,1)$, and the offset \mbox{$\koff = \big
    \lceil (\frac{2\myexp}{\beta(\epsilon,
      \constYoung,\delta)})^{1/(1-\myexp)} \big \rceil$,} we have
  \begin{align}
\label{EqnKrequiredPoly}
K(\epsilon,\delta,\constYoung) & \leq \big \lceil \big (
\frac{1-\myexp}{(1 - 2\constYoung)\beta(\epsilon, \constYoung,\delta)}
\log \frac{ 3 \TotLyapVal{1}}{\epsilon} + (1+\koff)^{(1-\myexp)} \big
)^{1/(1-\myexp)} -k_0 -1 \big \rceil .
  \end{align}
\end{subequations}
\end{enumerate}
\end{lemma}


\subsubsection{Controlling the distance to the boundary}

The delicacy with the iteration complexity in~\Cref{LemNGBoundPartial}
is that it is valid for the actual iterates of our process \emph{only
when}\footnote{To be explicit, the proof of ~\Cref{LemNGBoundPartial}
exploits the drift inequality~\eqref{EqnTotalDrift}, which is only
valid for our process when the iterates are $\delta$-good.}  they are
$\delta$-good at least up to iteration
$K(\epsilon,\delta,\constYoung)$.  Moreover, note that the stepsize
parameter $\beta$ from equation~\eqref{EqnKeyTauBeta} depends on the
target accuracy $\epsilon$, as well as the distance $\delta$ to the
boundary. As a degree of freedom, our results so far allow us to
choose $\delta \in (0,1)$ and $\constYoung \in (0, 1/2)$, and we do so
carefully in the following analysis.

Define $\Kgood(\beta,\delta)$ to be the maximum number of iterations of \cref{AlgPartial} with the stepsize parameter $\beta \in (0,1)$ that can
  be taken while ensuring that the iterates remain
  $\delta$-good (cf. equation~\eqref{EqnDeltaGood}), when the initial mixed strategies $\BStrat_1$ are uniform.

 For a given target accuracy $\epsilon \in (0, 1)$ and
convergence exponent $\offset > 0$, our proof strategy consists of the
following steps:
\begin{enumerate}[label=(\roman*),leftmargin=1.8em]
\item We first establish a lower bound on $\Kgood(\beta,\delta)$ for
  any $\beta,\delta \in (0,1)$. We apply this lower bound for the
  $\beta(\epsilon,\constYoung,\delta)$ specified
  by~\Cref{LemNGBoundPartial}.
\item Our next step is to make careful choice of $
  \delta(\epsilon,\offset)$ and $ \constYoung(\epsilon, \offset)$ as a
  function of $\epsilon$ and $\offset$---in particular,
  see~\Cref{LemUpperBound} to follow---such that the lower bound on the number of
  ``good'' iterations $\Kgood \big
  (\beta(\epsilon,\constYoung(\epsilon,
  \offset),\delta(\epsilon,\offset)),\delta(\epsilon,\offset) \big )$
  is larger than the upper bound on the iteration complexity $K \big (\epsilon,
  \delta(\epsilon, \offset),\constYoung(\epsilon,\offset) \big )$ shown in ~\Cref{LemNGBoundPartial}.  This inequality ensures
  that the iterates $\{\BStrat_k\}_{k \geq 1}$ remain
  $\delta(\epsilon,\offset)$-good until time $K \big (\epsilon,
  \delta(\epsilon, \offset),\constYoung(\epsilon,\offset) \big )$.
\item In the final step, we can then apply~\Cref{LemNGBoundPartial} to
argue that
\begin{align*}
  \E [\nashg(\BStrat_{K \big (\epsilon, \delta(\epsilon,
      \offset),\constYoung(\epsilon,\offset) \big )+1}] \leq
  \epsilon.
\end{align*}
\end{enumerate}

We now state the lemma that provides a lower bound on the number of
``good'' iterations (cf. step (iii) above).  In proving this result,
we choose the free parameters $\constYoung$ and $\delta$ as a function
of the exponent $\offset$ that appears in our final guarantee; see the
proof in~\Cref{SecLemUpperBoundSimple} for details.
\begin{lemma}
\label{LemUpperBound}
For any $\offset > 0$, define
\begin{align}
\delta(\epsilon,\offset) \defn \begin{cases} \big(
  \frac{\epsilon}{3\TotLyapVal{1}}\big )^{1+\offset} \frac{1}{\Amax} &
  \mbox{for the constant stepsize case, and} \\
  \big( \frac{\epsilon}{3 \TotLyapVal{1} }\big )^{1+\offset}
  \frac{\exp(1)}{\Amax} & \mbox{for the inverse polynomial case.}
    \end{cases}
  \end{align}
There exists a scalar \mbox{$\constYoung(\epsilon,\offset) \in
  (0,\frac{1}{2})$} such that the maximum number of ``good" iterations
$\Kgood \big (\beta(\epsilon,\constYoung(\epsilon, \offset),
\delta(\epsilon,\offset)), \delta(\epsilon,\offset) \big )$ is lower
bounded by the upper bounds in~\Cref{LemNGBoundPartial} for both
stepsizes.
\end{lemma}

\subsubsection{Proof of~\Cref{LemUpperBound}}
\label{SecLemUpperBoundSimple}
We give the proof for the constant stepsize case here; the proofs for
the other cases are analogous, and given in the supplementary
Section~\ref{SecLemUpperBound}.  In order to
prove~\Cref{LemUpperBound}, it suffices to show that for any
$\uplowratio(\offset) \in (1, 1 + \offset)$, there exists a choice
\mbox{$\constYoung(\offset) \in (0,\frac{1}{2})$} such that
\begin{subequations}
\begin{align}
  \label{EqnUpperSimple}
\Kgood \big (\beta(\epsilon, \constYoung(\epsilon,
\offset),\delta(\epsilon, \offset)), \delta(\epsilon,\offset) \big )
\geq \frac{\uplowratio(\offset) \log \big ( \frac{\epsilon}{3
    \TotLyapVal{1} }\big )}{\log\big(1- (1 - 2
  \constYoung(\offset))\beta(\epsilon, \constYoung(\offset),
  \delta(\epsilon,\offset))\big )},
\end{align}
where $\delta(\epsilon,\offset) \defn
\big(\frac{\epsilon}{3\TotLyapVal{1}}\big )^{1+\offset}
\frac{1}{\Amax}$.  Here the quantity $\uplowratio(\offset)>1$ is the
factor by which $\Kgood(\beta,\delta)$ is greater than
$K(\epsilon,\delta,\constYoung)$ (cf. \cref{LemNGBoundPartial}), for
the given choices of $\delta$ and $\constYoung$.

In order to prove the bound~\eqref{EqnUpperSimple}, we first show that
the iterates remain $\delta$-good for all iterates $K$ such that
\begin{align}
  \label{EqnKbound}
  K \leq \frac{\log(\Amax \delta)}{\log(1-\beta)}.
\end{align}
\end{subequations}
Indeed, from the strategy update dynamics~\eqref{EqnPartialStrategy},
we have the elementwise inequality $\pii_{k+1} \succeq \pii_k
(1-\beta_k)$.  Iterating this inequality for $k = 1, \ldots, K + 1$
yields \mbox{$\strat[i]_{K+1} \succeq \frac{1}{\Amax} \prod_{j=1}^{K}
  (1-\beta_j)$,} where we have made use of the lower bound $\BStrat_1
\succeq \frac{1}{\Amax} \boldsymbol{e}$, as guaranteed by our uniform
initialization.  Therefore, in order to ensure that $\strat[i]_k
\succeq \delta$ for all $k = 1, \ldots, K + 1$, it suffices to have
$\frac{1}{\Amax} \prod_{k=1}^{K} (1-\beta_k)  \geq \delta$, using the
fact that $\beta_k \in (0,1)$ for all $k$.  For the constant stepsize
$\beta_k \equiv \beta \in (0,1)$, this condition translates to the
bound $(1-\beta)^{K} \geq \Amax \delta$, \mbox{or equivalently} $K
\leq \frac{\log (\Amax \delta ) }{\log (1-\beta)}$. \\

We now use the bound~\eqref{EqnKbound} to prove the
claim~\eqref{EqnUpperSimple}.  Note that the latter bound is a lower
bound on $\Kgood(\beta,\delta)$. Now we use the specified choices of
$\delta$ and $\constYoung$ as a function of $\epsilon$ and $\offset$
to ensure that these lower bounds on $\Kgood(\beta,\delta) \geq
K(\epsilon, \delta, \constYoung)$. For any $\uplowratio \in
(1,1+\offset)$, by choosing $\constYoung(\offset) \in (0,\frac{1}{2})$
small enough, we can ensure that
\begin{align*}
  \frac{\log(1-\beta)}{\log \big ( 1 - \beta \big ( 1- 2
    \constYoung(\offset) \big ) \big ) } \leq \frac{1 +
    \offset}{\uplowratio}.
\end{align*}
Applying this bound yields
\begin{align*}
  \frac{\log(\Amax \ \delta(\epsilon,\offset))}{\log(1 - \beta)} & = (1
  + \offset) \frac{\log \big ( \frac{\epsilon}{3 \TotLyapVal{1}}\big
    )}{\log(1-\beta)} \geq \frac{\uplowratio \log \big (
    \frac{\epsilon}{3 \TotLyapVal{1}}\big )}{\log \big ( 1 - \beta \big
    ( 1- 2 \constYoung(\offset) \big ) \big )}.
\end{align*}
Therefore, the iterates $\{\BStrat_k\}_{k \geq 1}$ remain
$\delta(\epsilon,\offset)$-good up until the claimed
time~\eqref{EqnUpperSimple}.

\subsubsection{Completing the proof}
In order to complete the proof of~\Cref{ThmFullInfoPoly}, it suffices
to combine the previous pieces.  Summarizing, we have shown that:
\begin{itemize}[leftmargin=1.8em]
\item For any pair of scalars $\delta \in (0,1)$ and $\constYoung \in
  (0,\frac{1}{2})$,~\Cref{LemNGBoundPartial} specifies the number of
  iterations $K(\epsilon, \delta,\constYoung)$ sufficient to ensure
  that $\E [\nashg(\BStrat_{K(\epsilon,\delta,\constYoung)+1})] \leq
  \epsilon$.  (This guarantee is predicated upon the iterates
  $\{\BStrat_k\}_{k \geq 1}$ of \cref{AlgPartial} being $\delta$-good
  until time $K(\epsilon, \delta,\constYoung)$, and
  that~\Cref{AlgPartial} is initialized according to
  equations~\cref{EqnDefnTimescaleSep,EqnKeyTauBeta}).
\item For a given target Nash gap $\epsilon \in (0,1)$ and $\offset >
  0$, \Cref{LemUpperBound} specifies the pair $(\delta, \constYoung)$
  as a function of the pair $(\epsilon, \offset)$, so that we may
  write $\delta(\epsilon, \offset)$ and $\constYoung(\epsilon,
  \offset)$ to indicate this dependence.  Note
  that~\Cref{LemNGBoundPartial} applies to these choices of
  $\delta(\epsilon, \offset)$ and $\constYoung(\epsilon, \offset)$ as
  well.
\item For these choices of $(\delta, \constYoung)$,
  \Cref{LemUpperBound} shows that the iterates $\{\BStrat_k\}_{k \geq
    1}$ of~\Cref{AlgPartial} are $\delta(\epsilon,\offset)$-good for
  at least $K^\star(\epsilon,\offset)$ iterations, where
  $K^\star(\epsilon,\offset)$ is the upper bound on $K \big
  (\epsilon,\delta(\epsilon,\offset),\constYoung(\epsilon,\offset)
  \big )$ from~\Cref{LemNGBoundPartial}.
\end{itemize}

It follows from the results of~\Cref{LemNGBoundPartial,LemUpperBound}
that the Nash gap evaluated at time $K^\star(\epsilon,\offset)+1$ is
\mbox{at most $\epsilon$.} The upper bound on $K(\epsilon)$ given
in~\cref{ThmMinInfoPoly} follows by bounding
$K^\star(\epsilon,\offset)$ from above.

\myunderparagraph{Remarks.} In~\Cref{LemNGBoundPartial}, the stepsize
parameter is chosen as a function of $\delta$.  It could be desirable
to obtain an anytime version of~\Cref{ThmMinInfoPoly}, for which the
Nash gap continues to be bounded by $\epsilon$ for all time steps
greater than $K^\star(\epsilon,\offset)$.  Such a guarantee can be
achieved by running $K^\star(\epsilon,\offset)$ iterations, and then
annealing the stepsize parameter as $\delta$ goes to zero, to ensure
that the iterates continue to be $\delta$-good
and~\Cref{LemNGBoundPartial} holds.



\section{Discussion}

In this paper, we studied best-response type learning dynamics that
arise from the discretization of continuous-time best-response
dynamics applied to zero-sum polymatrix games. We analyzed these
dynamics in both the full information setting as well as the more
challenging setting of minimal information, in which each player
observes only their random payoff.  In the latter context, we provided
the first polynomial-scaling finite-sample guarantees for these
best-response type dynamics in the \infosetuplower setting without
additional exploration. Our results also exhibited an interesting
dependence on the underlying polymatrix game through the parameter
$\opnorm{\BRMatrix}$. The analysis involved some new ideas, including
careful tracking of variance of the stochastic updates as they
approach the boundary of the probability simplex.

In terms of open questions, we suspect that our current results are
not sharp in terms of dependence on $\epsilon$ and $\Amax$ for
last-iterate convergence to a Nash equilibrium.  For instance, in the
simpler setting of full information, the optimal rate is known to
scale as \mbox{$\mathcal{O}(\frac{1}{\epsilon})$}; we are not yet
aware if our dynamics can achieve this rate, or if the
\mbox{$\mathcal{O}(\frac{1}{\epsilon^2})$} guarantees that we have
provided are, in fact, unimprovable.  In terms of the dependence on
the underlying polymatrix game, it would be interesting to see if a
finer dependence on the underlying pairwise games can be
elicited. Additionally, our analysis in this paper was based upon a
fixed temperature parameter $\tau$. Studying dynamics with a
time-varying temperature $\{\tau_k\}_{k \geq 1}$ might help eliminate
the constant smoothing bias in the current bounds on the Nash
gap.

\noindent \section*{Acknowledgments}

MJW was partially supported by ONR grant N00014-21-1-2842 and NSF
DMS-2311072, AO was supported by ONR grant N00014-25-1-2296 and FZF was supported by the MIT EECS Thriving Stars Energy Fellowship.

\printbibliography

\appendix

\section{The full information case}
In this section, we first complete the proof of~\Cref{ThmFullInfoPoly} and~\Cref{CorFullInfoPoly} for the inverse linear and inverse polynomial stepsize schedules. We then provide the proof of~\Cref{LemVprops} in~\Cref{SecLemVprops}.


\subsection{The two stepsizes}
\label{AppFullStep}

In the main text, we used the drift
inequality~\eqref{EqnDriftStratLyap} to prove the bound
in~\Cref{ThmFullInfoPoly} for constant stepsizes.  Here we cover the
two remaining cases.

\subsubsection{Inverse linear stepsizes}

Consider the inverse linear stepsizes $\beta_k = \frac{\beta}{k}$ for
some $\beta > 0$. Solving the inequality~\eqref{EqnDriftStratLyap}
gives that $\StratLyapVal{k} \defn \stratLyap\big( \BStrat_k \big )$
satisfies the equation
\begin{multline}
	\label{EqnBetaKlinear}
	\StratLyapVal{K+1} \leq \prod_{k=1}^{K} \big( 1- \beta_k\big)
	\StratLyapVal{1} +\frac{\NumPlayers\opnorm{\BRMatrix}^2}{\tau}
	\sum_{k=1}^{K} \beta_k^2 \prod_{j=k+1}^{K} (1-\beta_j) \\+ \NumPlayers
	\tau \log\Amax \sum_{k=1}^{K} \beta_k \prod_{j=k+1}^{K} (1-\beta_j).
\end{multline}

\myunderparagraph{Bounding the first term:}
\begin{align*}
	\prod_{k=1}^{K} (1-\beta_k) \; = \; \exp \big( \sum_{k=1}^{K}
	\log(1-\beta_k) \big) & \overset{(a)}{\leq} \exp \big(- \sum_{k=1}^{K}
	\beta_k \big) \\
	& \overset{(b)}{\leq} \exp \big( - \beta \int_{1}^{K+1} \frac{1}{x}
	\ dx \big) \\
	& = \big ( \frac{1}{K+1} \big )^{\beta},
\end{align*}
where in step (a), we used the bound $\log(1-x) \leq -x$ for $x \in
(0,1)$ and in step (b), we bounded the Riemann sum.

\myunderparagraph{Bounding the second term:} Using a series of
arguments similar to the upper bound for the first term yields
\begin{align*}
	\prod_{j=k+1}^{K} (1-\beta_j) &= \exp \big( \sum_{j=k+1}^{K}
	\log(1-\beta_j) \big)\\ & \leq \exp \big(- \sum_{j=k+1}^{K} \beta_j
	\big) \\
	& \leq \exp \big( - \beta \int_{k+1}^{K+1} \frac{1}{x} \ dx \big) \\
	& = \big ( \frac{k+1}{K+1} \big )^{\beta}. \numberthis\label{EqnFullInfoStepsizeCalc}
\end{align*}
Substituting this bound into the expression for the second term yields
\begin{align*}
	\frac{\NumPlayers\opnorm{\BRMatrix}^2}{\tau} \sum_{k=1}^{K}
	\beta_k^2 \prod_{j=k+1}^{K} (1-\beta_j) &\leq
	\frac{\NumPlayers\opnorm{\BRMatrix}^2}{\tau}\sum_{k=1}^{K} \beta_k^2
	\big ( \frac{k+1}{K+1} \big )^{\beta} \\ &\leq
	\frac{4\NumPlayers\opnorm{\BRMatrix}^2 \beta^2}{\tau (K+1)^{\beta}}
	\sum_{k=1}^{K} \frac{1}{(k+1)^{2-\beta}}.
\end{align*}

We now consider five cases.

\subsubsection*{Case 1} If $\beta=1$, then we have
\begin{align*}
	\sum_{k=1}^{K} \frac{1}{(k+1)^{2-\beta}} &\leq \int_{1}^{K+1}
	\frac{1}{x} \ dx = \log (K+1).
\end{align*}
Therefore,
\begin{align*}
	\frac{4 \NumPlayers \opnorm{\BRMatrix}^2 \beta^2}{\tau
		(K+1)^{\beta}} \sum_{k=1}^{K} \frac{1}{(k+1)^{2-\beta}}
	&\leq \frac{4 \NumPlayers \opnorm{\BRMatrix}^2 \log
		(K+1)}{\tau (K+1)} .
\end{align*}

\subsubsection*{Case 2}
If $\beta > 2$, then we can write
\begin{align*}
	\sum_{k=1}^{K} \frac{1}{(k+1)^{2-\beta}} &= \sum_{k=2}^{K+1}
	k^{\beta-2} \leq \int_{2}^{K+2} x^{\beta-2} dx \leq
	\frac{(K+2)^{\beta-1}}{\beta-1}.
\end{align*}
Therefore,
\begin{align*}
	\frac{4 \NumPlayers \opnorm{\BRMatrix}^2 \beta^2}{\tau
		(K+1)^{\beta}} \sum_{k=1}^{K} \frac{1}{(k+1)^{2-\beta}} &\leq
	\frac{4 \NumPlayers \opnorm{\BRMatrix}^2\beta^2
		(K+2)^{\beta-1}}{\tau(\beta-1) (K+1)^{\beta}} \\
	& \leq \frac{4 \NumPlayers \opnorm{\BRMatrix}^2\beta^2 2^{\beta-1}
	}{\tau(\beta-1)(K+1)} .
\end{align*}

\subsubsection*{Case 3} If $\beta = 2$, then
\begin{align*}
	\sum_{k=1}^{K} \frac{1}{(k+1)^{2-\beta}} = K.
\end{align*}
Therefore,
\begin{align*}
	\frac{4 \NumPlayers \opnorm{\BRMatrix}^2 \beta^2}{\tau
		(K+1)^{\beta}} \sum_{k=1}^{K} \frac{1}{(k+1)^{2-\beta}}
	&\leq \frac{4 \NumPlayers \opnorm{\BRMatrix}^2 \beta^2 K}{\tau
		(K+1)^{\beta}} \\ &\leq \frac{16 \NumPlayers
		\opnorm{\BRMatrix}^2 }{\tau (K+1)}.
\end{align*}

\subsubsection*{Case 4} If $\beta \in (1,2)$, then
\begin{align*}
	\sum_{k=1}^{K} \frac{1}{(k+1)^{2-\beta}} = \sum_{k=2}^{K+1}
	\frac{1}{k^{2-\beta}} \leq \int_{1}^{K+1} \frac{1}{x^{2-\beta}}
	dx \leq \frac{K^{\beta-1}}{{\beta-1}}.
\end{align*}
Therefore, we have
\begin{align*}
	\frac{4 \NumPlayers \opnorm{\BRMatrix}^2 \beta^2}{\tau
		(K+1)^{\beta}} \sum_{k=1}^{K} \frac{1}{(k+1)^{2-\beta}} &\leq
	\frac{4 \NumPlayers \opnorm{\BRMatrix}^2 \beta^2
		K^{\beta-1}}{\tau(\beta-1) (K+1)^{\beta}} \\
	& \leq \frac{4 \NumPlayers \opnorm{\BRMatrix}^2 \beta^2
	}{\tau(\beta-1)K}.
\end{align*}

\subsubsection*{Case 5} Finally, if $\beta \in (0,1)$, then we have
\begin{align*}
	\sum_{k=1}^{K} \frac{1}{(k+1)^{2-\beta}} = \sum_{k=2}^{K+1}
	\frac{1}{k^{2-\beta}} \leq \int_{1}^{K+1}
	\frac{1}{x^{2-\beta}} dx \leq \frac{1}{{1-\beta}}.
\end{align*}
Therefore,
\begin{align*}
	\frac{4\NumPlayers\opnorm{\BRMatrix}^2 \beta^2}{\tau
		(K+1)^{\beta}} \sum_{k=1}^{K} \frac{1}{(k+1)^{2-\beta}} &\leq \frac{4\NumPlayers\opnorm{\BRMatrix}^2 \beta^2}{\tau(1-\beta)
		(K+1)^{\beta}}.
\end{align*}

\myunderparagraph{Bounding the third term:} Using
equation~\eqref{EqnFullInfoStepsizeCalc} and a series of arguments
similar to the upper bound for the second term:
\begin{align*}
	\NumPlayers\tau \log \Amax \sum_{k=1}^{K} \beta_k \prod_{j=k+1}^{K}
	(1-\beta_j) &\leq \NumPlayers\tau \log \Amax \sum_{k=1}^{K} \beta_k
	\big ( \frac{k+1}{K+1} \big )^{\beta} \\ &\leq \frac{2\NumPlayers\tau
		\log \Amax \beta}{(K+1)^{\beta}} \sum_{k=1}^{K}
	\frac{1}{(k+1)^{1-\beta}}.
\end{align*}

\subsubsection*{Case 1}  If $\beta = 1$, then we have
\begin{align*}
	\sum_{k=1}^{K} \frac{1}{(k+1)^{1-\beta}} & = K.
\end{align*}
Therefore,
\begin{align*}
	\frac{2\NumPlayers\tau \log \Amax \beta}{(K+1)^{\beta}}
	\sum_{k=1}^{K} \frac{1}{(k+1)^{1-\beta}} &\leq {2\NumPlayers\tau
		\log \Amax}.
\end{align*}

\subsubsection*{Case 2}

If $\beta > 1$, then we can write
\begin{align*}
	\sum_{k=1}^{K} \frac{1}{(k+1)^{1-\beta}} &= \sum_{k=2}^{K+1}
	{k^{\beta-1}} \leq \int_{2}^{K+2} x^{\beta-1} dx \leq
	\frac{(K+2)^{\beta}}{\beta}.
\end{align*}
Therefore,
\begin{align*}
	\frac{2\NumPlayers\tau \log \Amax \beta}{(K+1)^{\beta}}
	\sum_{k=1}^{K} \frac{1}{(k+1)^{1-\beta}} &= \frac{4\tau \log \Amax
		(K+2)^\beta}{(K+1)^{\beta}} \leq {2\NumPlayers\tau 2^\beta \log
		\Amax}.
\end{align*}

\subsubsection*{Case 3}
If $\beta < 1$, then
\begin{align*}
	\sum_{k=1}^{K} \frac{1}{(k+1)^{1-\beta}} &= \sum_{k=2}^{K+1}
	\frac{1}{k^{1-\beta}} \; \leq \; \int_{1}^{K+1} \frac{1}{x^{1-\beta}} dx
	\leq \frac{(K+1)^{\beta}}{\beta }.
\end{align*}
Therefore,
\begin{align*}
	\frac{2\NumPlayers\tau \log \Amax
		\beta}{(K+1)^{\beta}} \sum_{k=1}^{K} \frac{1}{(k+1)^{1-\beta}}  &\leq {2\NumPlayers\tau \log \Amax}.
\end{align*}
The statement of \cref{ThmFullInfoPoly} for inverse linear stepsizes
is for the specific case when $\beta \in (1,2]$. To bound the
iteration complexity, we choose $\tau = \epsilon/(24 \NumPlayers \log
\Amax)$.


\subsubsection{Inverse polynomial stepsize}

Now suppose that $\beta_K = \frac{\beta}{(K+\koff)^\myexp}, \myexp \in
(0,1)$ and $\beta \in (0,1)$. Recall that $\stratLyapVal{K+1}$
satisfies the inequality
\begin{multline*}
	\StratLyapVal{K+1} \leq \prod_{k=1}^{K} \big( 1- \beta_k\big)
	\StratLyapVal{1} +\frac{\NumPlayers\opnorm{\BRMatrix}^2}{\tau}
	\sum_{k=1}^{K} \beta_k^2 \prod_{j=k+1}^{K} (1-\beta_j) \\+
	\NumPlayers \tau \log\Amax \sum_{k=1}^{K} \beta_k \prod_{j=k+1}^{K}
	(1-\beta_j).
\end{multline*}

\myunderparagraph{Bounding the first term:} We have
\begin{align*}
	\prod_{k=1}^{K} (1-\beta_k) &\leq \exp \big( - \beta
	\int_{1}^{K+1} \frac{1}{(x+\koff)^\myexp} \ dx \big) \\ &=
	\exp \big( - \frac{\beta}{1-\myexp} \big(
	(K+\koff+1)^{1-\myexp} - (1+\koff)^{1-\myexp} \big)\big ).
\end{align*}

\myunderparagraph{Bounding the second term:} We define the sequence
$\{u_k\}_{k \geq 1}$ via the recursion
\begin{align*}
	u_{k+1} = \big( 1 - \beta_k \big) u_k + \beta_k^2 \qquad \mbox{with
		initialization $u_1 = 0$.}
\end{align*}
Unwrapping this recursion, we find that $u_K = \sum_{k=1}^{K}
\beta_k^2 \prod_{j=k+1}^{K} (1-\beta_j)$.  It can be shown by
induction~(see p.36, \cite{chen2023lyapunov}) that \mbox{$u_K \leq
	2\beta_K$.} Therefore,
\begin{align*}
	\frac{\NumPlayers\opnorm{\BRMatrix}^2}{\tau} \sum_{k=1}^{K} \beta_k^2
	\prod_{j=k+1}^{K} (1-\beta_j) = \frac{\NumPlayers\opnorm{\BRMatrix}^2
		u_K }{\tau} \leq \frac{2 \NumPlayers\opnorm{\BRMatrix}^2
		\beta_K}{\tau} = \frac{2\NumPlayers\opnorm{\BRMatrix}^2\beta}{\tau
		(K+\koff)^\myexp}.
\end{align*}


\myunderparagraph{Bounding the third term:} We define the sequence
$\{u_k\}_{k \geq 1}$ via the recursion
\begin{align*}
	u_{k+1} = \big( 1 - \beta_k \big) u_k + \beta_k \quad \mbox{with
		initial value $u_1 = 0$.}
\end{align*}
Note that since $\beta_k \in (0,1)$, it follows that $u_k \leq 1$ for
all $k = 1, 2, \ldots$.  Expanding out the recursion yields $u_K =
\sum_{k=1}^{K} \beta_k \prod_{j=k+1}^{K} (1-\beta_j)$.  Combining with
the inequality $u_K \leq 1$, we find that
\begin{multline*}
	\stratLyapVal{K+1} \leq \exp \big( - \frac{\beta}{1-\myexp} \big( (K
	+ \koff + 1)^{1-\myexp} - (1 + \koff)^{1-\myexp} \big)\big )
	\stratLyapVal{1} + \frac{2 \NumPlayers\opnorm{\BRMatrix}^2
		\beta}{\tau(K+\koff)^\myexp} \\
	+ \NumPlayers \tau \log \Amax.
\end{multline*}

The bound on the iteration complexity in~\Cref{CorFullInfoPoly} then follows by setting $\tau =
\epsilon/(3 \NumPlayers \log \Amax)$.


\subsection{Proof of~\Cref{LemVprops}}
\label{SecLemVprops}

We split our proof into parts, corresponding to the two statements in
the lemma.  In both cases, we make use of Danskin's theorem~\cite{danskin2012theory} to compute the gradient
\begin{align}
	\label{EqnDanskinGradient}
	\nabla_{\strat[i]} \stratLyap(\BStrat) &= \sum_{j \in \nodes}
	{\RMatrix{j}{i}}^{\intercal}
	\sigma_{\tau}(\AvgPayoff{j}(\BStrat)).
\end{align}

\paragraph{Proof of part (a)}

In order to compute the Hessian, we take the derivative of the
gradient from equation~\eqref{EqnDanskinGradient}.  Doing so via chain
rule yields
\begin{align*}
	\nabla^2_{\strat[i], \strat[m]} \stratLyap(\BStrat)
	&=\frac{1}{\tau} \sum_{j \in
		\nodes}{\RMatrix{j}{i}}^{\intercal} \Sigma^j_\tau(\BStrat)
	{\RMatrix{j}{m}},
\end{align*}
where $\Sigma^j(\BStrat) \defn \text{diag} \big(
\sigma_{\tau}(\AvgPayoff{j}(\BStrat) \big) - \sigma_{\tau} \big(
\AvgPayoff{j}(\BStrat) \big) \sigma_{\tau} \big(
\AvgPayoff{j}(\BStrat) \big)^{\intercal}$ for each $j \in \nodes$. The
Hessian of $\stratLyap$ can now be decomposed in the following way:
\begin{align*}
	\nabla^2 \stratLyap(\BStrat) &= \frac{1}{\tau} \BRMatrix^\intercal
	\boldsymbol{\Sigma}(\BStrat) \BRMatrix,
\end{align*}
where $\boldsymbol{\Sigma}(\BStrat)$ is a block diagonal matrix with
the $i^{\text{th}}$ diagonal block $\Sigma^i(\BStrat)$ for $i \in
\nodes$.  By the sub-multiplicativity of the operator norm, we have
\begin{align*}
	\opnorm{ \nabla^2 \stratLyap(\BStrat)} & \leq \frac{1}{\tau}
	\opnorm{\BRMatrix}^2 \; \opnorm{\boldsymbol{\Sigma}(\BStrat)} .
\end{align*}

Turning to the matrix ${\Sigma^i}(\BStrat)$, we recognize it as the
covariance matrix of a multinomial random vector, so that it must be
positive semi-definite (PSD).  It follows that the matrix
$\boldsymbol{\Sigma}(\BStrat)$ is also PSD, so that its operator norm
can be written as $\opnorm{\boldsymbol{\Sigma}(\BStrat)} = \underset{j
	\in \nodes}{\max} \ \lammax(\Sigma^j(\BStrat))$.  Now observe that
\begin{align*}
	\Sigma^i(\BStrat) \preceq \text{diag} \big(
	\sigma_{\tau}(\AvgPayoff{i}(\strat^{-i})) \big),
\end{align*}
where $\preceq$ denotes the PSD order.  Thus, we have
\begin{align*}
	\lammax(\Sigma^i((\BStrat)) \leq \lammax \Big[ \text{diag}
	\big( \Sigma^i(\AvgPayoff{i}(\strat^{-i})) \big) \Big] \;
	\leq 1.
\end{align*}
Putting together the pieces yields the bound $\opnorm{\nabla^2
	\stratLyap(\BStrat)} \leq \frac{\opnorm{\BRMatrix}^2}{\tau}$, as
claimed.

\paragraph{Proof of part (b)}
Again making use of the gradient
representation~\eqref{EqnDanskinGradient} yields
\begin{multline*}
	\big \langle \nabla_{\strat[i]} \stratLyap(\BStrat),
	\sigma_{\tau}(\AvgPayoff{i}(\strat^{-i})) - \strat[i] \big \rangle =
	\big \langle \sum_{j \in \nodes} {\RMatrix{j}{i}}^\intercal
	\sigma_{\tau}(\AvgPayoff{j}(\BStrat)),
	\sigma_{\tau}(\AvgPayoff{i}(\strat^{-i})) - \strat[i] \big \rangle
	\\
	= \sum_{j \in \nodes}
	\sigma_{\tau}(\AvgPayoff{j}(\BStrat))^{\intercal} {\RMatrix{j}{i}}
	\sigma_{\tau}(\AvgPayoff{j}(\BStrat)) - \sum_{j \in \nodes}
	\sigma_{\tau}(\AvgPayoff{j}(\BStrat))^{\intercal} {\RMatrix{j}{i}}
	\strat[i]
\end{multline*}
for each $i \in \nodes$.  The zero-sum property~\eqref{EqnZeroSum}
ensures that
\begin{align*}
	\sum_{i \in \nodes }\sum_{j \in \nodes}
	\sigma_{\tau}(\AvgPayoff{j}(\BStrat))^{\intercal} {\RMatrix{j}{i}}
	\sigma_{\tau}(\AvgPayoff{i}(\strat^{-i})) = 0.
\end{align*}
Combined with the definition~\eqref{EqnNovelLyapunov} of the Lyapunov
function, we find that
\begin{align*}
	\sum_{i \in \nodes} \big \langle \nabla_{\strat[i]}
	\stratLyap(\BStrat), \sigma_{\tau}(\AvgPayoff{i}(\strat^{-i})) -
	\strat[i] \big \rangle & = - \sum_{i \in \nodes} \sum_{j \in \nodes}
	\sigma_{\tau}(\AvgPayoff{j}(\BStrat))^{\intercal} {\RMatrix{j}{i}}
	\strat[i] \\ &= - \sum_{j \in \nodes}
	\sigma_{\tau}(\AvgPayoff{j}(\BStrat))^{\intercal}
	\AvgPayoff{j}(\BStrat) \\ & = - \sum_{j \in \nodes}
	\sigma_{\tau}(\AvgPayoff{j}(\BStrat))^{\intercal}
	\AvgPayoff{j}(\BStrat) + \tau
	\EntFun(\sigma_{\tau}(\AvgPayoff{j}(\BStrat))) \\ & \hspace{15em}- \tau
	\EntFun(\sigma_{\tau}(\AvgPayoff{j}(\BStrat))) \\
	& \leq -\stratLyap(\BStrat) + \NumPlayers\tau \log \Amax,
\end{align*}
which completes the proof.


\section{The minimal information case}
\label{SecKeyTechnical}

We now collect together the proofs of the key technical lemmas that
were used in the proof of~\Cref{ThmMinInfoPoly}
from~\Cref{SecProofThmMinInfo}.  More specifically, we prove
Lemmas~\ref{LemCoupledDrift}, \ref{LemNGBoundPartial}
and~\ref{LemUpperBound} in Sections~\ref{SecLemCoupledDrift}, \ref{SecLemLowerBound},
and~\ref{SecLemUpperBound}, respectively.

\subsection{Proof of~\Cref{LemCoupledDrift}}
\label{SecLemCoupledDrift}
This section is devoted to the proof of the drift inequality stated
in~\Cref{LemCoupledDrift}.  It suffices to establish the following two
claims:
\begin{subequations}
	\begin{align}
		\E \stratLyapVal{k+1} & \leq \big (1-\beta_k(1-\constYoung) \big )
		\E \stratLyapVal{k} + \tfrac{4 \NumPlayers
			\opnorm{\BRMatrix}^2}{\tau} \beta_k^2 +
		\beta_k\tfrac{\opnorm{\BRMatrix}^2}{2 \constYoung \tau^3} \E
		\TDLyapVal{k} + 2 \NumPlayers \beta_k \tau \log \Amax, \quad
		\numberthis \label{EqnCoupledDriftStrat}
	\end{align}
	as well as
	\begin{multline}
		\label{EqnCoupledDriftQ}
		\E \TDLyapVal{k+1} \leq \big ((1-\alpha_k)^2 + \alpha_k^2 \tfrac{\Amax
			\NumPlayers }{\delta} \big ) \E \TDLyapVal{k} + {4 \NumPlayers
			\beta_k^2 \opnorm{\BRMatrix}^2}+ \tfrac{4 \NumPlayers^2 \Amax
			\opnorm{\BRMatrix}^2}{\delta} \alpha_k^2 \\ +
		(1-\alpha_k)\tfrac{2\beta_k}{ \tau \constYoung} \opnorm{\BRMatrix}^2
		\E \TDLyapVal{k} + (1-\alpha_k)\beta_k\constYoung \E\StratLyapVal{k},
	\end{multline}
\end{subequations}
for any choice of scalar $\constYoung \in (0, \frac{1}{2})$.

Given these claims, inequality~\eqref{EqnTotLyapDrift}
in~\Cref{LemCoupledDrift} follows by summing together the
inequalities~\eqref{EqnCoupledDriftStrat}
and~\eqref{EqnCoupledDriftQ}. The rest of this section is devoted to
proving each of these drift inequalities.

\paragraph{Proof of $\StratLyap$-inequality~\eqref{EqnCoupledDriftStrat}}

Our proof makes use of the previously stated result
from~\Cref{LemVprops} on the properties of the Lyapunov function
$\stratLyap$.  We also require some additional properties, which we
summarize in the following auxiliary result:
\begin{lemma}
	\label{LemVpropsPartial}
	The Lyapunov function $\stratLyap$ has the following properties:
	\begin{enumerate}[label=(\alph*)]
		\item We have the lower bound
		\begin{align*}
			\stratLyap(\BStrat) &\geq \frac{\tau}{2} \big \|
			\sigma_{\tau}(\BAvgPayoff(\BStrat)) - \BStrat \big \|_2^2 \qquad
			\mbox{for any mixed strategy $\BStrat \in \bigdelta$.}
		\end{align*}
		\item For any mixed strategy $\BStrat \in \bigdelta$ and
		$\boldsymbol{u}=(u^i)_{i \in \nodes} \in \otimes_{i \in \nodes}
		\R^{|\actions^i|} $, we have
		\begin{align*}
			\big \langle \nabla_{\BStrat} \stratLyap(\BStrat),
			\sigma_{\tau}(\boldsymbol{u}) -\sigma_{\tau}(\BRMatrix \BStrat)
			\big \rangle &\leq \NumPlayers\tau\log\Amax + \constYoung
			\stratLyap(\BStrat) + \frac{\opnorm{\BRMatrix}^2}{2\constYoung
				\tau^3} \| \BAvgPayoff(\BStrat) - \boldsymbol{u}\|^2_2,
		\end{align*}
		valid for any constant $\constYoung \in (0,1)$.
	\end{enumerate}
\end{lemma}
\noindent See~\cref{AppLemVdriftPartial} for the proof of this
auxiliary claim. \\

From the proof of~\Cref{ThmFullInfoPoly}, recall
inequality~\eqref{EqnFullInfoVDrift}, which ensures that
\begin{align}
	\label{eqn:V_taylor}
	\StratLyapVal{k+1} & \leq \StratLyapVal{k} + \beta_k \big \langle
	\nabla_{\BStrat} \stratLyap(\BStrat_k),
	\sigma_{\tau}(\BAvgPayoff(\BStrat_k)) - \BStrat_k \big \rangle +
	\frac{4 \NumPlayers \opnorm{\BRMatrix}^2}{\tau} \beta_k^2.
\end{align}
Next we modify the first-order terms by adding and subtracting
$\sigma_{\tau}(\BAvgPayoff(\BStrat_k))$ inside the inner product;
doing so yields
\begin{multline}
	\label{EqnVFirstOrderDecomp}
	\big \langle \nabla_{\BStrat} \stratLyap(\BStrat_k),
	\sigma_{\tau}(\BAvgPayoff(\BStrat_k)) - \BStrat_k \big \rangle = \big
	\langle \nabla_{\BStrat} \stratLyap(\BStrat_k),
	\sigma_{\tau}(\BAvgPayoff(\BStrat_k)) - \BStrat_k \big \rangle \\ +
	\big \langle \nabla_{\BStrat} \stratLyap(\BStrat_k),
	\sigma_{\tau}(\BAvgPayoff(\BStrat_k)) -
	\sigma_{\tau}(\BAvgPayoff(\BStrat_k)) \big \rangle.
\end{multline}
Applying~\cref{LemVprops}(b) and \cref{LemVpropsPartial}(b) gives
upper bounds on each of the terms in the RHS of
equation~\eqref{EqnVFirstOrderDecomp}---namely
\begin{subequations}
	\label{EqnVFirstOrderBounds}
	\begin{align}
		\big \langle \nabla_{\BStrat} \stratLyap(\BStrat_k),
		\sigma_{\tau}(\BAvgPayoff(\BStrat_k)) - \BStrat_k \big \rangle
		&\leq - \stratLyap\big( \BStrat_k \big) +\NumPlayers\tau \log
		\Amax, \quad \text{and} \\
		\big \langle \nabla_{\BStrat} \stratLyap(\BStrat_k),
		\sigma_{\tau}(\BAvgPayoff(\BStrat_k)) -
		\sigma_{\tau}(\BAvgPayoff(\BStrat_k)) \big \rangle &\leq \constYoung
		\stratLyap(\BStrat_k) + \NumPlayers \tau \log \Amax +
		\frac{\opnorm{\BRMatrix}^2}{2 \constYoung \tau^3} \TDLyap(\BStrat_k,
		\bq_k).
	\end{align}
\end{subequations}
Combining equations \eqref{eqn:V_taylor}
and~\eqref{EqnVFirstOrderBounds} yields the claimed
inequality~\eqref{EqnCoupledDriftStrat}.


\paragraph{Proof of the $\TDLyap$-inequality~\eqref{EqnCoupledDriftQ}}

Let us restate the TD update~\eqref{EqnTDUpdate} for convenience:
\begin{align*}
	\BAvgPayoff_{k+1} & = \BAvgPayoff_k + \alpha_k \:
	\frac{\boldsymbol{E}(\BAction_k)}{\BStrat_k(\BAction_k)} \big (
	\BRMatrix \ \boldsymbol{e}(\BAction_k) - \BAvgPayoff_k \big ) \qquad
	\mbox{for $k = 1, 2, \ldots$.}
\end{align*}
Adding and subtracting terms on each side of this equation yields
\begin{align*}
	\BAvgPayoff_{k+1} - \BRMatrix \BStrat_{k+1} &= \BAvgPayoff_k -
	\BRMatrix \BStrat_k + \BRMatrix(\BStrat_k - \BStrat_{k+1})+ \alpha_k
	\: \frac{\boldsymbol{E}(\BAction_k)}{\BStrat_k(\BAction_k)} \big (
	\BRMatrix \ \boldsymbol{e}(\BAction_k) - \BAvgPayoff_k \big )
	\\ &\hspace{15em} + \alpha_k
	\frac{\boldsymbol{E}(\BAction_k)}{\BStrat(\BAction_k)} \BRMatrix
	\BStrat_k - \alpha_k
	\frac{\boldsymbol{E}(\BAction_k)}{\BStrat(\BAction_k)} \BRMatrix
	\BStrat_k.
\end{align*}
After some rearranging, by combining the last three terms above, we
find that
\begin{align*}
	\Btderror{k+1} &= (1-\alpha_k)\Btderror{k} + \alpha_k \big
	(\boldsymbol{I} -
	\frac{\boldsymbol{E}(\BAction_k)}{\BStrat_k(\BAction_k)} \big )
	\Btderror{k} \\ &\hspace{10em}+ \alpha_k
	\frac{\boldsymbol{E}(\BAction_k)}{\BStrat_k(\BAction_k)} \BRMatrix
	(\boldsymbol{e}(\BAction_k) - \BStrat_k) + \BRMatrix(\BStrat_k -
	\BStrat_{k+1}), \numberthis\label{EqnTDError}
\end{align*}
where $\Btderror{k} \defn \BAvgPayoff_k - \BRMatrix \BStrat_k$ and
$\boldsymbol{I}$ is a diagonal block matrix of the same dimensions as
$\BRMatrix$ with the $i^{\text{th}}$ diagonal block being the identity
matrix in $\R^{|\actions^i| \times |\actions^i|}$ for $i \in
\nodes$. Denote by $\E_k$ the conditional expectation operator given
$\BStrat_k$ and $\Btderror{k}$.  Note that the second and third terms
in~\Cref{EqnTDError} have zero mean conditioned on $\BStrat_k$ and
$\Btderror{k}$. As a result, when we take the mean conditioned on
$\BStrat_k$ and $\Btderror{k}$ of the square of the second norm on
both sides, four out of six cross-terms on the right-hand side
immediately vanish.

Moreover, the expected value of one additional cross-term is also
zero; in particular, we have
\begin{align*}
	E_k \big \langle (\boldsymbol{I} -
	\frac{\boldsymbol{E}(\BAction_k)}{\BStrat_k(\BAction_k)} \big )
	\Btderror{k}, \frac{\boldsymbol{E}(\BAction_k)}{\BStrat_k(\BAction_k)}
	\BRMatrix & (\boldsymbol{e}(\BAction_k) - \BStrat_k) \big \rangle \\
	& = E_k \big \langle {
		\frac{\boldsymbol{E}(\BAction_k)}{\BStrat_k(\BAction_k)}}^{\intercal}(\boldsymbol{I}
	- \frac{\boldsymbol{E}(\BAction_k)}{\BStrat_k(\BAction_k)} \big )
	\Btderror{k}, \BRMatrix (\boldsymbol{e}(\BAction_k) - \BStrat_k) \big
	\rangle.
\end{align*}
In the above inner product, as we are conditioning on $\BStrat_k$ and
$\Btderror{k}$, the $i^{\text{th}}$ block of the first term only
depends on player $i$'s action $\action^i_k$, while the
$i^{\text{th}}$ block of the second term does not depend on
$\action^i_k$ as $\RMatrix{i}{i}=0$. As the actions played by each
player are independent, the expected value of the inner product above
is zero. Dropping all the vanishing cross-terms after taking the mean
conditioned on $\BStrat_k$ and $\Btderror{k}$ of the square of the
second norm on both sides of~\Cref{EqnTDError} leaves us with
\begin{multline}
	\label{EqnW2Expand}
	\E_k [ \| \Btderror{k+1}\|_2^2 ] \leq (1-\alpha_k)^2 \ \E_k \big [ \|
	\Btderror{k}\|_2^2 \big ] + \alpha_k^2 \ \E_k [ \|(\boldsymbol{I} -
	\frac{\boldsymbol{E}(\BAction_k)}{\BStrat_k(\BAction_k)} \big )
	\Btderror{k}\|_2^2] \\
	+ \opnorm{\BRMatrix}^2 \ \E_k [ \|\BStrat_k - \BStrat_{k+1} \|_2^2 ] +
	\alpha_k^2 \ \E_k \big [ \|
	\frac{\boldsymbol{E}(\BAction_k)}{\BStrat_k(\BAction_k)} \BRMatrix
	(\boldsymbol{e}(\BAction_k) - \BStrat_k) \|_2^2 \big ] \\
	+ (1 - \alpha_k) \langle \Btderror{k}, \BRMatrix (\BStrat_k -
	\BStrat_{k+1}) \rangle .
\end{multline}
Our next step is to examine each of these terms and upper bound them
appropriately. \\

\myunderparagraph{Bounding the second term:} The submultiplicativity
of the operator norm gives
\begin{align*}
	\E_k \big [ \|(\boldsymbol{I} -
	\frac{\boldsymbol{E}(\BAction_k)}{\BStrat_k(\BAction_k)} \big )
	\Btderror{k}\|_2^2 \big ] &\leq E_k \big [ \opnorm{(\boldsymbol{I} -
		\frac{\boldsymbol{E}(\BAction_k)}{\BStrat_k(\BAction_k)} \big )
	}^2 \big ] \| \Btderror{k}\|_2^2 \\ & \leq E_k \big [
	\frac{1}{\min_{i \in \nodes} \strat[i]_k(\action^i_k)^2}\big ] \; \:
	\| \Btderror{k}\|_2^2 \\
	& \leq \frac{\Amax \NumPlayers \| \Btderror{k}\|_2^2}{\delta}.
\end{align*}
The final bound follows by upper bounding each term in the expectation
and noting that there at most $\Amax \NumPlayers$ number of terms in
the expectation.

\myunderparagraph{Bounding the third term:} We bound this term via
H\"{o}lder's inequality, which leads to
\begin{align*}
	\opnorm{\BRMatrix}^2 \ \E_k [ \|\BStrat_k - \BStrat_{k+1} \|_2^2 ]
	&= \beta_k^2 \opnorm{\BRMatrix}^2 \|\sigma_{\tau}(\BAvgPayoff_k) -
	\BStrat_k \|_2^2 \\ &\leq {4 \NumPlayers \beta_k^2
		\opnorm{\BRMatrix}^2}.
\end{align*}

\myunderparagraph{Bounding the fourth term:} Note that $ \| \BRMatrix
(\boldsymbol{e}(\BAction_k) - \BStrat_k)\|_{\infty} \leq 2\NumPlayers$
due to the assumption that each element of $\BRMatrix$ lies in
$[-1,1]$.  Using the submultiplicativity of the matrix norm, we have
\begin{align*}
	\E_k \big [ \|
	\frac{\boldsymbol{E}(\BAction_k)}{\BStrat_k(\BAction_k)} \BRMatrix
	(\boldsymbol{e}(\BAction_k) - \BStrat_k) \|_2^2 \big ] & \leq \E_k
	\big [
	\opnorm{\frac{\boldsymbol{E}(\BAction_k)}{\BStrat_k(\BAction_k)}}^2
	\opnorm{\BRMatrix}^2 \| \boldsymbol{e}(\BAction_k) - \BStrat_k
	\|_2^2 \big ] \\
	& \overset{(i)}{\leq} 4 \NumPlayers \opnorm{\BRMatrix}^2 \E_k \big [
	\opnorm{\frac{\boldsymbol{E}(\BAction_k)}{\BStrat_k(\BAction_k)}}^2
	\big ] \\
	& \overset{(ii)}{\leq} \frac{4 \NumPlayers^2 \Amax
		\opnorm{\BRMatrix}^2}{\delta}.
\end{align*}
In inequality (i), we use H\"{o}lder's inequality the same way as
in equation~\eqref{EqnStratOneStep}; inequality (ii) follows from the way we
bounded the second term above.

\myunderparagraph{Bounding the fifth term:} Plugging in the strategy
update equation~\eqref{EqnPartialStrategy} and performing some algebra
yields
\begin{align*}
	\langle \Btderror{k}, \BRMatrix (\BStrat_k - \BStrat_{k+1}) \rangle &=
	-\beta_k \langle \Btderror{k}, \BRMatrix(\sigma_{\tau}(\BAvgPayoff_k)
	- \BStrat_k) \rangle \\
	& = -\beta_k \langle \Btderror{k},
	\BRMatrix(\sigma_{\tau}(\BAvgPayoff_k) -
	\sigma_{\tau}(\BRMatrix\BStrat_k)) \rangle \\ & \hspace{10em}- \beta_k \langle
	\Btderror{k}, \BRMatrix(\sigma_{\tau}(\BRMatrix\BStrat_k) - \BStrat_k)
	\rangle.
\end{align*}
We now use the Cauchy-Schwarz inequality to bound each of these
terms. Using the fact that the function $u \mapsto \sigma_{\tau}(u)$
is $1/\tau$-Lipschitz gives
\begin{align*}
	-\langle \Btderror{k}, \BRMatrix(\sigma_{\tau}(\BAvgPayoff_k)
	- \sigma_{\tau}(\BRMatrix\BStrat_k)) \rangle &\leq
	\frac{\|\Btderror{k}\|_2^2 \opnorm{\BRMatrix} }{\tau}.
\end{align*}
Again using the Cauchy-Schwarz inequality results in
\begin{align*}
	-\langle \Btderror{k},
	\BRMatrix(\sigma_{\tau}(\BRMatrix\BStrat_k) - \BStrat_k)
	\rangle &\leq \| \Btderror{k} \|_2 \opnorm{\BRMatrix} \|
	\sigma_{\tau}(\BRMatrix\BStrat_k) - \BStrat_k\|_2.
\end{align*}
Finally, applying Young's inequality for any $c>0$ and parameter
$\constYoung \in (0,1)$ and \cref{LemVpropsPartial}(a) gives
\begin{align*}
	\langle \Btderror{k}, \BRMatrix(\sigma_{\tau}(\BRMatrix\BStrat_k) -
	\BStrat_k) \rangle &\leq \frac{1}{2c \constYoung}\| \Btderror{k}
	\|_2^2 \opnorm{\BRMatrix}^2 + \frac{c \constYoung}{2}\|
	\sigma_{\tau}(\BRMatrix\BStrat_k) - \BStrat_k\|_2^2 \\
	& \leq \frac{1}{2c \constYoung}\| \Btderror{k} \|_2^2
	\opnorm{\BRMatrix}^2 + \frac{c \constYoung
		\stratLyap(\BStrat_k)}{\tau}.
\end{align*}
Choosing $c=\tau$ gives
\begin{align*}
	\langle \Btderror{k}, \BRMatrix (\BStrat_k - \BStrat_{k+1})
	\rangle &\leq \frac{2\beta_k}{ \tau \constYoung}\|
	\Btderror{k} \|_2^2 \opnorm{\BRMatrix}^2 + \beta_k\constYoung
	\stratLyap(\BStrat_k).
\end{align*}


\myunderparagraph{Combining the upper bounds:}
Putting together the pieces yields
the claimed $\TDLyap$-inequality~\eqref{EqnCoupledDriftQ}---that is
\begin{multline*}
	\E\TDLyapVal{k+1} \leq (1-\alpha_k)^2 \E\TDLyapVal{k} + \alpha_k^2
	\frac{\Amax \NumPlayers \E\TDLyapVal{k} }{\delta} + {4 \NumPlayers
		\beta_k^2 \opnorm{\BRMatrix}^2}+ \alpha_k^2 \frac{4 \NumPlayers^2
		\Amax \opnorm{\BRMatrix}^2}{\delta} \\
	+ (1-\alpha_k)\frac{2\beta_k}{ \tau \constYoung} \opnorm{\BRMatrix}^2
	\E\TDLyapVal{k} + (1-\alpha_k)\beta_k\constYoung
	\E\stratLyap(\BStrat_k) .
\end{multline*}


\subsection{Proof of~\Cref{LemNGBoundPartial}}
\label{SecLemLowerBound}
We begin with the drift inequality~\eqref{EqnTotalDrift} on the total
Lyapunov values $\{\TotLyapVal{k}\}_{k \geq 1}$. By solving this drift
inequality, we obtain bounds on the expected Nash gap at time $K+1$
for a trajectory $\{\BStrat_k\}_{k \geq 1}$ that remains $\delta$-good
until time $K$.
\begin{lemma}
	\label{LemDriftPartial}
	Suppose that the iterates $\{\BStrat_k\}_{k \geq 1}$ are $\delta$-good
	for times $k = 1, \ldots, K$.  Then the following properties hold:
	\begin{enumerate}[label=(\alph*)]
		\item For the constant stepsize $\beta_k \equiv \beta$, we have
		\begin{multline*}
			\E \nashg(\BStrat_{K+1}) \leq (1-\beta(1-2\constYoung))^K
			\TotLyapVal{1} + \frac{2\NumPlayers}{1-2\constYoung} \tau \log
			\Amax \\+ \frac{5 \NumPlayers^2 \Amax \opnorm{\BRMatrix}^2
				\beta}{(1-2\constYoung )\timescalesep^2(\constYoung,\tau)
				\delta}.
		\end{multline*}
		\item For the inverse polynomial stepsize $\beta_k = \frac{
			\beta}{(k+\koff)^\myexp}$ for some exponent $\myexp \in (0,1)$ and
		offset \mbox{$\koff = \big \lceil
			(\frac{2\myexp}{\beta})^{1/(1-\myexp)} \big \rceil$,} we have
		\begin{multline*}
			\E \nashg(\BStrat_{K+1}) \leq
			\frac{\exp{-\frac{(1-2\constYoung)\beta}{(1-\myexp)} (K + \koff +
					1)^{1-\myexp}}}{\exp{-\frac{(1-2\constYoung)\beta}{(1 -
						\myexp)} (1 + \koff)^{\myexp}} } \TotLyapVal{1}+
			\frac{2\NumPlayers}{1-2\constYoung} \tau \log \Amax \\
			+ \frac{10 \NumPlayers^2 \Amax \opnorm{\BRMatrix}^2
				\beta}{(1-2\constYoung)\timescalesep^2(\constYoung,\tau)
				\delta (K+\koff)^\myexp } .
		\end{multline*}
	\end{enumerate}
\end{lemma}
The results claimed in~\Cref{LemNGBoundPartial} then follow
from~\Cref{LemDriftPartial} by setting each term in the upper bounds
on the Nash gap to be equal to $\mbox{$\frac{\epsilon}{3}$}$.

\begin{proof}[Proof of \cref{LemDriftPartial}]
	
	We split our analysis into different cases, depending on the choice of
	stepsizes.
	
	\paragraph{Constant stepsizes}  For constant stepsizes
	$\beta_k \equiv \beta \in (0,1)$, by iterating the drift inequality
	\eqref{EqnTotalDrift}, we find that
	\begin{multline}
		\label{Eqn:ConstStepTkBoundPartial}
		\E \TotLyapVal{K+1} \leq \big (1-\beta \big (1-2\constYoung \big )
		\big )^{K} \E\TotLyapVal1 + {2 \NumPlayers} \tau \log \Amax
		\sum_{j=1}^K \beta \big ( 1-(1-2\constYoung)\beta\big )^{K-j} \\
		+ \frac{5 \Amax \NumPlayers^2
			\opnorm{\BRMatrix}^2}{\timescalesep^2(\constYoung,\tau)
			\delta}\sum_{j=1}^K \beta^2 \big ( 1-(1-2\constYoung) \beta
		\big)^{K-j}.
	\end{multline}
	
	
	\myunderparagraph{Bounding the second term:} Consider the sequence
	$\{u_k\}_{k \geq 1}$ defined with the initialization $u_1 = 0$ and the
	recursion \mbox{$u_{k+1} = \big (1-(1-2\constYoung)\beta\big ) u_k +
		(1-2\constYoung) \beta$.} Then for $K = 1, 2, \ldots$, we have
	\begin{align*}
		u_{K + 1} = \sum_{j=1}^K (1-2\constYoung) \beta \big( 1 -( 1 - 2
		\constYoung) \beta \big)^{K-j}.
	\end{align*}
	Since $\beta \leq 1$ and $\constYoung<0.5$, it follows from an
	inductive argument that $u_k \leq 1$ for all $k = 1, 2, \ldots$.
	Using this fact, we have
	\begin{align*}
		2 \NumPlayers \tau \log \Amax \sum_{j=1}^K \beta \big( 1- (1 - 2
		\constYoung) \beta \big)^{K-j} \leq \frac{2 \NumPlayers}{1 - 2
			\constYoung} \tau \log \Amax .
	\end{align*}
	
	
	\myunderparagraph{Bounding the third term:} Note that
	\begin{align*}
		\sum_{j=1}^K \big ( 1 -(1 - 2 \constYoung) \beta \big )^{K-j} =
		\sum_{j=0}^{K} \big (1 -( 1 -2 \constYoung) \beta \big )^j \leq
		\frac{1}{(1 - 2\constYoung) \beta}.
	\end{align*}
	Using this fact, it follows that
	\begin{align*}
		\sum_{j=1}^K \beta^2 \big ( 1 - (1 - 2\constYoung)\beta\big )^{K-j}
		\leq \frac{\beta}{1-2\constYoung}.
	\end{align*}
	Thus, equation~\eqref{Eqn:ConstStepTkBoundPartial} can be simplified
	to the form in~\Cref{LemDriftPartial}(a).
	
	\paragraph{Inverse polynomial stepsizes}
	By assuming that the iterates $\BStrat_k$ are $\delta$-good for all $k = 0, 1,
	\ldots, K$, we have
	\begin{align*}
		\TotLyapVal{K+1} &\leq \prod_{i = 1}^K \big ( 1 -(1 -2
		\constYoung)\beta_i\big ) \TotLyapVal1 + 2\NumPlayers \tau \log \Amax
		\sum_{i=1}^K \beta_i \prod_{j = i+1}^K \big ( 1-(1 - 2
		\constYoung)\beta_j\big ) \\
		& \hspace{13em} + \frac{5 \NumPlayers^2 \Amax
			\opnorm{\BRMatrix}^2}{\timescalesep^2(\constYoung, \tau) \delta}
		\sum_{i=1}^K \beta_i^2 \prod_{j=i+1}^K \big ( 1 -(1 - 2\constYoung)
		\beta_j\big ).
	\end{align*}
	
	
	\myunderparagraph{Bounding the first term:}
	In this case, we have
	\begin{align*}
		\prod_{i=1}^K \big ( 1-(1-2\constYoung)\beta_i\big ) &= \exp \big (
		\sum_{i=1}^K \log \big ( 1-(1-2\constYoung)\beta_i\big ) \big )
		\\ &\overset{(i)}{\leq} \exp \big ( -\sum_{i=1}^K
		(1-2\constYoung)\beta_i \big ) \\ &\overset{(ii)}{\leq} \exp \big (
		- (1-2\constYoung)\beta \int_{1}^{K+1} \frac{1}{(x+\koff)^\myexp} dx
		\big ) \\ &= \exp \big ( - \frac{(1-2\constYoung)\beta}{(1-\myexp)}
		\big ( (K+\koff+1)^{1-\myexp} - (1+\koff)^{\myexp}\big ) \big ).
	\end{align*}
	where step (i) follows from the inequality $\log(1-x) \leq - x$ for $x
	\in [0,1)$, whereas step (ii) follows from bounding the Riemann sum.
	
	\myunderparagraph{Bounding the second term:} Consider the sequence
	$\{u_k\}_{k \geq 1}$ defined by the initialization $u_1=0$ followed by
	the recursion \mbox{$u_{k+1} = \big (1-(1-2\constYoung)\beta_k\big )
		u_k + (1-2\constYoung)\beta_k$.}  Then for $K = 1, 2, \ldots$, we
	have
	\begin{align*}
		u_{K+1} = \sum_{i=1}^K (1-2\constYoung)\beta_i \prod_{j=i+1}^K
		\big( 1-(1-2\constYoung)\beta_j \big).
	\end{align*}
	Note that $\beta_k \leq 1$ for $k = 1, 2 \ldots$ and $\constYoung <
	\frac{1}{2}$ by assumption; therefore, it follows from an inductive
	argument that $u_k \leq 1$ for all $k = 1, 2, \ldots$.  Using this
	fact, we have the bound
	\begin{align*}
		2 \NumPlayers \tau \log \Amax \sum_{i=1}^K \beta_i \prod_{j=i+1}^K
		\big( 1 -(1 - 2 \constYoung)\beta_j \big) \leq \frac{2 \NumPlayers
		}{1 - 2 \constYoung} \tau \log \Amax.
	\end{align*}
	
	\myunderparagraph{Bounding the third term:} Consider the sequence
	$\{u_k\}_{k \geq 1}$ defined by the initialization $u_1=0$ and the
	recursion \mbox{$u_{k+1} = \big (1-(1-2\constYoung)\beta_k\big ) u_k +
		\beta_k^2$.} Then for $K = 1, 2, \ldots$, we have
	\begin{align*}
		u_{K+1} = \sum_{i=1}^K \beta_i^2 \prod_{j=i+1}^K \big( 1 - (1 - 2
		\constYoung) \beta_j \big).
	\end{align*}
	It can be shown by induction \cite[see p.36]{chen2023lyapunov} that
	$u_k \leq \frac{2}{(1-2\constYoung)} \beta_k$ for $k = 1, 2, \ldots$,
	from which it follows that
	\begin{align*}
		\frac{5 \NumPlayers^2\Amax\opnorm{\BRMatrix}^2}{
			\timescalesep^2(\constYoung, \tau) \delta} \sum_{i=1}^K \beta_i^2
		\prod_{j=i+1}^K \big ( 1-(1 - 2 \constYoung)\beta_j\big ) \leq
		\frac{10 \NumPlayers^2 \Amax \opnorm{\BRMatrix}^2 \beta}{(1 - 2
			\constYoung)\timescalesep^2(\constYoung,\tau) \delta (K +
			\koff)^\myexp } .
	\end{align*}
\end{proof}


\subsection{Proof of~\Cref{LemUpperBound}}
\label{SecLemUpperBound}

We prove the following more general claim:
\begin{lemma}
	\label{LemUpperBoundFull}
	Suppose that~\Cref{AlgPartial} is initialized with uniform initial
	strategies $\BStrat_1$ and the stepsize parameter initialized
	according to equation~\eqref{EqnKeyTauBeta}. The quantity
	$\uplowratio(\offset)$ denotes the ratio between the lower bound on
	$\Kgood(\beta,\delta)$ and $K(\epsilon,\delta,\constYoung)$.
	\begin{enumerate}[label=(\alph*)]
		\item For any $\uplowratio(\offset) \in (1, 1 + \offset)$, there
		exists a choice \mbox{$\constYoung(\offset) \in (0,\frac{1}{2})$}
		such that~\Cref{AlgPartial} with the constant stepsize $\beta_k
		\equiv
		\beta(\epsilon,\constYoung(\offset),\delta(\epsilon,\offset))$ and
		$\delta(\epsilon,\offset) \defn
		\big(\frac{\epsilon}{3\TotLyapVal{1}}\big )^{1+\offset}
		\frac{1}{\Amax}$ satisfies
		\begin{subequations}
			\begin{align*}
				\Kgood \big (\beta(\epsilon, \constYoung(\epsilon,
				\offset),\delta(\epsilon, \offset)), \delta(\epsilon,\offset)
				\big ) \geq \frac{\uplowratio(\offset) \log \big (
					\frac{\epsilon}{3 \TotLyapVal{1} }\big )}{\log\big(1- (1 - 2
					\constYoung(\offset))\beta(\epsilon, \constYoung(\offset),
					\delta(\epsilon,\offset))\big )}.
			\end{align*}
			\item Consider the inverse polynomial stepsize $\beta_k = \frac{
				\beta(\epsilon, \constYoung,
				\delta(\epsilon,\offset))}{(k+\koff)^\myexp}$ for some exponent
			$\myexp \in (0,1)$, the offset \mbox{$\koff = \big \lceil
				(\frac{2\myexp}{\beta})^{1/(1-\myexp)} \big \rceil$}, and
			$\delta(\epsilon,\offset) \defn \big( \frac{\epsilon}{3
				\TotLyapVal{1}}\big )^{1+\offset }\frac{e}{\Amax}$.  Then for any
			choice of $\uplowratio(\offset) \in (1,\sqrt{1+\offset})$ such that
			$\log(1-\beta_1) \geq - \uplowratio(\offset) \beta_1$, there exists
			some \mbox{$\constYoung(\offset) \in (0,\frac{1}{2})$} such that
			\begin{multline*}
				\Kgood \big (\beta(\epsilon,\constYoung(\epsilon,
				\offset),\delta(\epsilon,\offset)),\delta(\epsilon,\offset) \big )
				\\ \geq \big [ \frac{\uplowratio(\offset)}{1 - 2 \constYoung}
				\frac{(1-\myexp)}{\beta(\epsilon,\constYoung(\offset),\delta(\epsilon,\offset))}
				\log \big ( \frac{3 \TotLyapVal{1}}{\epsilon}\big ) +
				(1+\koff)^{1-\myexp} \big]^{\frac{1}{1-\myexp}} - k_0 .
			\end{multline*}
		\end{subequations}
	\end{enumerate}
\end{lemma}
In the above lemma, the quantity $\uplowratio(\offset)>1$ is the
factor by which $\Kgood(\beta,\delta)$ is greater than
$K(\epsilon,\delta,\constYoung)$ (cf. \cref{LemNGBoundPartial}), for
the specific choices of $\delta$ and $\constYoung$ mentioned in the
lemma above. Our proof of \cref{LemUpperBoundFull} exploits the
following auxiliary result:
\begin{lemma}
	\label{LemUpperBound_general}
	For any $\beta,\delta \in (0,1)$, either of the following two
	conditions are sufficient to ensure that the iterates
	$\{\BStrat_k\}_{k \geq 1}$ of \cref{AlgPartial} with timescale
	constant $\beta$ are $\delta$-good up to time $K$.
	\begin{enumerate}[label=(\alph*)]
		\item For the constant stepsize $\beta_k \equiv \beta$, it suffices
		to have \mbox{$K \leq \frac{\log(\Amax \delta)}{\log(1-\beta)}$.}
		\item For \mbox{$\beta \in (0,\frac{1}{2})$}, suppose that $\beta_k
		= \frac{ \beta}{(k+\koff)^\myexp}$ where $\myexp \in (0,1)$,
		$\koff = \big(\frac{2\myexp}{\beta}\big)^{1/(1-\myexp)}$, and
		$\uplowratio >1$ is chosen to ensure that $\log(1-\beta_1) \geq -
		\uplowratio \beta_1$. Then it suffices to have
		\begin{align*}
			(K + \koff)^{1-\myexp} &\leq \frac{(1-\myexp)}{\uplowratio
				\beta} \log \big ( \frac{e}{\Amax \delta }\big ) +
			(1+\koff)^{1-\myexp},
		\end{align*}
	\end{enumerate}
\end{lemma}
\noindent See~\Cref{AppLemUpperBound_general} for the
proof. Equivalently, the upper bounds on $K$ in
\cref{LemUpperBound_general} are lower bounds on
$\Kgood(\beta,\delta)$.\\

Now, we choose $\delta$ and $\constYoung$ as a function of $\epsilon$
and $\offset$ to ensure that these lower bounds on
$\Kgood(\beta,\delta)$ for the choice of $\beta$ in
\cref{LemNGBoundPartial} are bigger than $K(\epsilon, \delta,
\constYoung)$. The statement in \cref{LemUpperBound} follows for the
specific choice of $\delta(\epsilon,\offset) = \big( \frac{\epsilon}{3
	\TotLyapVal{1}}\big )^{1+\offset}\frac{1}{\Amax}$ for the case of
constant stepsizes and $\delta(\epsilon,\offset) = \big(
\frac{\epsilon}{3 \TotLyapVal{1}} \big)^{1+\offset}\frac{e}{\Amax}$
for the case of sublinearly decaying stepsizes. The exact value of
$\constYoung$ will be discussed below for each stepsize.

\paragraph{Constant stepsizes}
In this case, we have $\beta_k \equiv \beta$ for any $\beta \in
(0,1)$. Let $\uplowratio \in (1,1+\offset)$. There exists a small
enough $\constYoung(\offset) \in (0,\frac{1}{2})$ such that
\begin{align*}
	\frac{\log(1-\beta)}{\log \big ( 1 - \beta \big ( 1- 2
		\constYoung(\offset) \big ) \big ) } \leq \frac{1 +
		\offset}{\uplowratio}.
\end{align*}
Applying this bound yields
\begin{align*}
	\frac{\log(\Amax \ \delta(\epsilon,\offset))}{\log(1 - \beta)} & = (1
	+ \offset) \frac{\log \big ( \frac{\epsilon}{3 \TotLyapVal{1}}\big
		)}{\log(1-\beta)} \geq \frac{\uplowratio \log \big (
		\frac{\epsilon}{3 \TotLyapVal{1}}\big )}{\log \big ( 1 - \beta \big
		( 1- 2 \constYoung(\offset) \big ) \big )}.
\end{align*}
Therefore, the iterates $\{\BStrat_k\}_{k \geq 1}$ remain
$\delta(\epsilon,\offset)$-good up until time
\begin{align*}
	\frac{\uplowratio \log \big ( \frac{\epsilon}{3 \TotLyapVal{1}}\big
		)}{\log \big ( 1 - \beta \big ( 1- 2 \constYoung(\offset) \big )
		\big )}.
\end{align*}
The statement of \cref{LemUpperBoundFull} follows for $\beta =
\beta(\epsilon,\delta(\epsilon,\offset),\constYoung(\offset))$.


\paragraph{Inverse polynomial stepsizes}

Let $\uplowratio(\offset) \in (1,\sqrt{1+\offset})$. If $\beta_k =
\frac{ \beta(\epsilon,\delta,\constYoung)}{(k+\koff)^\myexp}$ where
$\myexp \in (0,1)$, $\koff =
\big(\frac{2\myexp}{\beta}\big)^{1/(1-\myexp)}$, as
$\beta(\epsilon,\delta,\constYoung) \rightarrow 0$ as $\constYoung
\rightarrow 0$ (equation \eqref{EqnKeyTauBeta}), there will be a small
enough $\constYoung(\offset)$ that satisfies
$\log(1-\beta(\epsilon,\delta,\constYoung(\offset))) \geq -
\uplowratio(\offset) \beta(\epsilon,\delta,\constYoung(\offset))$.
Find a small enough $\constYoung(\offset) \in (0,\frac{1}{2})$ such
that we also have
\begin{align}
	\label{EqnConstYoungDecayStep}
	1 + \offset \geq
	\frac{\uplowratio^2(\offset)}{1-2\constYoung(\offset)} .
\end{align}
Since $\uplowratio^2(\offset) < 1 + \offset$ by assumption, such a
choice for $\constYoung(\offset)$ is possible.  With such a choice, we
have
\begin{align*}
	\frac{1-\myexp}{\uplowratio(\offset)\beta(\epsilon,\delta(\epsilon,\offset),\constYoung(\offset))}
	\log \big ( \frac{e}{\Amax \delta(\epsilon,\offset) } \big) & = \big
	(1 + \offset \big )
	\frac{(1-\myexp)}{\uplowratio(\offset)\beta(\epsilon,\delta(\epsilon,\offset),\constYoung(\offset))}
	\log \big ( \frac{3 \TotLyapVal{1}}{\epsilon}\big ) \\
	& \geq
	\frac{\uplowratio(\offset)(1-\myexp)}{\beta(\epsilon,\delta(\epsilon,\offset),\constYoung(\offset))(1-2\constYoung(\offset))}
	\log \big ( \frac{3 \TotLyapVal{1}}{\epsilon}\big ).
\end{align*}
Therefore, the iterates $\{\BStrat_k\}_{k \geq 1}$ are
$\delta(\epsilon,\offset)$-good for all $K$ such that
\begin{align*}
	(K + \koff)^{1-\myexp} &\leq
	\frac{\uplowratio(\offset)}{1-2\constYoung(\offset)}
	\frac{(1-\myexp)}{\beta(\epsilon,\delta(\epsilon,\offset),\constYoung(\offset))}
	\log \big ( \frac{3 \TotLyapVal{1}}{\epsilon}\big ) +
	(1+\koff)^{1-\myexp}.
\end{align*}


\section{Auxiliary results for~\Cref{ThmMinInfoPoly}}

In this appendix, we collect together the proofs of various auxiliary
results used in proving~\Cref{ThmMinInfoPoly}.


\subsection{Proof of \cref{LemVpropsPartial}}
\label{AppLemVdriftPartial}

We prove each of the two parts in turn.

\paragraph{Proof of part (a)}

For any joint mixed strategy $\BStrat$, define $\BStrat^{-i}$ to be
the collection of mixed strategies of all players except player
$i$. Note that the average payoff function for player $i$,
$\AvgPayoff{i}(\strat^{-i})$ is a function of $\BStrat^{-i}$. Define
the function
\begin{align*}
	F^i(\strat^i, \BStrat^{-i}) = \underset{\pihat^i \in
		\simplex{i}}{\max} \big\{ \big ({{}\pihat^i} - \strat[i] \big
	)^{\intercal} \AvgPayoff{i}(\strat^{-i}) + \tau \EntFun({{}\pihat^i})
	- \tau \EntFun({{}{\pi}^i}) \big\}, \quad i=1,2,
\end{align*}
and observe that $F^i(\BStrat) \geq 0$. Additionally,
$F^i(\sigma_{\tau}(\AvgPayoff{i}(\strat^{-i})),
\BStrat^{-i})=0$. Therefore, by the first-order optimality condition
and because the minimizer is in the relative interior of
$\simplex{i}$, we must have
\begin{align*}
	\inprod{\nabla_{\strat[i]}
		F^i(\sigma_{\tau}(\AvgPayoff{i}(\strat^{-i}))), \BStrat^{-i})}{\pi_1
		- \pi_2} = 0 \quad \mbox{for all mixed strategies $\pi_1, \pi_2 \in
		\simplex{i}$.}
\end{align*}
Note that as the negative of the Shannon entropy is $1$-strongly
convex on the probability simplex; consequently, the function $F^i$ is
$\tau$-strongly convex with respect to $\strat[i]$ uniformly for all
$\BStrat^{-i}$. Therefore, we have
\begin{align*}
	F^i(\strat[i], \BStrat^{-i}) & = F^i(\strat[i], \BStrat^{-i}) -
	F^i(\AvgPayoff{i}(\strat^{-i})), \BStrat^{-i}) \\
	& {\geq} \; \inprod{\nabla_{\strat[i]}
		F^i(\sigma_{\tau}(\AvgPayoff{i}(\strat^{-i})), \BStrat^{-i})}
	{\strat[i] - \sigma_{\tau}(\AvgPayoff{i}(\strat^{-i}))} +
	\frac{\tau}{2} \big \| \sigma_{\tau}(\AvgPayoff{i}(\strat^{-i})) -
	\strat[i] \big \|_2^2 \\
	& = \frac{\tau}{2} \big \| \sigma_{\tau}(\AvgPayoff{i}(\strat^{-i})) -
	\strat[i] \big \|_2^2.
\end{align*}
Using this lower bound, we can write
\begin{align*}
	\frac{\tau}{2} \| \sigma_{\tau}(\BAvgPayoff(\BStrat)) - \BStrat
	\|_2^2 = \frac{\tau}{2} \sum_{i \in [\NumPlayers]} \big \|
	\sigma_{\tau}(\AvgPayoff{i}(\strat^{-i})) - \strat[i] \big \|_2^2 \;
	\leq \; \sum_{i \in \NumPlayers} F^i(\strat[i], \BStrat^{-i}) &
	\stackrel{(i)}{\leq} \stratLyap(\BStrat),
\end{align*}
where step (i) follows from the non-negativity of the Shannon entropy.
This completes the proof.

\newcommand{\piopti}{\ensuremath{\strat[i]_\star}}

\paragraph{Proof of part (b)}
Let $\boldsymbol{v} = (v_i)_{i \in [\NumPlayers]} \in \otimes_{i \in
	[\NumPlayers]} \R^{|\actions^i|}$. For each $i \in [\NumPlayers]$,
consider the constrained optimization problem
\begin{align*}
	\max_{\strat^i \in \R^{|\A_i|}_{+}} \big \{ \inprod{\strat[i]}{v^i} + \tau
	\EntFun(\strat[i]) \big \} \qquad \mbox{such that $\sum_{\action \in \actions^i}
		\strat[i](\action) = 1$.}
\end{align*}
The gradient $\|\nabla \EntFun(\strat[i])\|_2$ diverges whenever any
element of $\strat[i]$ approaches zero, so that we can argue that the
optimum will be achieved at a vector $\piopti$ with strictly positive
co-ordinates.  We introduce a Lagrange multiplier $\lambda$ for the
sum-constraint.  The associated KKT conditions imply the optimum
$\piopti$ satisfies the condition
\begin{align*}
	v^i + \tau \nabla \EntFun(\piopti) + \lambda \boldsymbol{e} = 0,
\end{align*}
where $\boldsymbol{e}$ is a vector of all ones, and $\lambda$ is the
Lagrange multiplier.  Solving for the optimum and using the fact that
$\lambda$ is chosen to ensure that $\piopti$ satisfies the
normalization constraint yields $\piopti = \sigma_{\tau}(v^i)$.
Equivalently, we have shown that $v^i + \tau \nabla
\EntFun(\sigma_{\tau}(v^i)) + \lambda \boldsymbol{e} = 0$, from which
it follows that
\begin{align*}
	\inprod{v^i + \tau \nabla \EntFun(\sigma_{\tau}(v^i))}{\strat_1 -
		\strat_2} & = 0 \qquad \mbox{for any pair $\strat_1, \strat_2 \in
		\simplex{i}$.}
\end{align*}
The above result can also be expressed as
\begin{align}
	\label{EqnPartialOptTrick}
	\big \langle \boldsymbol{v} + \tau \nabla
	\EntFun(\sigma_{\tau}(\boldsymbol{v})), \BStrat_1 - \BStrat_2 \rangle
	= 0 \qquad \mbox{for any pair $\BStrat_1, \BStrat_2 \in \BigSimplex$,}
\end{align}
where $\nabla \EntFun(\sigma_{\tau}(\boldsymbol{v})) = ( \nabla
\EntFun(\sigma_{\tau}(v^i)))_{i \in [\NumPlayers]}$. We will now use
this result to prove \Cref{LemVpropsPartial}(b).

Using the gradient we evaluated using Danskin's
Theorem~\eqref{EqnDanskinGradient}, we can write
\begin{align*}
	\big \langle \nabla_{\strat} \stratLyap(\BStrat), \sigma_{\tau}(
	\boldsymbol{u}) - \sigma_{\tau}(\BAvgPayoff(\BStrat)) \big \rangle &
	= \big \langle \BRMatrix^{\intercal}
	\sigma_{\tau}(\BAvgPayoff(\BStrat)), \sigma_{\tau}(\boldsymbol{u})
	-\sigma_{\tau}(\BAvgPayoff(\BStrat)) \big \rangle
	\\ &\overset{(i)}{=} \big \langle
	\sigma_{\tau}(\BAvgPayoff(\BStrat)), \BRMatrix
	\ \sigma_{\tau}(\boldsymbol{u}) \big \rangle \\ &\overset{(ii)}{=}
	\big \langle \sigma_{\tau}(\BAvgPayoff(\BStrat)) -
	\sigma_{\tau}(\boldsymbol{u}), \BRMatrix
	\ \sigma_{\tau}(\boldsymbol{u}) \big \rangle,
\end{align*}
where steps (i) and (ii) follow from the zero-sum
property~\eqref{EqnZeroSum}. Next, choosing $\boldsymbol{v} =
\boldsymbol{u}$ and applying the earlier
result~\eqref{EqnPartialOptTrick} gives
\begin{align*}
	\big \langle \sigma_{\tau}(\BAvgPayoff(\BStrat)) -
	\sigma_{\tau}(\boldsymbol{u}), \BRMatrix
	\ \sigma_{\tau}(\boldsymbol{u}) \big \rangle &= \big \langle
	\sigma_{\tau}(\boldsymbol{u}) - \sigma_{\tau}(\BAvgPayoff(\BStrat)),
	\tau \nabla \EntFun(\sigma_{\tau}(\boldsymbol{u})) \big \rangle
	\\ & \hspace{5em}+ \big \langle \sigma_{\tau}(\BAvgPayoff(\BStrat))
	- \sigma_{\tau}(\boldsymbol{u}), \BRMatrix
	\ \sigma_{\tau}(\boldsymbol{u}) - \boldsymbol{u} \big \rangle\\ &=
	\nablaVBoundPartialFirst + \nablaVBoundPartialSecond,
\end{align*}
where we define
\begin{align*}
	\nablaVBoundPartialFirst
	& \defn \big \langle   \sigma_{\tau}(\boldsymbol{u}) - \sigma_{\tau}(\BAvgPayoff(\BStrat)), \tau \nabla \EntFun(\sigma_{\tau}(\boldsymbol{u}))  \big \rangle , \quad \mbox{and} \\
	\nablaVBoundPartialSecond & \defn \big \langle
	\sigma_{\tau}(\BAvgPayoff(\BStrat)) - \sigma_{\tau}(\boldsymbol{u}),
	\BRMatrix \ \sigma_{\tau}(\boldsymbol{u}) - \boldsymbol{u} \big
	\rangle.
\end{align*}
In order to complete the proof, it suffices to show that
\begin{align}
	\label{EqnTwoBounds}
	\nablaVBoundPartialFirst \stackrel{(a)}{\leq} \NumPlayers \tau \log
	\Amax, \quad \mbox{and} \quad \nablaVBoundPartialSecond
	\stackrel{(b)}{\leq} \constYoung \stratLyap(\BStrat) +
	\frac{\opnorm{\BRMatrix}^2}{2\constYoung \tau^3} \|
	\BAvgPayoff(\BStrat) - \boldsymbol{u}\|^2_2.
\end{align}
We prove each of these two claims in turn.

\paragraph{Proof of the bound~\eqref{EqnTwoBounds}(a)}
Since $\EntFun$ is a concave function, the first-order tangent bound
implies that
\begin{align*}
	\nablaVBoundPartialFirst &= \sum_{i \in [\NumPlayers]}\big \langle
	\sigma_{\tau}({u}^i) - \sigma_{\tau}(\AvgPayoff{i}(\strat^{-i})),
	\tau\nabla \EntFun(u^i) \big \rangle \\ &\leq \sum_{i \in
		[\NumPlayers]} \tau \EntFun(\sigma_{\tau}(u^i)) - \tau
	\EntFun(\sigma_{\tau}(\AvgPayoff{i}(\strat^{-i}))\\ &
	\stackrel{(i)}{\leq} \sum_{i \in [\NumPlayers]} \tau
	\EntFun(\sigma_{\tau}(u^i)) \\ &\stackrel{(ii)}{\leq}
	\NumPlayers\tau \log \Amax.
\end{align*}
In the above argument, step (i) follows from the non-negativity of
$H$, whereas step (ii) follows from the fact that the discrete entropy
is at most log cardinality of the space.

\paragraph{Proof of the bound~\eqref{EqnTwoBounds}(b)}
We first begin by re-arranging as follows
\begin{align*}
	\nablaVBoundPartialSecond &= \big \langle
	\sigma_{\tau}(\BAvgPayoff(\BStrat)) - \sigma_{\tau}(\boldsymbol{u}),
	\BRMatrix \ \sigma_{\tau}(\boldsymbol{u}) - \boldsymbol{u} \big
	\rangle \\
	& = \big \langle \sigma_{\tau}(\BAvgPayoff(\BStrat)) -
	\sigma_{\tau}(\boldsymbol{u}), \BRMatrix
	\ \sigma_{\tau}(\boldsymbol{u}) - \BRMatrix \ \sigma_{\tau}(\BRMatrix
	\BStrat) \big \rangle \\
	&\hspace{3em}+ \big \langle \sigma_{\tau}(\BAvgPayoff(\BStrat)) -
	\sigma_{\tau}(\boldsymbol{u}), \BRMatrix \ \sigma_{\tau}(\BRMatrix
	\BStrat) - \BRMatrix \BStrat \big \rangle + \big \langle
	\sigma_{\tau}(\BAvgPayoff(\BStrat)) - \sigma_{\tau}(\boldsymbol{u}),
	\BRMatrix \BStrat - \boldsymbol{u} \big \rangle \\
	&= \big \langle \sigma_{\tau}(\BAvgPayoff(\BStrat)) -
	\sigma_{\tau}(\boldsymbol{u}), \BRMatrix \ \sigma_{\tau}(\BRMatrix
	\BStrat) - \BRMatrix \BStrat \big \rangle + \big \langle
	\sigma_{\tau}(\BAvgPayoff(\BStrat)) - \sigma_{\tau}(\boldsymbol{u}),
	\BRMatrix \BStrat - \boldsymbol{u} \big \rangle,
\end{align*}
where the final step follows from the zero-sum
property~\eqref{EqnZeroSum}. Using the fact that
$\sigma_{\tau}(\cdot)$ is $1/\tau$-Lipschitz, the Cauchy-Schwarz
inequality and the submultiplicativity of the matrix norm, we have
\begin{align*}
	\nablaVBoundPartialSecond &\leq \frac{1}{\tau} \opnorm{\BRMatrix} \|
	\sigma_{\tau}(\BRMatrix \BStrat) - \BStrat \|_2\|
	\BAvgPayoff(\BStrat) - \boldsymbol{u}\|_2 + \frac{\|
		\BAvgPayoff(\BStrat) - \boldsymbol{u}\|^2_2}{\tau}.
\end{align*}
Applying the inequality $ab \leq \frac{a^2c + \frac{b^2}{c}}{2}$ with
the choices $a = \| \sigma_{\tau}(\BRMatrix \BStrat) - \BStrat \|_2 $,
$b = \frac{1}{\tau} \opnorm{\BRMatrix} \| \BAvgPayoff(\BStrat) -
\boldsymbol{u}\|_2$ and an arbitrary $c > 0$ combined
with~\Cref{LemVpropsPartial}(a) yields
\begin{align*}
	\nablaVBoundPartialSecond &\leq \frac{1}{2} c \|
	\sigma_{\tau}(\BRMatrix \BStrat) - \BStrat \|_2^2 +
	\frac{\opnorm{\BRMatrix}^2}{2c \tau^2} \| \BAvgPayoff(\BStrat) -
	\boldsymbol{u}\|^2_2 \\
	& \leq \frac{c}{\tau} \stratLyap(\BStrat) +
	\frac{\opnorm{\BRMatrix}^2}{2c \tau^2} \| \BAvgPayoff(\BStrat) -
	\boldsymbol{u}\|^2_2.
\end{align*}
Choosing our free parameter as $c = { \constYoung \tau}$ for some
$\constYoung \in (0,1)$, we find that
\begin{align*}
	\nablaVBoundPartialSecond & \leq \constYoung \stratLyap(\BStrat) +
	\frac{\opnorm{\BRMatrix}^2}{2\constYoung \tau^3} \|
	\BAvgPayoff(\BStrat) - \boldsymbol{u}\|^2_2,
\end{align*}
as claimed.


\subsection{Proof of~\Cref{LemUpperBound_general}}
\label{AppLemUpperBound_general}

From the dynamics~\eqref{EqnPartialStrategy}, we have the elementwise
inequality $\pii_{k+1} \succeq \pii_k (1-\beta_k)$.  Iterating this
inequality for $k = 1, \ldots, K + 1$ yields
\begin{align*}
	\strat[i]_{K+1} & \succeq \frac{1}{\Amax} \prod_{j=1}^{K} (1-\beta_j).
\end{align*}
where we have made use of the lower bound $\BStrat_1 \succeq
\frac{1}{\Amax} \boldsymbol{e}$, as guaranteed by our uniform
initialization.

Therefore, in order to ensure that $\strat[i]_k \succeq \delta$ for
all $k = 1, \ldots, K + 1$, it suffices to have
\begin{align}
	\label{eqn:pik_lb_condition}
	\frac{1}{\Amax} \prod_{k=1}^{K} (1-\beta_k) & \geq \delta.
\end{align}
Here we make use of the fact that $\beta_k \in (0,1)$ for all $k$.  We
now analyze the condition~\eqref{eqn:pik_lb_condition} for our two
different choices of stepsizes.

\paragraph{Constant stepsizes}
For the constant stepsizes $\beta_k \equiv \beta \in (0,1)$,
condition~\eqref{eqn:pik_lb_condition} is equivalent to
\begin{align*}
	(1-\beta)^{K} \geq \Amax \delta, \quad \mbox{or equivalently} \quad K
	\leq \frac{\log (\Amax \delta ) }{\log (1-\beta)}.
\end{align*}

\paragraph{Inverse polynomial stepsizes}
Now consider the inverse polynomial stepsizes $\beta_k = \frac{
	\beta}{(k+\koff)^\myexp}$ for some exponent $\myexp \in (0,1)$ and
offset $\koff =
\big(\frac{2\myexp}{\beta}\big)^{1/(1-\myexp)}$. Choose $\uplowratio
>1$ to be the smallest possible number such that $\log(1-\beta_1) \geq
- \uplowratio \beta_1$. For $\beta < \frac{1}{2}$, we also have
$\uplowratio \beta < 1$. Now, to find an elementwise lower bound for
$\BStrat_{K+1}$, we write
\begin{align*}
	\prod_{i=1}^{K} (1-\beta_i) &= \exp \big ( \sum_{i=1}^K
	\log(1-\beta_i) \big ) \\
	& \overset{(i)}{\geq} \exp \big ( -\sum_{i=1}^K
	\frac{\uplowratio\beta}{(i+\koff)^\myexp} \big ) \\ &= \exp \big ( -
	\frac{\uplowratio\beta}{(1+\koff)^\myexp} \big ) \exp \big (
	-\sum_{i=2}^K \frac{\uplowratio\beta}{(i+\koff)^\myexp} \big )
	\\
	& \overset{(ii)}{\geq} \exp \big(
	-\frac{\uplowratio\beta}{(1+\koff)^\myexp} \big ) \exp \big ( -
	\int_1^K \frac{\uplowratio\beta}{(x+\koff)^\myexp} dx \big ) \\
	& = \exp \big ( - \frac{\uplowratio\beta}{(1+\koff)^\myexp} \big )
	\exp\big ( - \frac{\uplowratio\beta}{1-\myexp} \big(
	(K+\koff)^{1-\myexp} - (1+\koff)^{1-\myexp} \big ) \big) \\
	& \overset{(iii)}{\geq} \exp \big( - 1 \big) \exp\big ( -
	\frac{\uplowratio\beta}{1-\myexp} \big( (K+\koff)^{1-\myexp} -
	(1+\koff)^{1-\myexp} \big) \big),
\end{align*}
where step (i) follows from the inequality $\log(1-\beta_k) \geq
-\uplowratio \beta_k$ for $k \in \N$; step (ii) follows from bounding
the Riemann sum; and step (iii) follows from the fact that
\begin{align*}
	\frac{\uplowratio \beta}{(1+\koff)^\myexp} \leq 1,
\end{align*}
since $\uplowratio \beta <1 $ for a small enough $\beta$. Now,
\begin{align*}
	\exp \big( -1 \big) \exp\big ( -\frac{\uplowratio \beta}{1-\myexp}
	\big ( (K + \koff)^{1-\myexp} - (1 + \koff)^{1-\myexp} \big) \big)
	&\geq \Amax \delta \\
	\implies (K+\koff)^{1-\myexp} - (1 + \koff)^{1-\myexp} \leq
	\frac{1-\myexp}{\uplowratio \beta} \log \big( \frac{e}{\Amax \delta}
	\big).
\end{align*}


\section{Analyzing $\opnorm{\BRMatrix}$ for a $k$-regular ring graph}
\label{SeckRegularRMatrix}

Recall that the weighted adjacency matrix $\Amat$ is a circulant
matrix for a $k$-regular graph. Let the first row of $\Amat$ be
$(\cir_1, \ldots, \cir_{\NumPlayers})$, where $\cir_m \in \{-1,0,1\}$
for $m \in \nodes$. We know that $\cir_1=0$ as the diagonal elements
are zero for a weighted adjacency matrix. Since the graph $\graph$ is
$k$-regular, exactly $k$ of the $\cir_m$'s are non-zero. For the
matrix $\Amat$ to be skew-symmetric, we also need the condition
$\cir_m = - \cir_{\NumPlayers-m+2}$ for $m \in \nodes \setminus
\{1\}$. Using the skew-symmetry, we have $\Amat^\top \Amat = -
\Amat^2$, as a result of which $\opnorm{\Amat}$ is equal to the
largest absolute eigenvalue of $\Amat$.  By standard Fourier analysis,
the eigenvalues of a circulant matrix are of the form
\begin{align*}
	\lambda_j = \sum_{m=1}^{N} r_m \omega_N^{j(m-1)}, \quad \mbox{for $j
		= 0, 1, \ldots, N - 1$,}
\end{align*}
where $\omega_N = e^{2 \pi i/N}$ is the $N^{\text{th}}$ root of
unity. The $(m-1)^{\text{th}}$ and $(N-m+1)^{\text{th}}$ powers of
$\omega_N$ are symmetric around the real line, i.e., $e^{2 \pi i
	j(m-1)\NumPlayers} = e^{ - 2 \pi i j(N-m+1)/N}$. From the
skew-symmetric property $\cir_{m} = - \cir_{\NumPlayers-m+2}$, for $m
\leq \lceil N/2 \rceil $ we have
\begin{align*}
	r_m \omega_N^{j(m-1)} + r_{N-m+2} \omega_N^{j(N-m+1)} = 2 i r_m
	\sin(2 \pi j(m-1)/\NumPlayers).
\end{align*}
It follows that the absolute value of any of the eigenvalues can be lower bounded as
\begin{align*}
	|\lambda_j| &\geq \Big | \sum_{1< m \leq \big \lceil
		\tfrac{N}{2} \big \rceil } r_m \sin(2 \pi j (m-1)/N) \Big |.
\end{align*}
This can be lower bounded further depending on whether $k/4$ is even
or odd. Let
\begin{align*}
	\lambda_{\text{lb}}^{(j)} &\defn \lfloor \tfrac{k}{4} \rfloor
	\min_{1 < m<\ell \leq \big \lceil \tfrac{N}{2} \big \rceil} | \sin(2
	\pi j (m-1)/N) - \sin(2 \pi j (\ell-1)/N) |.
\end{align*}
We then have the lower bound
\begin{align*}
	|\lambda_j| &\geq \begin{cases} \lambda_{\text{lb}}^{(j)} , &
		\tfrac{k}{4} \text{ even}, \\ \lambda_{\text{lb}}^{(j)} + \min_{1
			< m \leq \big \lceil \tfrac{N}{2} \big \rceil } | \sin(2 \pi j
		(m-1)/N) |, & \tfrac{k}{4} \text{ odd},
	\end{cases}
\end{align*}
thereby showing that $\opnorm{\Amat} = \Omega(k)$. As we already have
the bound $\opnorm{\Amat} \leq \dmax$, it follows that
$\opnorm{\BRMatrix}$ grows linearly with $\dmax$ for this graph.

\section{Consequences for two-player zero-sum games}
\label{sec:results_zs_matrix}

Our results for general polymatrix games imply more specific
guarantees for the canonical class of two-player zero-sum matrix
games, which we describe here. We begin in~\Cref{SecFull2P} with
guarantees for the dynamics in~\Cref{AlgFull} that apply to the full
information setting. In~\Cref{SecMin2P}, we turn to the two timescale
updates in~\Cref{AlgPartial} that apply to the \infosetuplower
setting.

A two-player zero-game is fully defined by the payoff matrix $\PayOff$
for player $1$, and payoff matrix $-\PayOff^\top$ for player $2$.
Consequently, we can obtain upper bounds on the iteration complexity
by setting $\NumPlayers = 2$ in the upper bounds
from~\Cref{ThmFullInfoPoly,ThmMinInfoPoly} and bounding the spectral
norm $\opnorm{\BRMatrix}^2$ appropriately.  In this very special case,
we have $\BRMatrix = \left [ \begin{matrix} 0 & \PayOff
	\\ -\PayOff^\top & 0
\end{matrix}\right ] $,
from which it follows that $\opnorm{\BRMatrix} = \opnorm{\PayOff}$.


\subsection{Two-player zero-sum games with full information}
\label{SecFull2P}

We begin by stating finite-sample guarantees on the full information
procedure summarized in~\Cref{AlgFull}.  The shorthand $\Vinit \defn
\stratLyap(\BStrat_1)$ for the initial value of the Lyapunov function
$\stratLyap$ from equation~\eqref{EqnNovelLyapunov} used here
corresponds to the case with $\NumPlayers=2$.

\begin{corollary}[Nash gap finite-sample guarantees]
	\label{ThmFullInfo}
	Consider the full information learning dynamics (\Cref{AlgFull})
	initialized with \mbox{$\tau = \newlcc{} \epsilon/\Amax$}.  Then the
	iteration complexity $K(\epsilon)$ is bounded as follows:
	\begin{enumerate}[label=(\alph*)]
		\item For the constant stepsizes $\beta_k \equiv \beta \defn
		\epsilon^2/(\newlcc{} \opnorm{\PayOff}^2 \log \Amax)$,
		\begin{align*}
			K(\epsilon) \leq \frac{\newlcc{} \opnorm{\PayOff}^2 \log
				\Amax}{\epsilon^2}\log \big( \frac{\StratLyapVal{1}}{\epsilon}
			\big).
		\end{align*}
		\item For the inverse linear stepsize $\beta_k = \beta/k$ for some
		$\beta \in (1,2]$,
		\begin{equation*}
			K(\epsilon) \leq \frac{\newlcc{}
				\opnorm{\PayOff}^2 \beta^2 \stratLyapVal{1}
				\log \Amax }{(\beta-1)\epsilon^2}.
		\end{equation*}
		\item For the inverse polynomial stepsize $\beta_k =
		\beta/(k+\koff)^\myexp$ where $\beta \in (0,1)$, $\myexp \in (0,1)$
		and $\koff \geq \big ( \frac{1-\myexp}{\beta}\big
		)^{\frac{1}{1-\myexp}}$,
		\begin{align*}
			K(\epsilon) \leq \newlcc{}( k_0^{1-\myexp} + \frac{(1-\myexp)}{\beta}
			\log \frac{\stratLyapVal{1}}{\epsilon} \big )^{1/(1-\myexp)} + \big (\frac{ \newlcc{} \opnorm{\PayOff}^2 \beta
				\log\Amax}{\epsilon^2} \big )^{1/\myexp} .
		\end{align*}
	\end{enumerate}
\end{corollary}
The bounds on the iteration complexity in~\Cref{ThmFullInfo} scales as
$\mathcal{O}(\opnorm{\PayOff}^2/\epsilon^2)$.


\subsection{Two-player zero-sum games with minimal information}
\label{SecMin2P}

We now turn to bounds applicable to~\Cref{AlgPartial} that applies to
the minimal information setting.  In this case, we give explicit
bounds on the \textit{iteration complexity} $K(\epsilon)$, meaning the
minimum number of rounds required to ensure that $\E
\nashg(\BStrat_{K(\epsilon) + 1}) \leq \epsilon$.  We show the
existence of a polynomial iteration complexity with a scaling of the
order $(1/\epsilon)^{8 + \offset}$, where the rate parameter $\offset
> 0$ can be chosen arbitrarily close to zero.  As was the case with
zero-sum polymatrix games, the price of taking $\offset \rightarrow
0^+$ manifests in the growth of certain pre-factors; we use the
notation $g(\offset)$ and variants thereof to indicate terms of this
type.

Our result applies to~\Cref{AlgPartial} where the temperature and the
timescale separation constant are set as
\begin{equation}
	\label{eqn:tau_c_init_twoP}
	\tau = \frac{g_{\tau}(\offset)\epsilon}{ \log \Amax} \quad \text{and}
	\quad \timescalesep = \frac{g_{\alpha,\beta}(\offset)
		\tau^3}{\opnorm{\PayOff}^2} \quad \text{, respectively.}
\end{equation}
\begin{corollary}[Nash gap finite-sample guarantees]
	\label{ThmMinInfo}
	Suppose that~\Cref{AlgPartial} is run with the
	parameters~\eqref{eqn:tau_c_init_twoP}, and with the initial mixed
	strategies $\BStrat_1$ being uniform.  Then the iteration complexity
	$K(\epsilon)$ is bounded as follows:
	\begin{enumerate}[label=(\alph*)]
		\item For the constant stepsize $\beta_k \equiv \beta \defn
		\frac{g_1(\offset)\epsilon^{8+\offset}}{\Amax^7\opnorm{\PayOff}^6}$, we
		have
		\begin{align*}
			K(\epsilon) \leq K^\star(\epsilon,\offset) \defn
			\frac{\Amax^7\opnorm{\PayOff}^6}{g_2(\offset)\epsilon^{8+\offset}}
			\log \big ( \frac{3 \TotLyapVal{1}}{\epsilon}\big ) ,
		\end{align*}
		\item For the inverse polynomial stepsize $\beta_k =
		\frac{\beta}{(k+\koff)^\myexp}$ for some exponent $\myexp \in
		(0,1)$, offset \mbox{$\koff = \big \lceil
			(\frac{2\myexp}{\beta})^{1/(1-\myexp)} \big \rceil$,} and
		$\beta=\frac{g_1(\offset)\epsilon^{8+\offset}}{ \Amax^7
			\opnorm{\PayOff}^6}$, we have
		\begin{align*}
			K(\epsilon) \leq K^\star(\epsilon,\offset) & \defn \Big \{
			\frac{(1-\myexp)\Amax^7
				\opnorm{\PayOff}^6}{g_3(\offset)\epsilon^{8+\offset}} \log
			\big ( \frac{3 \TotLyapVal{1}}{\epsilon}\big ) \Big
			\}^{\frac{1}{1 - \myexp}}.
		\end{align*}
	\end{enumerate}
\end{corollary}

\section{Connections among Lyapunov functions}
\label{sec:lyapunov_remarks}

It is worthwhile comparing the Lyapunov functions used in our analysis
with those from related work.  For tracking the strategy updates, the
proofs of both~\Cref{ThmFullInfoPoly,ThmMinInfoPoly} make use of the
Lyapunov function
\begin{subequations}
	\begin{align}
		\label{EqnRepeat}
		\stratLyap(\BStrat) & \defn \sum_{i=1}^{\NumPlayers}
		\underset{\hat{\strat} \in \simplex{i}}{\max} \left\{
		{\hat{\strat}}^\top \AvgPayoff{i}(\strat^{-i}) + \tau
		\EntFun(\hat{\strat})\right\},
	\end{align}
	where $H$ is the Shannon entropy~\eqref{EqnDefnShannon}.  (In
	addition, the proof of~\Cref{ThmMinInfoPoly} also exploits an
	additional Lyapunov function for tracking the $q$-updates.)
	
	As noted previously, our Lyapunov function $\stratLyap$ is an
	entropically-regularized \mbox{$\NumPlayers$-player} version of the
	Lyapunov function used by Harris~\cite{harris1998} for two-player
	zero-sum matrix games. Using the notation of our paper, this
	$\NumPlayers$-player extension of Harris' function takes the form
	\begin{align}
		\label{EqnDefnHarris}
		\HarrisLyap(\BStrat) & \defn \sum_{i=1}^{\NumPlayers} \underset{\pihat
			\in \simplex{i}}{\max} \ {\hat{\strat}}^\top
		\AvgPayoff{i}(\strat^{-i}).
	\end{align}
	The work of Chen et al.~\cite{zaiwei2023} makes use of an alternative
	Lyapunov function for two-player zero-sum matrix games; its
	$\NumPlayers$-player extension is given by
	\begin{align}
		\label{EqnDefnChen}
		\valt(\BStrat) & \defn \sum_{i=1}^\NumPlayers \Big[ \underset{\pihat
			\in \simplex{i}}{\max} \big\{ {\hat{\strat}}^\top
		\AvgPayoff{i}(\strat^{-i}) + \tau \EntFun(\pihat)\big\} - \tau
		\EntFun(\pii) \Big].
	\end{align}
\end{subequations}
Observe that this function differs from our Lyapunov
function~\eqref{EqnRepeat} by the subtraction of the additional
entropy terms (i.e., the term $-\tau \EntFun(\pii)$).  As we discuss
below, for zero-sum games, this Lyapunov function~\eqref{EqnDefnChen}
has a natural interpretation in terms of the Kullback--Leibler
divergence. \\

\noindent \myparagraph{Connections to $\HarrisLyap$} Let us first
compare and contrast $\stratLyap$ with the Harris Lyapunov function
$\HarrisLyap$ from equation~\eqref{EqnDefnHarris}.  The Lyapunov
function \HarrisLyap \, is directly related to Nash equilibria, since
we have $\HarrisLyap(\BStrat) = 0$ whenever $\BStrat$ is a Nash
equilibrium. In contrast, since the Shannon entropy is non-negative,
our Lyapunov function~\eqref{EqnRepeat} need not be zero at a Nash
equilibrium nor at a $\tau$-regularized Nash equilibrium.  Although we
cannot hope to drive down $\stratLyap(\BStrat)$ to zero as $\BStrat$
approaches a Nash equilibrium (as in a standard Lyapunov analysis),
our function $\stratLyap$ does have three key properties that enable
our analysis:
\begin{enumerate}[label=(\alph*)]
	\begin{subequations}
		\item First, from equation~\eqref{eqn:nash_bound}, we have the
		bound $\nashg(\BStrat) \leq \stratLyap(\BStrat)$, so that
		any strategy $\BStrat$ with $\stratLyap(\BStrat) \leq
		\epsilon$ has Nash gap at most $\epsilon$.
		\item Second, we have the approximation error bound
		\begin{align}
			\big | \stratLyap(\BStrat) - \HarrisLyap(\BStrat) \big |
			\leq \NumPlayers \tau \log \Amax,
		\end{align}
		so that $\stratLyap$ is a good approximation to $\HarrisLyap$ for
		small $\tau$.
		\item Third, the Lyapunov function $\stratLyap$ is differentiable and
		$L \defn \frac{\opnorm{\BRMatrix}^2}{\tau}$-smooth with respect to
		the Euclidean norm on $\BigSimplex$, meaning that
		\begin{align}
			\| \nabla \stratLyap(\BStrat) - \nabla \stratLyap(\tilde{\BStrat})
			\|_2 & \leq L \| \BStrat - \tilde{\BStrat}\|_2 \qquad \mbox{for all
				$\BStrat, \tilde{\BStrat} \in \BigSimplex$.}
		\end{align}
	\end{subequations}
\end{enumerate}

\myparagraph{Connections to $\Valt$} The $(1/\tau)$-scaling of the
smoothness constant in item (c) is crucial to our analysis, and
underlies our choice of Lypunov function $\stratLyap$, as opposed to
the functions $\HarrisLyap$ and $\Valt$ used in past work.  Our
analysis based on $\stratLyap$ enables us to prove a finite-sample
guarantee with polynomial scaling in $(1/\epsilon)$. In contrast, for
the Lyapunov function~\eqref{EqnDefnChen} used by Chen et
al.~\cite{zaiwei2023}, the smoothness constant increases exponentially
in the inverse temperature $(1/\tau)$.  Since obtaining an
$\epsilon$-Nash-gap requires reducing the temperature, this scaling
means that the finite-sample guarantees in the paper~\cite{zaiwei2023}
also scale exponentially in $(1/\epsilon)$.

As a side-remark, we note that the Lyapunov
function~\eqref{EqnDefnChen}---despite its less desirable smoothness
properties---has a very natural interpretation in terms of the
Kullback--Leibler (KL) divergence $\operatorname{KL}(p , q)$ between
two discrete distributions $p$ and $q$.  With this notation, consider
the Lyapunov function
\begin{align}
	\label{eqn:KL_tau_br}
	\vkl(\BStrat) & \defn \sum_{i=1}^{\NumPlayers}
	\operatorname{KL}(\strat[i],
	\sigma_{\tau}(\AvgPayoff{i}(\strat^{-i}))).
\end{align}
By standard properties of the KL divergence, it can be seen that this
function is zero if and only if $\strat[i] =
\sigma_{\tau}(\AvgPayoff{i}(\strat^{-i})))$ for $i \in \nodes$, so its
minima correspond to the set of $\tau$-regularized Nash
equilibria~\eqref{EqnTauRegNash}.  In fact, for any zero-sum game, the
KL-based function~\eqref{eqn:KL_tau_br} is proportional to $\Valt$,
and we can use it to show that $\tau$-regularized Nash equilibria are
unique.  We summarize these facts in the following:

\begin{prop}
	\label{PropLyapKL}
	For any zero-sum polymatrix game and any $\tau > 0$, we have the
	equivalence
	\begin{align}
		\label{EqnEquivalence}
		\valt(\BStrat) \equiv \tau \; \vkl(\BStrat).
	\end{align}
	Moreover, there is a unique $\tau$-regularized Nash
	equilibrium~\eqref{EqnTauRegNash}, corresponding to the unique global
	minimum of the function $\vkl$.
\end{prop}
\noindent See~\Cref{AppPropLyapKL} for the proof of this claim. \\

The proof of the equivalence~\eqref{EqnEquivalence} hinges crucially
on the zero-sum nature of the polymatrix game.  Since $\valt$ is a
strictly convex function by inspection, it implies that the KL-based
function $\vkl$ is strictly convex for any zero-sum polymatrix game.
(Again, this strict convexity need not be true in general.)  The
global minima of $\vkl$ are equivalent to $\tau$-regularized Nash
equilibria, and since $\vkl$ is strictly convex, the claimed
uniqueness condition follows.


\subsection{Proof of Proposition~\ref{PropLyapKL}}
\label{AppPropLyapKL}

We prove the claimed equivalence in equation~\eqref{EqnEquivalence} as a
consequence of the following auxiliary result: the KL-based function
$\vkl$ can be written as
\begin{align}
	\label{EqnAuxiliary}
	\vkl(\BStrat) = \sum_{i=1}^\NumPlayers \big \{ - \EntFun(\pii) +
	\Gfun_{\tau}(\AvgPayoff{i}(\strat[-i]))\big \},
\end{align}
where $\Gfun_\tau(\theta) \defn \log \big( \sum_{\action \in \A^i}
e^{\theta(\action)/\tau} \big)$ for any vector $\theta \in
\R^{|\A^i|}$.

Taking this auxiliary result~\eqref{EqnAuxiliary} as given, let us
complete the proof of the equivalence~\eqref{EqnEquivalence}.  By
standard results on exponential families (e.g.,~\cite{WaiJor08}), the
function $\Gfun_\tau$ is convex, and its Legendre dual (up to a
rescaling by $\tau$) is the Shannon entropy.  Consequently, we can
write
\begin{align*}
	\Gfun_{\tau}(\AvgPayoff{i}(\strat[-i])) = \max_{\pihat \in \simplex{i}} \big\{ \big
	\langle \pihat, \frac{1}{\tau} \AvgPayoff{i}(\strat[-i]) \big \rangle +
	\EntFun(\pihat)\big\} = \frac{1}{\tau} \underset{\pihat \in
		\simplex{i}}{\max} \big\{ \big \langle \pihat, \AvgPayoff{i}(\strat[-i]) \big \rangle +
	\tau \EntFun(\pihat) \big\}.
\end{align*}
Substituting this equivalence into the auxiliary
claim~\eqref{EqnAuxiliary} and re-arranging yields the claimed
equivalence $\valt(\BStrat) \equiv \tau \vkl(\BStrat)$.

It remains to prove the auxiliary claim~\eqref{EqnAuxiliary}.  By
definition of the function $\vkl$, we have
\begin{align*}
	\vkl(\BStrat) &= \sum_{i=1}^\NumPlayers \operatorname{KL}\big(\strat[i] ,
	\sigma_\tau\big(\AvgPayoff{i}(\strat[-i])\big)\big) \\
	& \stackrel{(i)}{=} \tau
	\sum_{i=1}^\NumPlayers\big\{-\EntFun\big(\strat[i]\big)-\big\langle\strat[i], \log
	\sigma_\tau \big(\AvgPayoff{i}(\strat[-i])\big)\big\rangle\big\} \\
	& \stackrel{(ii)}{=}
	\sum_{i=1}^\NumPlayers\big\{-\EntFun\big(\strat[i]\big)-\big\langle\strat[i],
	\frac{1}{\tau} \AvgPayoff{i}(\strat[-i])-G_\tau\big(\AvgPayoff{i}(\strat[-i])\big) \mathbf{1}
	\big \rangle\big\} \\
	& = \sum_{i=1}^\NumPlayers \big \{-\EntFun\big(\strat[i]\big) +
	\Gfun_\tau\big(\AvgPayoff{i}(\strat[-i])\big) \big\} -
	\frac{1}{\tau} \sum_{i=1}^\NumPlayers \big \langle \strat[i],
	\AvgPayoff{i}(\strat[-i]) \big \rangle \\
	& \stackrel{(iii)}{=} \sum_{i=1}^\NumPlayers \big \{ - \EntFun
	\big(\strat[i]\big)+G_\tau\big(\AvgPayoff{i}(\strat[-i]) \big) \big\},
\end{align*}
where step (i) follows from the definition of the KL divergence; step
(ii) follows from the definition of $\Gfun_{\tau}$; and step (iii)
follows from the zero-sum property. \\

Finally, we establish uniqueness of the $\tau$-regularized Nash
equilibrium~\eqref{EqnTauRegNash}.  Note that any $\tau$-regularized
NE is a global minimizer of $\vkl$ over the set $\bigdelta$.  Since
$\vkl$ is continuous and the set $\bigdelta$ is compact, the global
minimum is achieved.  Since $\valt$ is a strictly convex function, the
equivalence~\eqref{EqnEquivalence} ensures that $\vkl$ is also
strictly convex.  Consequently, the function $\vkl$ has a unique
global minimum, meaning that there is a unique $\tau$-regularized NE.



\end{document}